\documentclass[mnsc,nonblindrev]{informs3} 

\OneAndAHalfSpacedXI

\usepackage{natbib}
 \bibpunct[, ]{(}{)}{,}{a}{}{,}%
 \def\bibfont{\small}%
\usepackage{subcaption}
\usepackage[hypertexnames=false]{hyperref} 
\usepackage{multirow}
\hypersetup{
    colorlinks=true,
    linkcolor=red,
    filecolor=magenta,      
    urlcolor=cyan,
    citecolor = blue,
}
\usepackage[capposition=top]{floatrow}
\usepackage{algorithm}
\usepackage{algorithmic}

\TheoremsNumberedThrough     
\ECRepeatTheorems

\EquationsNumberedThrough    

\MANUSCRIPTNO{MS-0001-1922.65}

\def\jp#1{{\color{black}#1}}

\begin{document}



\RUNTITLE{Transfer Policies for Parallel Queues}

\TITLE{Dynamic Transfer Policies for Parallel Queues}

\ARTICLEAUTHORS{%
\AUTHOR{Timothy C. Y. Chan, Jangwon Park, Vahid Sarhangian}
\AFF{Department of Mechanical and Industrial Engineering, University of Toronto, ON, CANADA} 

} 
\ABSTRACT{%
We consider the problem of load balancing in parallel queues by transferring customers between them at discrete points in time. Holding costs accrue as customers wait in the queue, while transfer decisions incur both fixed (setup) costs and variable costs \jp{that increase with the number of transfers and travel distance, and vary by transfer direction.} Our work is primarily motivated by inter-facility patient transfers \jp{to address imbalanced congestion and inequity in access to care during surges in hospital demand}. Analyzing an associated fluid control problem, we show that under general assumptions, including time-varying arrivals and convex holding costs, the optimal policy partitions the state-space into a well-defined \emph{no-transfer region} and its complement, implying that transferring is optimal if and only if the system is sufficiently imbalanced. In the absence of fixed transfer costs, an optimal policy moves the state to the no-transfer region’s boundary; in contrast, with fixed costs, the state is moved to its relative interior. \jp{Leveraging our structural results, we propose a simulation-based approximate dynamic programming (ADP) algorithm to find effective transfer policies for the stochastic system. We investigate the performance and robustness of the \jp{fluid and ADP policies} in a case study calibrated using data during the COVID-19 pandemic in the Greater Toronto Area, which demonstrates that transferring patients between hospitals could result in up to 27.7\% reduction in total cost with relatively few transfers.}}



\maketitle

%
\section{Introduction}\label{sec:intro}
The problem of load balancing in parallel queues has applications in various areas including computing and networking, service, and healthcare operations. Most studies in the literature focus on the routing decisions, i.e., which of the queues a newly arriving customer should be routed to upon arrival. See for instance \cite{vanderboorScalableLoadBalancing2022} for a survey focusing on applications in communication networks, and \cite{chenSurveySkillbasedRouting2020} focusing on applications in service and healthcare operations. In this paper, we are concerned with settings where load balancing is conducted through \emph{transfers} between queues, i.e., after customers have joined a queue. 

The primary motivation for our study is the use of inter-facility patient transfers between hospitals to \jp{address surges in hospitalization demand during events such as pandemics, mass casualty events, and natural disasters. For example, during the COVID-19 pandemic, patient transfers were frequently used as a means to address the geographical mismatch between demand for hospitalization and the available hospital capacity, particularly in Intensive Care Units (ICUs), as seen in Canada \citep{chanOptimizingInterHospitalPatient2023a}, the U.S. \citep{henryInterfacilityPatientTransfers2024}, Australia \citep{ciniInterhospitalTransferClinical2023}, and the Netherlands \citep{dijkstraDynamicFairBalancing2023}. In contrast to load balancing through ambulance diversion (e.g., \citealt{dolanHospitalLoadBalancing2022}), inter-facility transfers were conducted \textit{after} the arrival and accommodation of patients in the hospital or the Emergency Department (ED) and often covered larger distances. For example, while ambulance diversion typically increases transport distance modestly, by 1.7--7 minutes \citep{ongAmbulanceDiversionIts2025}, transfers during the pandemic took nearly two hours on average \citep{tienCriticalCareTransport2020}.} In addition, transfer decisions were made and implemented at a much slower time-scale (e.g., weekly or daily) compared to arrivals of new patients, as they required significant coordination and information sharing between hospitals. 
 
The problem of inter-facility patient transfers poses several new operational features and constraints that have not been considered previously in the context of load balancing in parallel queues. \jp{First, the decision maker can directly control the number of customers in different queues (i.e., patients waiting in the ED or inpatient wards to be admitted to ICU) through transfers (while preserving the total number) but incurs \textit{transfer costs} in doing so. This includes a variable component that is proportional to the number of transfers, varies with the direction of transfers, and scales with the distances between queues.} In addition, it includes a fixed component to capture the effort associated with coordination and information sharing required to do even a single transfer. Second, decisions are typically made at discrete times (e.g., every morning or once a week) and over much longer time-scales relative to that of arrivals, service completions, and the time taken to complete transfers. This is in contrast to continuous-time control where decisions are made at arrival or service completion epochs. Third, a transient control formulation is more appropriate because transfer decisions typically arise in response to a ``shock" to the system that has pushed the system to an undesirable state and in presence of non-stationary arrivals. We note that, while our work is primarily motivated by inter-facility patient transfers, these features are also present in other service and telecommunication systems. In cloud computing, for instance, a central load balancer seeks to distribute user requests optimally among data centers in the presence of time-varying demand. \cite{luoSpatioTemporalLoadBalancing2015} considers control at discrete time intervals by incurring an energy cost proportional to the amount of control; while the authors optimize routing decisions from workload queues to data centers, an alternative formulation can involve direct transfers of user requests among workload queues. See also \cite{kumarIssuesChallengesLoad2019} for a survey of transfer policies in cloud computing for balancing tasks across nodes.
 
To capture these new characteristics and gain insights into the structure of optimal transfer policies, we consider a general network of parallel queues. Each queue receives dedicated arrivals according to independent non-stationary Poisson processes. Service times are exponentially distributed with queue-dependent rates. Customers incur holding costs in queues according to queue-dependent convex non-decreasing functions. At each discrete control epoch, a central decision maker can transfer customers between queues to balance holding costs, but incurs variable and fixed transfer costs in doing so. The objective is to minimize the total expected holding and transfer costs over a finite horizon. 

Optimal load balancing through transfers has been studied in the literature (e.g., \citealt{downDynamicLoadBalancing2006}) but focusing on two-queue settings under stationary dynamics and continuous-time control (see Section \ref{ssec:lit_review} for a detailed discussion). When the decision maker can control the system in continuous-time, transfers are made to queues with higher holding costs only when they are empty. But in discrete-time, determining optimal transfers requires a careful balancing of holding costs, transfer costs, and idleness. The type of control we consider (i.e., instantaneous state changes) connects our work to the literature on impulse/singular control, typically studied for one-dimensional diffusion processes and under stationary dynamics. In contrast, we consider a multi-dimensional fluid control problem with non-stationary arrivals; see Section \ref{ssec:lit_review} for additional discussion.

Our main contributions and results can be summarized as follows.
\begin{itemize}
    \item \textbf{Parallel queueing model with impulse control}: We formulate the problem of dynamic transfers as a discrete-time stochastic control problem for a general parallel queueing network with non-stationary arrivals, fixed and variable (linear) transfer costs, and convex holding costs. We propose an associated fluid control problem that allows us to characterize and gain insights into the structure of the optimal policy.
    \item \textbf{Structure of the optimal policy}: We characterize the structure of the optimal fluid policy under fairly general assumptions including time-varying arrivals and convex non-decreasing holding costs. We show that the optimal policy partitions the state-space into a single well-defined \emph{no-transfer region} and its complement, such that transferring is optimal if and only if the state of the system is sufficiently imbalanced. We further establish that when transferring is optimal and there are no fixed costs, it is optimal to move the state to the boundary of the no-transfer region. In contrast, with fixed costs, it is optimal to move the state to the relative interior of the region. Practically, this implies that the optimal policy tends to transfer larger numbers of customers at a time and less frequently in the presence of setup costs. When specialized to a two-queue system, this structure reduces to a state-dependent $(s,S)$ policy commonly arising in inventory control.
    \item \jp{\textbf{Approximate dynamic programming (ADP) algorithm}: We leverage the structural results to develop a simulation-based approximate policy iteration (API) algorithm for the original stochastic control problem. The algorithm directly approximates the no-transfer region via a classifier that labels each state as inside or outside the region starting with the fluid solution, and iteratively refines the classification while provably preserving the structural properties of the region. The algorithm further utilizes Common Random Numbers (CRN) to reduce variance and bypass computing future costs through coupling. It is also applicable to other problems where a region-of-inaction policy is optimal. }
    \item \textbf{Numerical results and case study}: \jp{To further motivate the API algorithm, we numerically confirm that the optimal policy for the stochastic problem has the same structure as that established for the fluid control problem. Using simulation experiments, we examine the performance of the API policy and show that it consistently outperforms other benchmark policies, including the fluid policy, especially in more critically loaded settings.} Lastly, we conduct a case study calibrated using real data from a network of four intensive care units (ICU) in the Greater Toronto Area during the COVID-19 pandemic. We also relax some modeling assumptions such as exponential service times and known arrival rates, and add additional application-relevant constraints, e.g., an upper bound on the number of permissible transfers. \jp{We demonstrate that the API policy can improve the total expected system cost by up to 27.7\% over a one-week horizon, reducing the number of patient days over ICU capacity by 46. This is achieved by transferring an average of 2.6 patients per day within the network.}
\end{itemize}

\textbf{Organization of the paper}. In Section \ref{ssec:lit_review} we provide a brief review of the related literature. We describe the stochastic control problem and its associated fluid control problem in Section \ref{sec:problem_definition}. We present our main results on the structure of the optimal fluid policy in Section \ref{sec:characterization} \jp{and the ADP algorithm in Section \ref{sec:ADP}}. Section \ref{sec:numerical_experiments} summarizes our numerical experiments and the case study. We conclude the paper in Section \ref{sec:conclusion}. All proofs are provided in the Online Appendix. 

\textbf{Notation}. We denote the non-negative real line using $\mathbb{R}_+$ and the $N$-dimensional non-negative Euclidean space by $\mathbb{R}_+^N$. We use $1\{\cdot\}$ to denote an indicator function. Given two matrices $U, V \in \mathbb{R}_+^{N \times N}$, we define $U \cdot V \equiv \sum_{i\in \mathcal{N}} \sum_{j \in \mathcal{N}} U_{ij}V_{ij}$ where $\mathcal{N} \equiv \{1, \ldots, N\}$. We let $(x)^+ \equiv \max(0,x)$. If $x\in \mathbb{R}^N$, $(x)^+$ is a vector where the $i$th component is equal to $(x_i)^+$. We use $x^\top$ to indicate the transpose of $x$ \jp{and $\| x \|$ the Euclidean norm of $x$}. The vector of all ones is denoted by $\mathrm{e}$, whose dimension should be clear from the context. 

\section{Related Literature} \label{ssec:lit_review}

\textbf{Load balancing in parallel queues.} There is a large literature on dynamic load balancing for telecommunications and distributed computing networks. In that context, a single load balancer or a dispatcher directs arrivals dynamically to one of many parallel servers at the point of entry to the system. Ideally, jobs are routed to the shortest queue, but sampling all queues can be expensive. As such, a large body of literature focuses on large-server regimes and the power of sampling only two queues; see, e.g., \cite{sitaraman2001power} for a survey. Routing decisions have also been studied in service  and healthcare operations, sometimes jointly with scheduling decisions. Examples include routing calls to different server pools in contact centers, e.g., \cite{armony2005dynamic,armony2010fair}, and joint routing and scheduling of patients to hospital wards \citep{chen2023optimal}.

Closer to our work are studies that allow load balancing \textit{after} arrival of customers. \cite{heTwoQueuesTransfers2002} study a two-queue system under a threshold policy whereby if the difference in queue-length between the two queues exceeds the threshold, a fixed number of customers is transferred. Customers incur holding costs as they wait in the queues and transfers incur a variable cost.  \cite{downDynamicLoadBalancing2006} study the stability of a general parallel queueing network with transfer of customers at general, possibly random points in time. They characterize certain properties of the optimal policy for a two-queue system under general arrival and service processes, and partially characterize the structure of the optimal policy for a two-queue Markovian system under continuous-time control. For systems with more than two queues, they propose a heuristic policy. \cite{caudillo-fuentesSimpleHeuristicLoad2010a} extend their analysis and propose heuristic policies for a two-queue system with general, heavy-tailed service distributions. Our work significantly expands the structural results on the optimal policy for a much more general system with multiple queues, time-varying arrivals, and convex holding costs. We establish the structural results for the fluid control problem, but provide numerical evidence that the same structure holds for the stochastic problem as well.

\textbf{Impulse control}. With respect to the type of control, our work relates to the large body of literature on impulse control. Impulse control finds applications in diverse settings such as inventory control \citep{bensoussanOptimalityPolicyCompound2005, ormeci2008impulse, benkheroufOptimalityPolicyCompound2009, daiBrownianInventoryModels2013, daiBrownianInventoryModels2013a}, finance and economics \citep{kornApplicationsImpulseControl1999, cadenillasClassicalImpulseStochastic2000, mitchell2014impulse}, and internet congestion control \citep{avrachenkovInfiniteHorizonOptimal2015}. However, this body of literature predominantly focuses on single-dimensional control. In contrast, transfer as a control mechanism is inherently multi-dimensional because of the coupling constraint that it must preserve the total number of customers in the system. For some applications, the absence of this constraint allows one to consider a single-dimensional problem without loss of generality. Furthermore, our work differs from much of the literature by considering a transient (finite horizon) problem with non-stationary dynamics and restricting control to the beginning of discrete time intervals.

Our work contributes to the literature on multi-dimensional impulse control by establishing the structure of the optimal policy in the presence of fixed costs and in the presence of queuing dynamics. Examples of multi-dimensional impulse control problems are found in ride-hailing platforms, where the objective is to minimize the expected lost sales (or maximize profit) by repositioning the inventory such as cars or bikes among geographic locations. \cite{heRobustRepositioningVehicle2020} consider relocation decisions at discrete epochs using a distributionally robust optimization approach in which the decisions are approximated as linear functions of uncertain customer demands. For a two-location problem, the authors characterize the optimal policy as a threshold-type policy. \cite{benjaafarDynamicInventoryRepositioning2022} extend the results to a general $N$-location problem by considering a stochastic DP formulation. They characterize the optimal policy as a region-of-inaction type policy with the optimal policy moving the state to the boundary when it lies outside of the region. While the structure of the optimal policy in our problem shares similarities with these works, neither of these works consider queueing dynamics or fixed costs. Furthermore, both works consider a closed network, for which the region-of-inaction only needs to be estimated for a fixed value of the total number of vehicles.  \cite{ataDynamicDispatchCentralized2020} consider the joint decision of dispatching cars to customers and centrally relocating cars between geographic areas by considering a closed stochastic processing network \citep{harrison2003broader} and investigating an associated Brownian control problem. They consider continuous-time control and preclude fixed costs. As we show in this work, considering fixed costs fundamentally changes the structure of the optimal policy. Specifically, it moves the state to the relative interior of the no-transfer region (or region-of-inaction), rather than the boundary. \jp{Lastly, we note that fixed costs have been considered in other multi-dimensional control settings, such as make-to-order systems \citep[e.g.,][]{sunDynamicControlMaketoOrder2025}; our problem is distinct in that control decisions always affect multiple queues, which leads to a different cost and optimal policy structure.}



\textbf{Transient queueing control}: Transient queueing control problems are often challenging due to the complexity of characterizing transient dynamics, even for simple queueing models. As such, fluid and diffusion approximations are often used to derive asymptotically optimal policies as well as insights into the structure of the optimal policy. Our approach relies on a fluid approximation of the queueing dynamics arising from the \jp{\emph{conventional} scaling}. Fluid approximations (both under conventional and many-server regimes) have been leveraged in the literature to study complex scheduling and routing control problems; see for example, \cite{meyn1997stability}, \cite{maglaras2000discrete}, and \cite{bauerle2000asymptotic} for fluid-based policies for control of general queueing networks, and \cite{zychlinski2023applications} for a recent review. Most studies focus on continuous-time control and leverage optimal control theory (see, e.g., \citealt{sethi2000economic}) to characterize the structure of the optimal policy. For example, \cite{huOptimalSchedulingProactive2022} study proactive scheduling in the presence of customer deterioration and improvement. \cite{chen2023optimal} study routing and scheduling in parallel queues with time-varying arrivals. \cite{zychlinski2023managing} examine scheduling policies when customers may need multiple servers using a discrete-time model with Bernoulli arrivals and Geometric service times. We also consider a discrete-time control problem, but account for continuous-time queueing dynamics between decision epochs. \cite{chanDynamicServerAssignment2021} also consider a discrete-time control but focus on server assignment.  Our control problem differs from routing and scheduling problems both in terms of the type of control and the cost components. In particular, compared to routing/scheduling problems which focus on minimizing holding costs, transfer policies must also balance the benefits holding cost reduction with variable transfer and fixed (setup) costs. As such, our characterization of the optimal policy relies on showing (multi-dimensional) $K$-convexity \citep{gallegoKConvexityRn2005a} of the value function of a discrete-time dynamic programming (DP) formulation of the fluid control problem. \cite{chanOptimizingInterHospitalPatient2023a} develop a numerical approach for guiding patient transfers in a network of hospitals modeled as two-stage tandem queues. In contrast, here we focus on characterizing the structure of the optimal transfer policy. 

\jp{\textbf{ADP for load balancing and inventory repositioning}:}
\jp{Relevant to our work is the literature on ADP applications to routing and inventory repositioning in service and healthcare operations. Examples include ambulance redeployment \citep{maxwellTuningApproximateDynamic2013}, patient overflow management \citep{daiInpatientOverflowApproximate2019a}, and vehicle repositioning in on-demand rental networks \citep{benjaafarDynamicInventoryRepositioning2022}. A common approach in this literature has been to approximate the value function. One line of research focuses on using a linear combination of basis functions, which can be informed by the limiting fluid models or other approximate models \citep[e.g.,][]{moallemiApproximateDataDrivenDynamic, chenApproximateDynamicProgramming2009, daiInpatientOverflowApproximate2019a}. Other approaches include cutting-plane methods that leverage properties such as convexity; see \cite{benjaafarDynamicInventoryRepositioning2022} and references therein. As we show in Section \ref{sec:structure}, however, the value function is generally non-convex in the presence of fixed costs. In contrast, our approach is based on approximating the policy directly. This connects our work to policy gradient methods from the reinforcement learning literature, which search for an optimal stochastic policy typically represented using a neural network; see, e.g., \cite{dai2022queueing,sun2024inpatient} for applications to queueing problems. In contrast, we leverage the structural properties of the optimal fluid policy to search for a connected region-of-inaction policy. This enables us to solve practical instances of the problem.} 


\section{Problem Formulation}\label{sec:problem_definition}
Consider $N$ parallel single-server, First-Come, First-Served (FCFS) queues indexed by $i \in \mathcal{N}\equiv \{1, \ldots, N\}$. Customers arrive to queue $i$ according to a non-stationary Poisson process with rate $\lambda_i(t)$ and \jp{have exponentially distributed service times with rate $\mu_i$.} Decisions are made over a finite horizon of length $T$ \jp{divided into $M$ periods (discrete epochs), indexed by $m\in \mathcal{M} \equiv \{0, \ldots, M-1\}$, with each period having a fixed length $\tau$.} At the beginning of each decision epoch, the decision maker can transfer customers between queues. 


\jp{Let $X^\pi(t)=(X^\pi_1(t),\ldots,X^\pi_N(t))$ denote the process tracking the number of customers in each queue under a (transfer) policy $\pi$, and let $U^{\pi}(t_m)$ denote the transfer decision matrix at time $t_m$, where $U^{\pi}_{ij}(t_m)$ represents the number of customers transferred from queue $i$ to $j$. A policy $\pi$ is \textit{admissible} if it is non-anticipating, $U^{\pi}_{ij}(t_m) \geq 0$, and $\sum_{j \in \mathcal{N}} U_{ij}^\pi(t_m) \leq X_i(t_m^{-})$ for all $i\in \mathcal{N}$,  $m\in \mathcal{M}$. }
 
For each $i\in \mathcal{N}$, let $\{A_i(t);t\geq 0\}$ denote a unit-rate independent Poisson process corresponding to arrivals, \jp{and let $\{D_i(t);t\geq 0\}$ denote the same for service completions.} The sample paths of $X^\pi$ satisfy the following for all $m$ and $t \in [t_m,t_{m+1})$: 
\begin{equation}\label{eq:sp}
    X^\pi_i(t) = X^\pi_i(t_m^{-}) + \sum_{j \in \mathcal{N} }( U^\pi_{ji}(t_m) - U^\pi_{ij}(t_m) ) + A_i\left(\int_{t_m}^{t} \lambda_i(s)ds \right) -D_i\left(\int_{t_m}^{t} \mu_i 1\{X^\pi_i(s)>0\}ds\right),
\end{equation}
where $X_i^{\pi}(0^-)\equiv X_i^{\pi}(0)$ and $1\{X^\pi_i(s)>0\}=1$ if $X^\pi_i(s)>0$ and 0 otherwise. The terms on the right-hand-side of \eqref{eq:sp} correspond respectively to \jp{the pre-transfer queue length at queue $i$,} the \emph{net} number of customers transferred into queue $i$ \jp{(possibly negative)}, the number of new arrivals into queue $i$ up to time $t \in [t_m,t_{m+1})$, and the number of departures up to time $t \in [t_m,t_{m+1})$. 

\jp{The transfer decisions $U^{\pi}$ incur a fixed transfer cost (setup cost) of $\tilde{\kappa}(U^\pi)$. In the most general case, $\tilde{\kappa}(U) = \sum_{i \in \mathcal{N}}\sum_{j \in \mathcal{N}}\tilde{K}_{ij}1\{U_{ij} > 0\}$, which accumulates $\tilde{K}_{ij}$ for any \jp{positive} number of customers transferred from queue $i$ to $j$. There is also a variable transfer cost of $r_{ij}$ per transferred customer from queue $i$ to $j$, and a holding cost at rate $h(X(t))$, where $h(\cdot)$ is a convex, \jp{non-decreasing} function.} 
The objective is then to find an admissible policy that minimizes the total expected cost over the horizon starting at $X(0)$:
\begin{equation}
    \mathbb{E} \left[ \sum_{m \in \mathcal{M}} \int_{t_m}^{t_{m+1}} h(X^{\pi}(s)) ds  + r \cdot U^{\pi}(t_m) + \tilde{\kappa}(U^{\pi}(t_m)) \right].
\end{equation}
It is natural to think of $X(0)$ as a large and imbalanced initial state just prior to making any transfer decisions, possibly after a ``shock" to the system, and the number of periods (horizon length) to be large enough so that the effect of the shock can subside during the horizon. 

Finally, we note that from \eqref{eq:sp} it is clear that the sample path dynamics for each queue only depend on the \emph{net-transfer} $\tilde{U}_i(t_m) \equiv \sum_{j \in \mathcal{N} }( U^\pi_{ji}(t_m) - U^\pi_{ij}(t_m) ), i \in \mathcal{N}$. Hence, by picking the lowest-cost transfers $U$ that achieves a given net-transfer $\tilde{U}$ in each period, we can express the problem using the lower-dimensional control $\Tilde{U}$, or equivalently the \textit{post-transfer state} $X^\pi_i(t_m^{-}) + \tilde{U}_i(t_m)$. We leverage this observation when considering the dynamic programming formulation of the fluid control problem in the next section.

\subsection{The Fluid Control Problem} \label{sec:DP}

\jp{The fluid control problem is obtained by approximating the queueing dynamics during each period with a deterministic fluid approximation justified by a Functional Law of Large Numbers (FLLN) \citep{mandelbaum1995strong}. Specifically, consider a sequence of stochastic systems indexed by $\eta$, such that the $\eta$th system has parameters $\lambda^\eta_i(t)=\eta\lambda_i(t)$, $\mu^\eta_i = \eta\mu_i$, $\forall i \in \mathcal{N}$, and initial condition $X^\eta (0) = \eta x^0$. The scaled process $\eta^{-1}X^\eta(t)$ converges to a deterministic fluid trajectory in the limit as $\eta \rightarrow \infty$ uniformly on compact sets (u.o.c.) and with probability 1. In formulating the fluid control problem, because the number of customers and hence the size of transfers is increasing, we view the fixed cost as scaled such that $\tilde{\kappa}^{\eta}(\cdot) = \eta \tilde{\kappa}(\cdot)$ for the $\eta$th system, while the holding and variable transfer costs remain unscaled.}
\jp{
\begin{remark}In practice, we may only have access to predicted arrival rates subject to prediction errors. Assume that the $\eta$-th system has arrival rate $\eta \lambda_i(t) + \epsilon^{\eta}_i(t)$, where  $\epsilon^{\eta}_i(t)$ is the estimation error in the $\eta$-th system. Then, assuming that $\eta^{-1} \epsilon^{\eta}_i(t) \to 0$ u.o.c. with probability 1, i.e., the uncertainty of the arrival rate vanishes under fluid scaling, the fluid dynamics in \eqref{eq:fluid_dynamics} remain valid; see, e.g., \cite{chen2023optimal}. In our simulation experiments in Section \ref{sec:case_study}, we numerically evaluate the performance of our proposed policies under arrival rate estimation error. 
\end{remark}
\begin{remark}
    In the above scaling, the arrival and service rates uniformly increase while  the number of servers remains fixed. As such, the fluid dynamics in \eqref{eq:fluid_dynamics} also serve as an approximation for multiserver queues, after multiplying the service rates by the number of servers. In Section \ref{sec:case_study} we illustrate the performance of our policies for multiserver queues. 
\end{remark}
}



\jp{Let $x(t) \in \mathbb{R}_+^N$ denote the fluid state at time $t \geq 0$. We use $x[m] \equiv x(m\tau^{-})$ to denote the state of the system at the beginning of period $m\in \mathcal{M}$ \emph{before} the transfer decision is made. Further, let $u[m] \in \mathbb{R}_+^{N \times N}$ be the fluid transfer matrix in period $m$. The post-transfer fluid state then satisfies $y[m]=x[m] + \left(u[m]^\top-u[m] \right)\mathrm{e}$. Denote by $f^m:\mathbb{R}_+^N \times \mathbb{R}_+ \to \mathbb{R}_+^N$ the state transition function that returns the system state at a given time during period $m$, starting from a given (post-transfer) state. Then $f^m(y,t)$ is the solution to the following initial value problem starting from $y$:}
\jp{\begin{align}
    \frac{d^{+}}{dt}x_i(t) &= \lambda_i(t) - \mu_i 1\{x_i(t) > 0\}, \quad i \in \mathcal{N} \text{ and } t \in [m\tau, (m+1)\tau),  \label{eq:fluid_dynamics}
\end{align}}
\jp{where $\frac{d^{+}}{dt}$ is the right-derivative (see, e.g., \citealt{meyn2008control}, Page 40.)}

The minimum transfer cost $C(y-x)$ associated with a given net-transfer $y-x$ is given by
\begin{equation}\label{eq:transfer_cost}
    \begin{split}
        C(y-x) = \min_{u} &\quad r \cdot u + \tilde{\kappa}(u) \\
        \mbox{s.t.} &\quad (u^\top - u) \mathrm{e} = y - x, \\
        &\quad u \geq 0. 
    \end{split}
\end{equation}
Let $H^m(y)$ denote the holding cost incurred in period $m$, starting from post-transfer state $y$. Then 
\begin{align} 
    H^m(y) = \int_{m\tau}^{(m+1)\tau}  h(f^m(y,s)) ds.\label{eq:holding_cost}
\end{align}
\jp{Finally, denote the fluid value function by $V^m: \mathbb{R}_+^N \xrightarrow{} \mathbb{R}_+$ for each $m \in \mathcal{M}$. Then $V^m(x)$ is the minimum cost-to-go starting from $x$ in period $m$, and the optimal fluid cost is given by $V^0(x^0)$. The fluid value function satisfies the optimality equation,
\begin{equation} \label{eq:DP}
    V^m(x) = \min_{y \in\Delta(\mathrm{e}^\top x)} \left[H^m(y) + C(y - x) + V^{m+1}(f^m(y,(m+1)\tau)) \right],
\end{equation}
with $V^M \equiv 0$, where $\Delta(n)=\{y \in \mathbb{R}_+^N: \mathrm{e}^\top y = n\}$ denotes the set of all feasible \textit{post-transfer} states, and $n$ is the total number of customers to be preserved at the time of decision.}

\jp{In general, the fluid control problem is a non-linear, non-convex problem due to the discontinuous objective function. Appendix \ref{appen:numerical_solution_approach} presents an equivalent formulation and numerical solution approach for solving this problem using a mixed-integer linear program, which can be used to compute optimal fluid policies for large problem instances.}


\jp{The solution of the fluid control problem can be directly translated to an admissible control for the stochastic problem using a \jp{rolling-horizon approach \citep{powell2007approximate}}. Specifically, denote by $u^*[m]$ an optimal transfer decision matrix corresponding to the initial condition $x^m=X(t_m^-)$ in period $m$. One can construct the transfer matrix $U(t_m) = \lfloor u^*[m] \rfloor$ for the stochastic system, where $\lfloor \cdot \rfloor$ is the floor function applied component-wise, and implement only the solution corresponding to the immediate period. The fluid control problem is then re-solved with the observed initial state at the start of the next period.}


\section{Characterization of the Optimal Fluid Policy}  \label{sec:characterization}
In this section, we characterize the structure of the optimal policy for the fluid control problem. We present these results in the general case of time-varying arrivals and convex holding costs. By considering the special case of stationary arrivals and linear holding costs, we provide further insights into the trade-off between the holding cost, transfer cost, and idleness.

To characterize the structure of the optimal fluid cost and policy, we make the following three assumptions about the system's arrival rates, the holding cost function per unit time $h_i(\cdot)$, and the variable transfer costs per customer $r_{ij}$.

\begin{assumption} \label{ass:time_varying_arrivals}
    For all $i \in \mathcal{N}$, the arrival rates $\{\lambda_i(t): t \geq 0\}$ are non-negative, piecewise monotone, and have finitely many pieces. 
\end{assumption}

\begin{assumption} \label{ass:non_linear_holding_cost}
    $h_i(\cdot)$ is convex, continuous, and non-decreasing for all $i \in \mathcal{N}$.
\end{assumption} 

\begin{assumption}\label{ass:triangular_inequality}
The unit variable transfer costs satisfy the triangle inequality, i.e.,
\begin{equation*}
    r_{ij} \leq r_{il} + r_{lj}, \quad \forall i,j,l \in \mathcal{N}. 
\end{equation*}
\end{assumption}
Assumption \ref{ass:time_varying_arrivals} allows for many widely-used time-varying arrival rate functions (e.g., piecewise-constant, piecewise-linear, sinusoidal). Assumption \ref{ass:non_linear_holding_cost} allows for convex increasing holding costs, suitable for practical settings. In healthcare, for instance, the impact of congestion on clinical outcomes can increase past a certain point in hospital occupancy \citep[e.g.,][]{kuntzStressWardEvidence2015, berryjaekerPointSpeedingNegative2017}, implying a convex, increasing cost structure. Lastly, Assumption \ref{ass:triangular_inequality} states that the transfer cost from one queue to another cannot be made smaller by going through an intermediary queue, and is common in the literature \citep[e.g.,][]{zengCostSharingCapacity2018, benjaafarDynamicInventoryRepositioning2022}.

\subsection{The Joint Setup Cost} \label{ssec:joint_setup_cost}
As we establish in the sequel, in the presence of fixed transfer costs, the value function is no longer convex. As such, we exploit the notation of $K$-convexity \citep{scarfOptimalityPoliciesDynamic1960a} and its extension to $\mathbb{R}^N$ proposed by \cite{gallegoKConvexityRn2005a}.

\begin{definition} \label{def:K_convexity}
Let $\tilde{\kappa}:\mathbb{R}^{N\times N}_+ \rightarrow \mathbb{R}_+$ be a generic setup cost function with parameter $\tilde{K} \in \mathbb{R}_+^{N\times N}$. A function $V:\mathbb{R}^N_+ \xrightarrow{} \mathbb{R}_+$ is $\tilde{K}$-convex if 
\begin{equation*}
    V(\theta x + (1-\theta)y) \leq \theta V(x) + (1-\theta)[V(y) + \tilde{\kappa}(u)],
\end{equation*}
for all $\theta\in[0,1]$ and all $x,y \in \mathbb{R}_+^N$ with $y\in \Delta(\mathrm{e}^\top x$), where $u$ is the minimum-cost transfer matrix that achieves the net-transfer $y-x$, i.e., solves \eqref{eq:transfer_cost}.
\end{definition}

We now specialize this definition to a particular kind of setup cost function, whereby a transfer between any pair of queues incurs a fixed cost of $K$ for the entire system. \jp{In Section \ref{sssec:general_setup}, we demonstrate numerically that the structural results we will show in the next section remain robust to more complex forms of setup cost functions.} For any given net-transfer $z \in \mathbb{R}^N$ and $K>0$, we define the \emph{joint setup cost} function as,
\begin{align}\label{eq:joint_setup_cost}
    \kappa(z) = K1\{z\neq0\} =
    \begin{cases}
        K, & \mbox{if } z\neq 0; \\
        0, & \mbox{otherwise}.
    \end{cases}
\end{align}
Despite its simplicity, the joint setup cost is practically relevant in applications involving a central decision maker, where there is a preference or necessity for less frequent interventions and where the initial cost of planning and preparing for transfers is significant. Additionally, it satisfies the following properties which are key for establishing the structure of the value function.
\begin{lemma}\label{lem:joint_setup_cost_properties}
    The joint setup cost function in \eqref{eq:joint_setup_cost} satisfies the following properties:
    \begin{enumerate}
        \item[(i)] (Subadditivity): For all $x,y \in \mathbb{R}^N$, we have $\kappa(x+y) \leq \kappa(x) + \kappa(y)$.
        \item[(ii)] (Homogeneous of degree 0): For all $x\in \mathbb{R}^N$ and $c\neq0$, we have $\kappa(cx) = \kappa(x).$
        In particular, $\kappa(-x)=\kappa(x)$, i.e., $\kappa(\cdot)$ is an even function.
        \item[(iii)] (Decomposition of total transfer cost): Denote by $R(y-x)$ the transfer cost in going from $x$ to $y\in\Delta(\mathrm{e}^\top x)$ without accounting for the joint setup cost, i.e.,
        \begin{equation}\label{eq:variable_cost}
            \begin{split}
                R(y-x) = \min_{u \geq 0} &\quad r \cdot u \\
                \mathrm{s.t.} &\quad (u^\top - u)\mathrm{e} = y - x,
            \end{split}
        \end{equation}
        Then we have $C(y-x) = R(y-x) + \kappa(y-x)$.
    \end{enumerate}
\end{lemma}
\jp{In particular, the third property allows us to decompose the total transfer cost into a convex component $R(y-x)$ and the setup cost. By isolating the setup cost, it allows us to invoke Definition \ref{def:K_convexity} and make use of additional properties of the value function described in Appendix \ref{appen:add_properties_V}.}

Note that we can restrict the domain of $\kappa(\cdot)$ to $\mathbb{R}^N$ because the joint setup cost function depends only on the net-transfer, as opposed to the entire transfer matrix. In the rest of the paper, we simply state that a function is $K$-convex when Definition \ref{def:K_convexity} is satisfied using the joint setup cost function $\kappa(\cdot)$ with parameter $K$.


\subsection{Structure of the Optimal Fluid Policy} \label{sec:structure}
We first establish the structural properties of the single-period holding cost function $\jp{H^m}(\cdot)$ and the value function $V^m(\cdot)$, which are key in characterizing the structure of the optimal policy.

\begin{lemma} \label{lem:properties_H}
    Under Assumptions \ref{ass:time_varying_arrivals} and \ref{ass:non_linear_holding_cost}, $\jp{H^m}(\cdot)$ is convex, continuous, and non-decreasing \jp{for all $m \in \mathcal{M}$}.
\end{lemma}

\begin{theorem} \label{thm:properties_V}
    Let $\kappa(\cdot)$ be the joint setup cost function in \eqref{eq:joint_setup_cost}. Under Assumptions \ref{ass:time_varying_arrivals} and \ref{ass:non_linear_holding_cost}, $V^m(\cdot)$ is $K$-convex, continuous, and non-decreasing for all $m\in \mathcal{M}$.
\end{theorem}
An important special case is when there are no setup costs ($K=0$). Then, the joint setup cost $\kappa(z)=0$ for all $z$ and Definition \ref{def:K_convexity} reduces to the standard definition of convexity. In this case, the value function $V^m(\cdot)$ is convex and the optimal policy can be obtained by solving a convex optimization problem.
\begin{corollary} \label{cor:properties_V_no_setup}
    Suppose $K=0$. Under Assumptions \ref{ass:time_varying_arrivals} and \ref{ass:non_linear_holding_cost}, $V^m(\cdot)$ is convex, continuous, and non-decreasing for all $m\in \mathcal{M}$.
\end{corollary}
We later highlight the impact of setup cost on the structure of the optimal policy.


The significance of the above results lies in their robustness under time-varying arrival rates and convex holding costs that satisfy Assumptions \ref{ass:time_varying_arrivals} and \ref{ass:non_linear_holding_cost}. The challenge in these cases is in obtaining the closed-form expression for the state transition function $f^m(y,t)$, which is difficult to characterize since the queue length process may be highly non-linear and may not stay at zero once (and if) it is reached. In the proof provided in Appendix \ref{appen:V_properties}, our argument uses a recursive expression for $f^m(y,t)$ within each period based on specific time points such that between two successive points, the queue length process is monotone. In the rest of this section, we will assume the joint setup cost \eqref{eq:joint_setup_cost} in our model and always assume that Assumptions \ref{ass:time_varying_arrivals}--\ref{ass:triangular_inequality} hold.


\jp{Before presenting the main result, we present an intermediary result which establishes the existence of an \emph{efficient} optimal policy that never transfers customers into and out of the same queue within the same period.
\begin{proposition}\label{prop:efficient_policy}
There exists an optimal policy such that when customers are transferred, no queues are both sending and receiving customers in the same period.
\end{proposition}
Proposition \ref{prop:efficient_policy} states that in any period, we can partition the set of queues into disjoint sets --- the senders, the receivers, and the non-participants --- and consequently reduce the search for an optimal policy to the set of policies under which each queue has a dedicated role.}

\jp{Our main characterization of the structure of the optimal policy is through partitioning of the state-space into the \emph{no-transfer region} and its complement. Let $n\geq0$ denote the total number of customers in the system at the beginning of period $m$ and before transfer decisions are made. We define the no-transfer region for all $m\in \mathcal{M}$ as follows:
\begin{align*}
    \jp{\Sigma^m(n)} = \{x \in\Delta(n): H^m(x) + V^{m+1}(f^m(x, \tau)) \leq H^m(y) + C(y - x) + V^{m+1}(f^m(y, \tau)), \forall y\in \Delta(n), y\neq x\}.
\end{align*}
The left-hand side of the inequality is the cost of \textit{staying} at a given state $x$ while the right-hand side is the cost of starting at another state $y$ plus the transfer and setup costs incurred in moving from $x$ to $y$. If the cost of staying at $x$ is less than or equal to moving to state $y$, then it is not optimal to move to $y$. A state $x$ belongs in the no-transfer region if this inequality holds for all other $y\in\Delta(n)$, and therefore, the optimal policy at $x$ is simply not to transfer any customers. Conversely, if a state does not belong in $\Sigma^m(\cdot)$, it is optimal to transfer customers at that state. Throughout this section, we refer to the post-transfer state under the optimal policy, i.e., the optimal solution to \eqref{eq:DP}, as the \emph{target state}, which may be unique, as shown next. We also denote by $\partial \Sigma^m(n)$ and $\mathrm{ri}(\Sigma^m(n))$ the \textit{boundary} and the \textit{relative interior} of $\Sigma^m(n)$, respectively. We now provide a characterization of the optimal policy.}
\begin{theorem} \label{thm:structure}
    In every period $m \in \mathcal{M}$, the no-transfer region $\jp{\Sigma^m(n)}$ is non-empty, compact, and connected for all $n \geq 0$. If $x \in \jp{\Sigma^m(n)}$, it is optimal not to move from $x$. Otherwise ($x \notin \Sigma^m(n)$):
    \begin{itemize}
        \item (No transfer and setup costs): if $\kappa(\cdot)=0$ and $r=0$, \jp{there exists a unique target state to which it is optimal to move}; 
        \item (No setup costs): if $\kappa(\cdot)=0$, it is optimal to move to a target state \jp{in $\partial \Sigma^m(n)$};
        \item (Joint setup cost): if $\kappa(\cdot)$ is the joint setup cost function \eqref{eq:joint_setup_cost}, it is optimal to move to a target state \jp{in $\mathrm{ri}(\Sigma^m(n))$}. 
    \end{itemize}
\end{theorem}

\jp{This result formally establishes the optimality of the \textit{region-of-inaction} policies. It further provides the structure of the optimal policy at increasing levels of complexity of the problem to highlight the impact of different cost components.} First, Theorem \ref{thm:structure} states that if there are no transfer costs, a constant (period-dependent) target state is optimal from \textit{any} initial condition with equal total customers. Since transferring does not cost anything, the result implies that the no-transfer region can be expressed as a singleton \jp{$\Sigma^m(n)=\{y^*\}$}, containing only the target state. However, in the presence of variable transfer costs only, the no-transfer region \jp{$\Sigma^m(n)$} expands to a compact, connected set of states at which it is (strictly) optimal not to transfer customers. This indicates that transferring becomes optimal if and only if the state of the system is sufficiently imbalanced. Moreover, target states exist on the boundary of the no-transfer region, \jp{$\partial \Sigma^m(n)$}, and generally depend on the initial condition. A boundary state implies that a small perturbation can induce the policy to switch from doing nothing to transferring, and upon transfer, return to a boundary state. Consequently, the optimal policy in this case tends to move customers frequently and in small numbers, and just enough to rectify excessive imbalance. In contrast, when the joint setup cost is included, target states are positioned in the relative interior of the no-transfer region, \jp{$\mathrm{ri}(\Sigma^m(n))$}. In particular, they cannot lie on the boundary of the no-transfer region, \jp{$\partial \Sigma^m(n)$}, which implies that the optimal policy will not switch to transferring unless the number of customers fall ``low enough'' at certain queues. Thus, the optimal policy tends to transfer less often and in larger numbers. We provide a numerical illustration of the structure in the presence and absence of setup costs in Section \ref{sec:illustrative_examples}.

Intuitively, $K$-convexity of the value function allows us to extend the structure because for any $x \notin \jp{\Sigma^m(n)}$ and its target state $y$, it implies that there is a range of values $\Theta \subset [0,1]$ such that for all $\theta \in \Theta$, the point $\theta y + (1-\theta)x$ also lies in the no-transfer region. Hence, there is a positive gap between the boundary of the no-transfer region and a target state. In Appendix \ref{appen:add_properties_V}, we also show that $K$-convexity of the value function allows us to characterize additional properties, which are key for proving compactness and connectedness of the no-transfer region in Appendix \ref{appen:proof_structure}.

In closing, we elaborate on the difficulty of extending Theorem \ref{thm:structure} to the general setup cost function $\tilde{\kappa}:\mathbb{R}^{N \times N}_+ \rightarrow \mathbb{R}_+$, defined as
\begin{align}\label{eq:general_setup_cost}
    \tilde{\kappa}(u) = \tilde{K}_01\{u \neq 0\}  + \sum_{i \in \mathcal{N}} \sum_{j \in \mathcal{N}} \tilde{K}_{ij}1\{u_{ij} > 0 \}.
\end{align}
The joint setup cost function is a special case with $\tilde{K}_0 > 0$ and $\tilde{K}_{ij}=0, \forall i,j$. \jp{Extending Theorem \ref{thm:structure} under \eqref{eq:general_setup_cost} introduces several challenges. The main difficulty is in performing the ``induction step'' in the proof, i.e., showing $\tilde{K}$-convexity of $V^{m}(\cdot)$ assuming $\tilde{K}$-convexity of $V^{m+1}(\cdot)$. This step hinges on the ability to decompose the total transfer cost that leaves a convex (variable) component and isolates the setup cost term prior to invoking Definition \ref{def:K_convexity}. Unlike the joint setup cost function, the optimal solutions to \eqref{eq:transfer_cost} and \eqref{eq:variable_cost} are generally not the same under the general setup cost function, thus violating Lemma \ref{lem:joint_setup_cost_properties} and losing the critical decomposition property as a result. Additionally, the setup cost may no longer be continuous on $\Delta(n) \setminus \Sigma^m(n)$, due to possible jump discontinuities following small changes in $x$, which further complicates the analysis. Nevertheless, we conjecture that the structure of the optimal policy continues to hold under the general setup cost structure and provide numerical evidence for this conjecture in Section \ref{sssec:general_setup}.}

\subsubsection{Special Case: The Two-Queue Model.} \label{sec:two_queue_model}
For a two-queue system, Theorem \ref{thm:structure} reduces to a state-dependent threshold policy. For any $n\geq0$, $\Delta(n)$ is simply a line segment connecting $(n,0)$ and $(0,n)$ in $\mathbb{R}^2_+$. Given its non-emptiness, compactness, and connectedness, the no-transfer region $\jp{\Sigma^m(n)}$ is the shorter line segment connecting $(s_1,n-s_1)$ and $(n-s_2,s_2)$; these two points correspond to the boundary of $\jp{\Sigma^m(n)}$. As a consequence of the target states belonging either to the boundary or the relative interior of $\jp{\Sigma^m(n)}$, there also exist $S_i$ satisfying $s_i \leq S_i$ such that when customers are transferred to queue $i$, its new state becomes $S_i$. We formalize this observation below. Again, for ease of exposition, we suppress the dependence of the parameters on the period.
\jp{\begin{proposition}\label{prop:optimal_two_queue_policy}
    Consider a two-queue system and let $n\geq0$ be the initial number of customers at a given period. In every period $m \in \mathcal{M}$, there exists an optimal policy characterized by $s_1(n),S_1(n),s_2(n),$ and $S_2(n)$ such that customers are only transferred from queue $i$ to $j$ for $x_j < s_j(n)$, and after transferring, the number of customers in queue $j$ is $S_j(n)$. Furthermore:
    \begin{itemize}
        \item (No transfer and setup costs): If $\kappa(\cdot)=0$ and $r_{12}=r_{21}=0$, then $s_1(n)=S_1(n)=s_2(n)=S_2(n)$;
        \item (No setup costs): If $\kappa(\cdot)=0$, then $s_i(n)=S_i(n)$ for $i=1,2$;
        \item (Joint setup cost): If $\kappa(\cdot)$ is the joint setup cost function \eqref{eq:joint_setup_cost}, then $s_i(n) < S_i(n)$ for $i=1,2$.
    \end{itemize}
\end{proposition}}
\jp{This result is the analogue of the classical $(s,S)$ policy in inventory control in our setting.} The proposition states that each queue has a pair of (period-dependent) parameters $(s_i(n),S_i(n)), i=1,2,$ representing the optimal ``re-order'' and ``order-up-to'' points, respectively. Therefore, customers are not transferred to a queue unless the number of customers in that queue falls below the re-order point, and when it does, it is replenished to the order-up-to point. \jp{Moreover, Proposition \ref{prop:optimal_two_queue_policy} presents a specialized structure in which the parameters $(s_i(n),S_i(n)), i=1,2$, are invariant with the initial state $x$, so long as the total number $n = x_1+x_2$ is fixed.} With three or more queues, the target state depends on the entire state vector and may be different for two initial conditions even when they have the same total number of customers. \jp{Finally, we can always find an optimal policy such that $s_1(n) + s_2(n) \leq n$ holds, as this is equivalent to the existence of an optimal policy where each queue is either receiving or sending customers, a result already established in Proposition \ref{prop:efficient_policy}. This also implies that $x_1 < s_1(n)$ and $x_2 < s_2(n)$ are never possible for a given initial condition $x$ under this policy.} We provide a numerical illustration of the structure in Section \ref{sec:illustrative_examples}.


\jp{We note that for two-queue systems, Proposition \ref{prop:optimal_two_queue_policy} provides the most general structure, provided that the cost parameters are symmetric, i.e., $\tilde{K}_{12}=\tilde{K}_{21}$. In this case, the joint setup cost is equivalent to the general setup cost in \eqref{eq:general_setup_cost}: since $u_{12}$ and $u_{21}$ are never positive at the same time under Proposition \ref{prop:efficient_policy}, \eqref{eq:general_setup_cost} can be reduced to the joint setup cost with parameter $K=\tilde{K}_0 + \tilde{K}_{12} = \tilde{K}_0 + \tilde{K}_{21}$.}

Proposition \ref{prop:optimal_two_queue_policy} is consistent with and extends the partial characterization of the optimal policy in \cite{downDynamicLoadBalancing2006}. They show that under continuous-time control, each queue has a constant optimal order-up-to point. The authors conjectured, but did not prove, that when the optimal policy does not move customers to queue $i$ at state $x_i < S_i(n)$, it should also not move customers at state $x_i+\delta < S_i(n)$ for $\delta>0$. Our results provide a complete characterization of the optimal policy under discrete-time control.

\subsection{On the Role of Idleness} \label{sec:trade_off}
\jp{When control is restricted to discrete points in time, avoiding idleness plays an important role in determining the optimal policy. To gain insights into the role of idleness and further characterize $\Sigma^m(n)$, we focus on a simpler model with stationary arrivals and linear holding costs.}

\jp{Our insights are characterized through what we call the \textit{non-idleness index}, $\tau(\mu_i - \lambda_i)^+$, which represents the number of customers required to avoid idleness at queue $i$ for one period. We show that non-idleness at certain queues serves as a sufficient condition for when not transferring is optimal. 
While it is more challenging to characterize such an index under more complex time-varying arrival rates, the same insight into the role of idleness continues to apply.}

With a slight abuse of notation, let $h=(h_1, \ldots, h_N)$ be the vector of unit holding cost per customer per unit time at each queue, where $h_N \geq \cdots \geq h_1$ without loss of generality.
\begin{proposition} \label{prop:trade_off}
    Let $\kappa(\cdot)$ be the joint setup cost function in \eqref{eq:joint_setup_cost} and let $x$ be an initial condition. 
    \begin{enumerate}
        \item[(i)] If $h_N \geq \cdots \geq h_1$, then for any $i,j$ with $i < j$, there exists an optimal policy that transfers customers from queue $i$ to $j$ only when $x_j < \tau(\mu_j - \lambda_j)^+$. 
        \item[(ii)] If $h_1 = \cdots = h_N$, then it is optimal not to transfer when $x \geq \tau(\mu - \lambda)^+$. 
    \end{enumerate}
\end{proposition}
\jp{If $x_j \geq \tau(\mu_j - \lambda_j)^+$, queue $j$ does not incur any idleness during the period. Therefore, the first part of the result states that it is optimal to transfer customers into a queue with a higher unit holding cost only to prevent idleness at it. If $\lambda_j \geq \mu_j$, an optimal policy never transfers customers into that queue.
Due to symmetry, if $h_1 = \ldots = h_N$, the second part follows directly and states that if we can guarantee non-idleness at all queues for the upcoming period, it is optimal not to transfer. Therefore, we must have $\jp{\Sigma^m(n)} \supseteq \{y \in \Delta(n): y \geq \tau(\mu - \lambda)^+\}$ \jp{for all $m$}, which implies that transferring is optimal only if there will be excessive idleness.}

\jp{In the absence of transfer costs, we can characterize $\Sigma^m(n)$ more explicitly.}
\begin{proposition} \label{prop:trade_off_no_transfer_setup_costs}
    Let $x$ be an initial condition \jp{such that $\sum_{i=1}^N x_i \geq \sum_{i=1}^N\tau(\mu_i - \lambda_i)^+$,} and suppose there are no transfer and setup costs, i.e., $\kappa(\cdot)=0$ and $r=0$.
    \begin{enumerate}
        \item[(i)] If $h_N \geq \cdots \geq h_1$, then a target state $y^m$ satisfies $y_i^m \leq \tau(\mu_i - \lambda_i)^+$ for all $i \geq 2$ and $m \in \mathcal{M}$. \label{prop:trade_off_no_transfer_setup_costs_part1}
        \item[(ii)] If $h_1 = \cdots = h_N$, then it is optimal not to transfer if and only if $x \geq \tau(\mu - \lambda)^+$. Moreover, any $y \geq \tau(\mu - \lambda)^+$ is a target state. \label{prop:trade_off_no_transfer_setup_costs_part2}
    \end{enumerate}
\end{proposition}
\jp{The first part can be viewed as an analogue of the classical $c\mu$ policy in scheduling, adapted to our setting: the optimal policy moves all customers to the ``cheapest queue,'' while leaving just enough to avoid excessive idleness elsewhere. If $h_1 = \cdots =h_N$, we further obtain the exact characterization $\jp{\Sigma^m(n)} = \{y \in \Delta(n): y \geq \tau(\mu - \lambda)^+\}$, and any non-idling policy is optimal in this case.}

\subsection{Illustrative Examples} \label{sec:illustrative_examples}
In this section, we use numerical examples to illustrate and provide additional observations on the structure of the \jp{optimal fluid policy} and the no-transfer region \jp{$\Sigma^m(\cdot)$}.

\subsubsection{Structure of the \jp{Optimal Fluid Policy}.}
\begin{figure}[h]
    \FIGURE
    {\begin{subfigure}{0.45\textwidth}
        \FIGURE
        {\includegraphics[width=\textwidth]{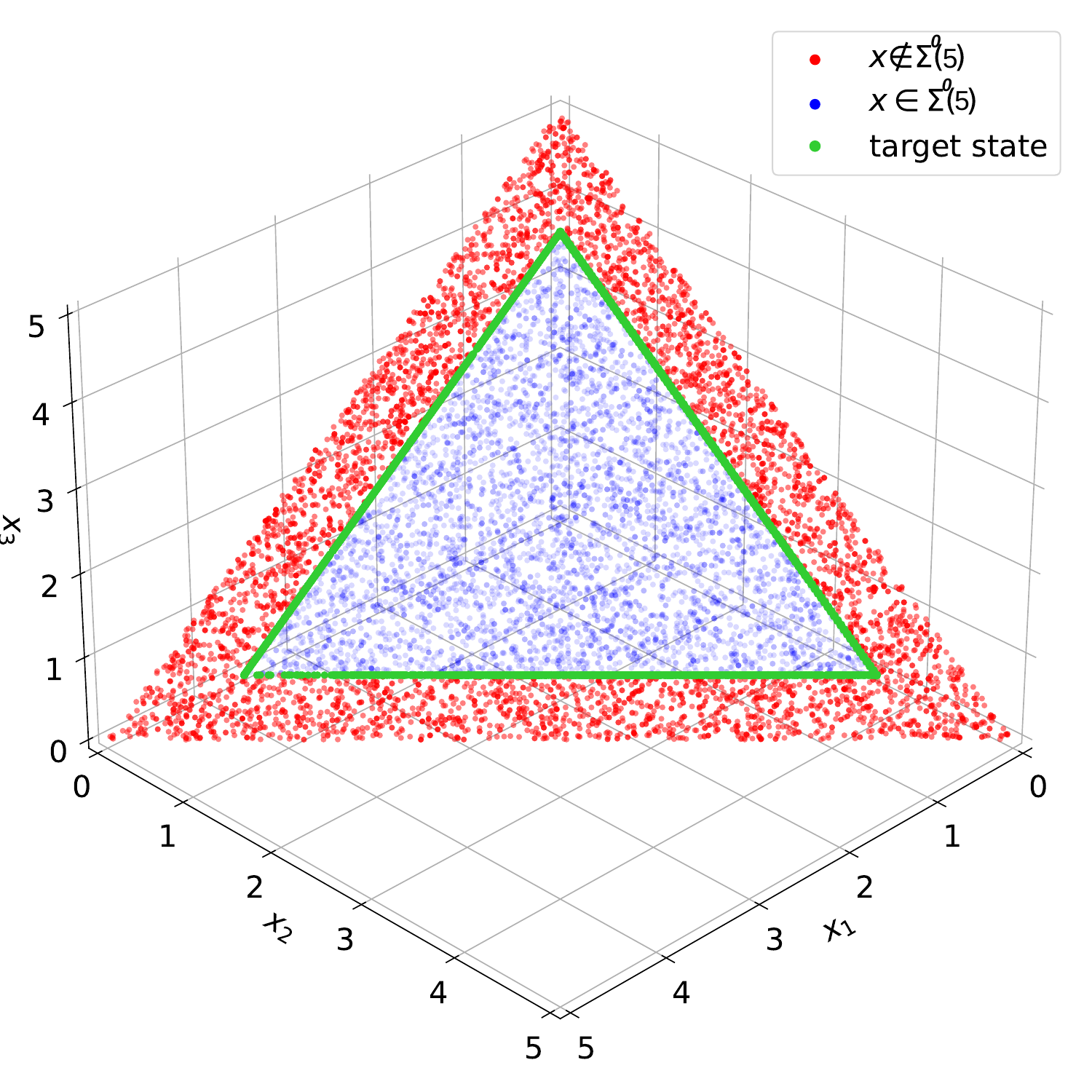}}
        {No setup cost ($K=0$) \label{fig:structure_no_setup}}{}
    \end{subfigure}
    \hspace{2em}
    \begin{subfigure}{0.45\textwidth}
        \FIGURE
        {\includegraphics[width=\textwidth]{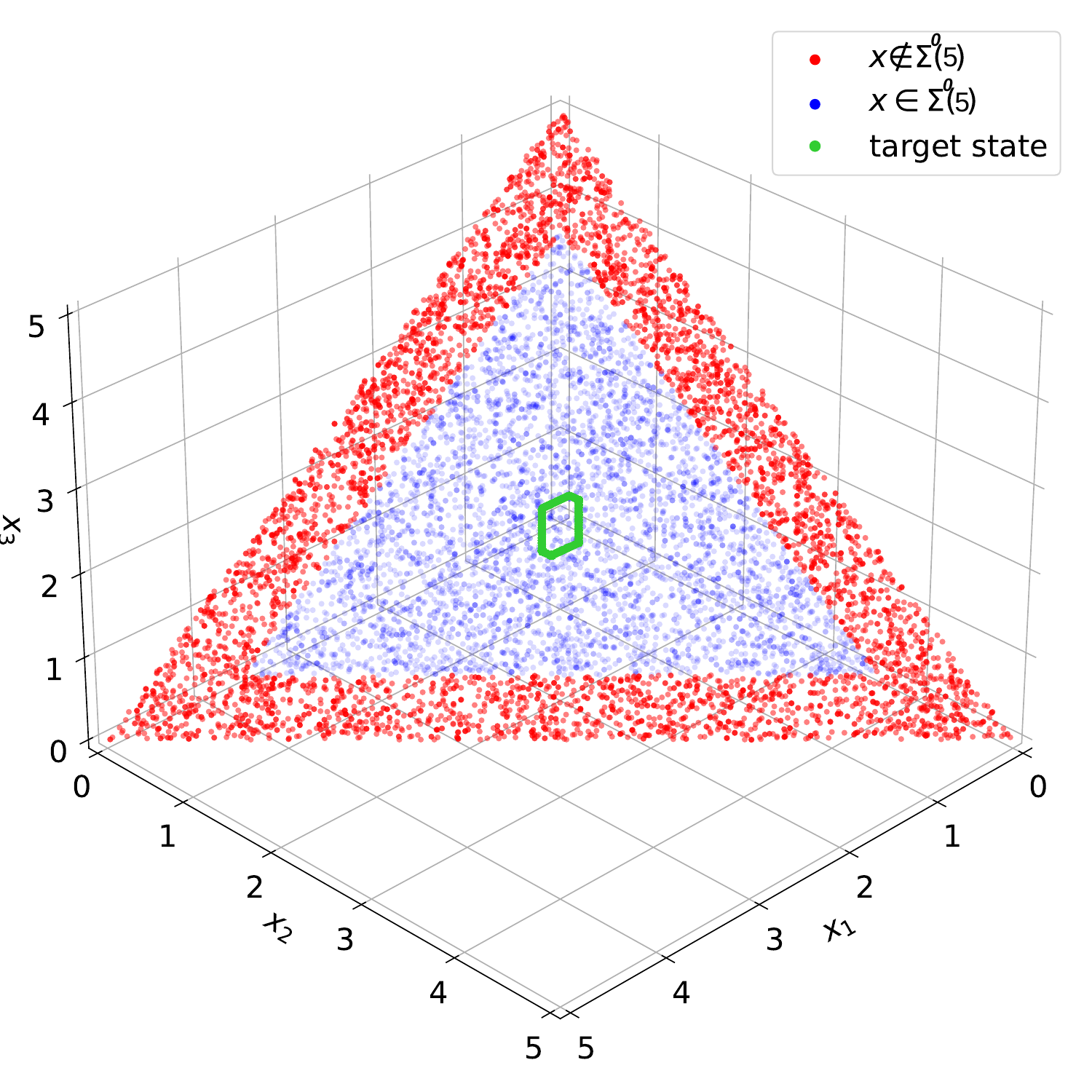}}
        {Joint setup cost ($K=5$) \label{fig:structure_joint_setup}}{}
    \end{subfigure}}
    {Structure of the \jp{three-queue optimal policy in period 0}. Red dots are states where transferring is optimal; green are target states; and blue are states where transferring is not optimal. \label{fig:structure_contrast}}
    {$\lambda_1=\lambda_2=\lambda_3=0.9, \mu_1=\mu_2=\mu_3=1, h_1=h_2=h_3=1, r_{12}=r_{21}=2, r_{13}=r_{31}=4, r_{23}=r_{32}=3, \tau=5, M=5$.}
\end{figure}

First, we illustrate and contrast the structure of the optimal policy established in Theorem \ref{thm:structure} with and without the joint setup cost. Figure \ref{fig:structure_contrast} illustrates the optimal policy through 10,000 randomly sampled initial conditions for a three-queue system. \jp{By identifying whether or not it is optimal to transfer at each of these states in period 0, we visualize the structure of the optimal policy (i.e., the no-transfer region).} We use parameters $\lambda_i=0.9$, $\mu_i=1$, $h_i=1$ for all $i$ and $r_{12}=r_{21}=2, r_{13}=r_{31}=4, r_{23}=r_{32}=3$, $\tau=5$, $M=5$, and the joint setup cost function with $K=0$ and $K=5$, respectively. The total number of customers in the system in all cases is equal to five. The target states are obtained by solving the associated fluid control problem \eqref{eq:fluid_obj}--\eqref{eq:non_negativity} with a long-enough horizon to empty the system, hence resulting in a stationary policy. \jp{We note that the choice of stationary arrivals and the stationary policy is illustrative and for simplicity, as the structure of the optimal policy is robust under more general arrival rate functions and any horizon length with at least one period.}

In Figure \ref{fig:structure_contrast}, the collection of blue points make up the no-transfer region. \jp{Note that, as established in Theorem \ref{thm:structure}, the region does not consist of multiple disjoint sub-regions.} When there is no setup cost (Figure \ref{fig:structure_no_setup}), we confirm that the target states belong to the boundary of the no-transfer region. In contrast, in the presence of the joint setup cost, target states lie in the relative interior of the no-transfer region (Figure \ref{fig:structure_joint_setup}). 



\begin{figure}[h]
    \FIGURE
    {\begin{subfigure}{0.45\textwidth}
        \FIGURE
        {\includegraphics[width=\textwidth]{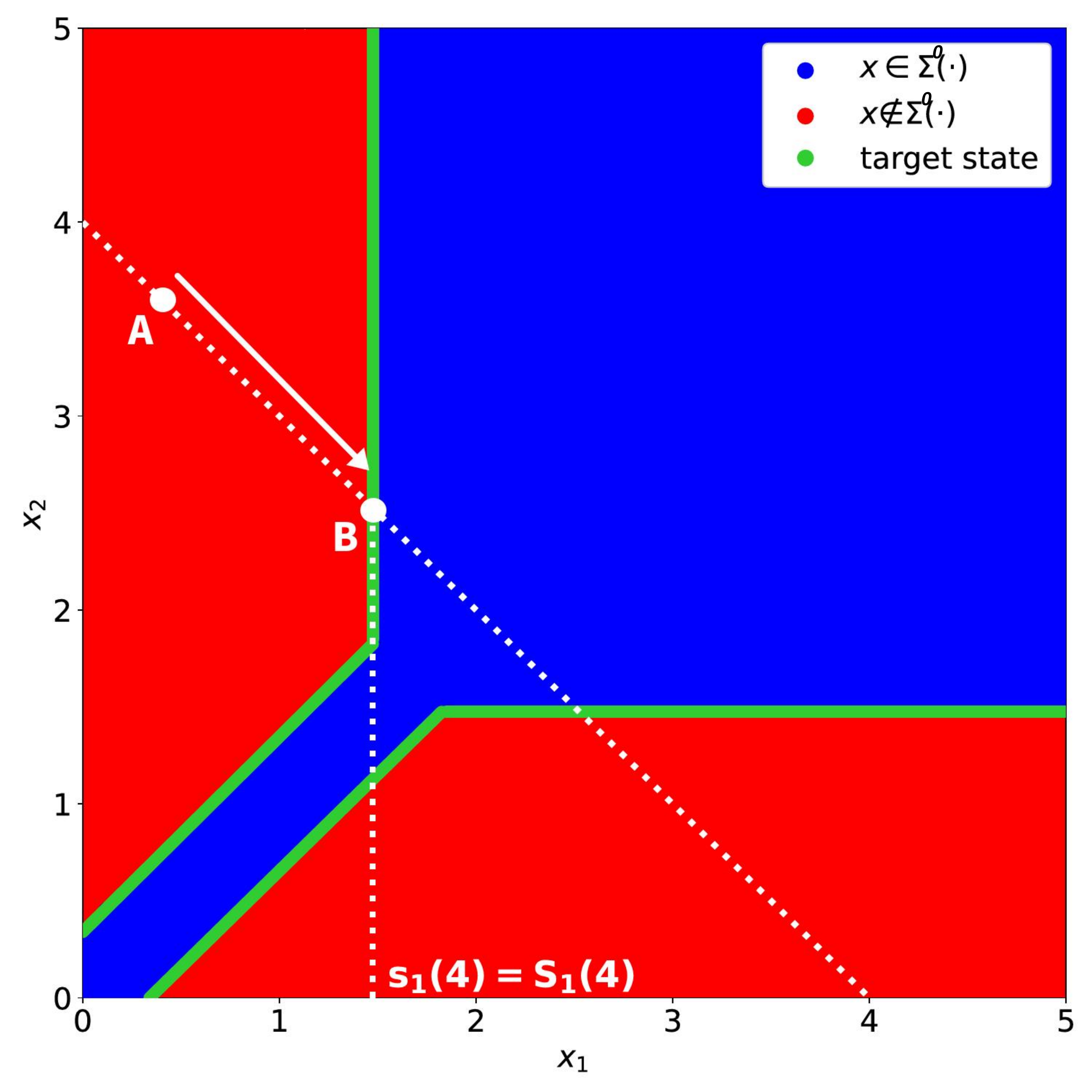}}
        {No setup cost ($K=0$) \label{fig:structure_two_queue_no_setup}}{}
    \end{subfigure}
    \hspace{2em}
    \begin{subfigure}{0.45\textwidth}
        \FIGURE
        {\includegraphics[width=\textwidth]{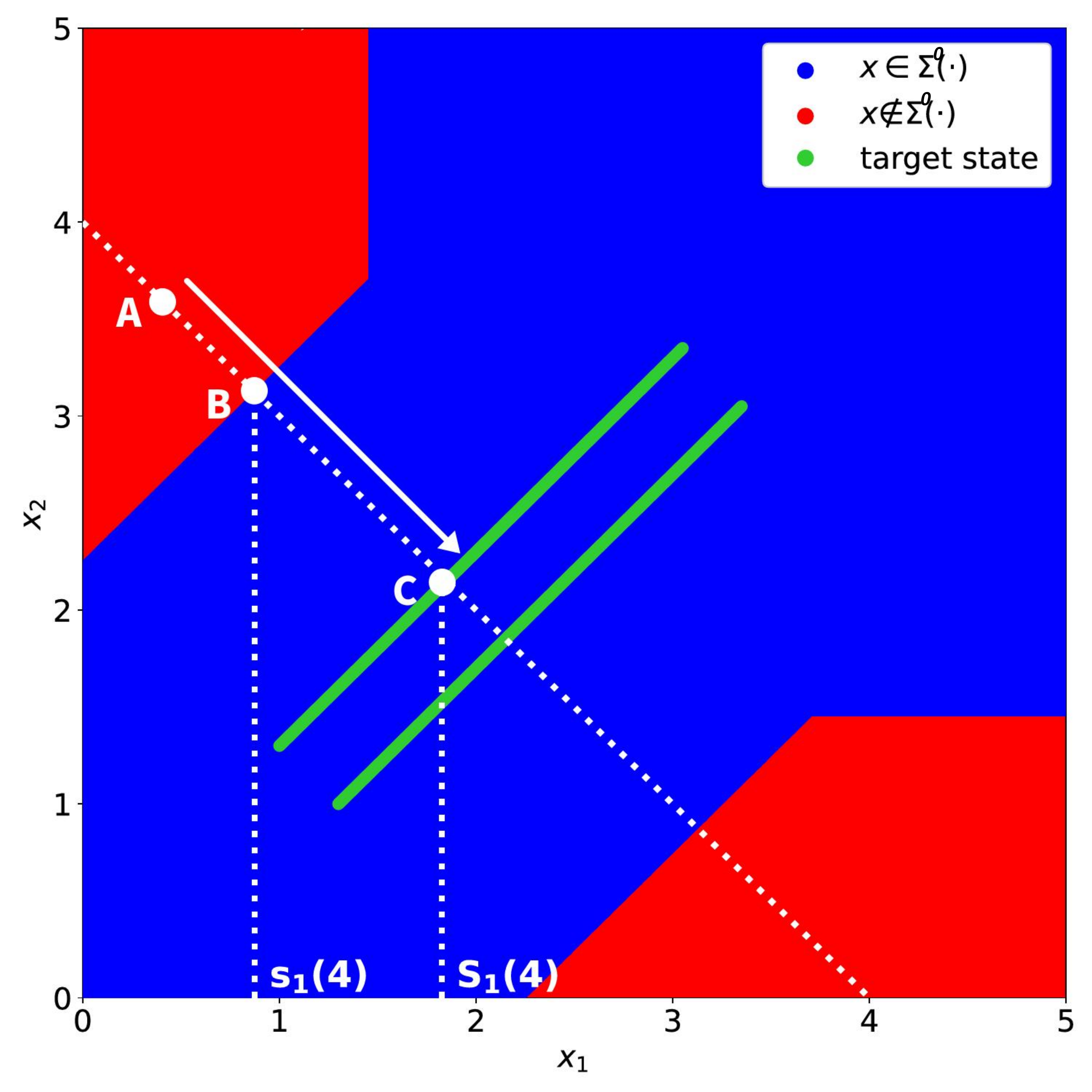}}
        {Joint setup cost ($K=3$) \label{fig:structure_two_queue_setup}}{}
    \end{subfigure}}
    {Structure of the optimal policy \jp{in period 0} for a two-queue system \label{fig:structure_two_queue}}
    {$\lambda_1=\lambda_2=0.7, \mu_1=\mu_2=1, h_1=h_2=1, r_{12}=r_{21}=1, \tau=5, M=10$.}
\end{figure}

Figure \ref{fig:structure_two_queue} presents the optimal structure for a two-queue system for initial conditions in $\{(x_1,x_2): x_1 + x_2 \leq 10, x_1,x_2\geq0\}$. \jp{We consider the subset of this state space for $n=4$ (white dotted line) to illustrate the four parameters $(s_i(n), S_i(n)), i=1,2$.} In Figure \ref{fig:structure_two_queue_no_setup}, at point A, we have $x_1 < s_1(4)$, implying that it is optimal to transfer customers from queue 2 to 1, or along the direction of the white arrow. The target state is point 
B, where $x_1=s_1(4)=S_1(4)$. In contrast, Figure \ref{fig:structure_two_queue_setup} shows that the target state is point C, as opposed to point B, for the same point A. Thus, with a positive setup cost, $s_i(n) < S_i(n)$ holds. Additionally, we verify Proposition \ref{prop:trade_off} in both cases: the non-idleness index $\tau(\mu_i-\lambda_i)$ equals 1.5 for both queues, and given $h_1=h_2$, the no-transfer region (blue) contains the set $\{(x_1,x_2): x_1, x_2 \geq 1.5\}$ (for any $n \geq 3$). (This is also true for Figures \ref{fig:structure_no_setup} and \ref{fig:structure_joint_setup}, but is more difficult to recognize under the current view angle.)

\subsubsection{General Setup Cost Function.} \label{sssec:general_setup}
\begin{figure}[h]
    \FIGURE
    {\begin{subfigure}{0.45\textwidth}
        \FIGURE
        {\includegraphics[width=\textwidth]{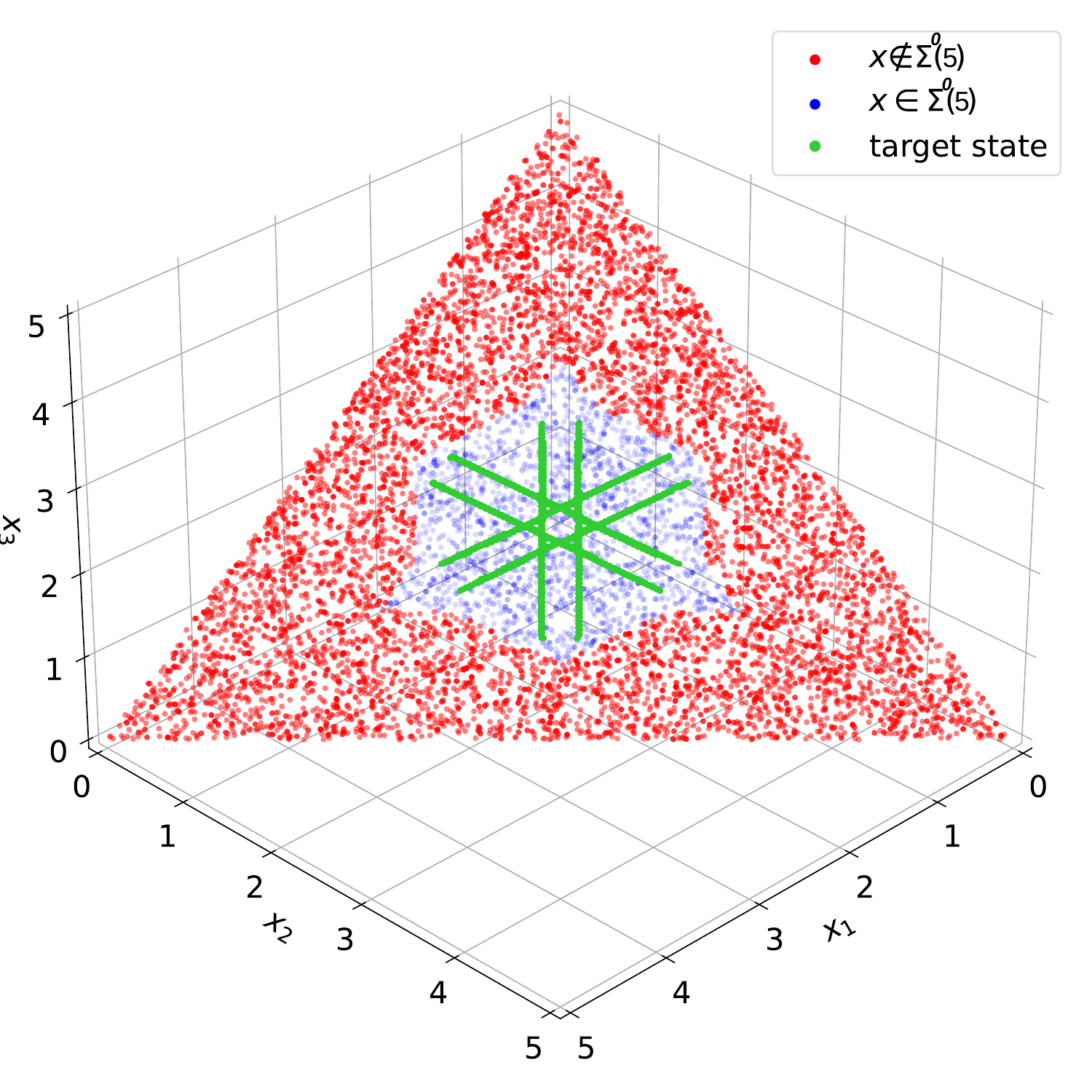}}
        {Pairwise setup cost function \label{fig:structure_pairwise_setup}}{}
    \end{subfigure}
    \hspace{2em}
    \begin{subfigure}{0.45\textwidth}
        \FIGURE
        {\includegraphics[width=\textwidth]{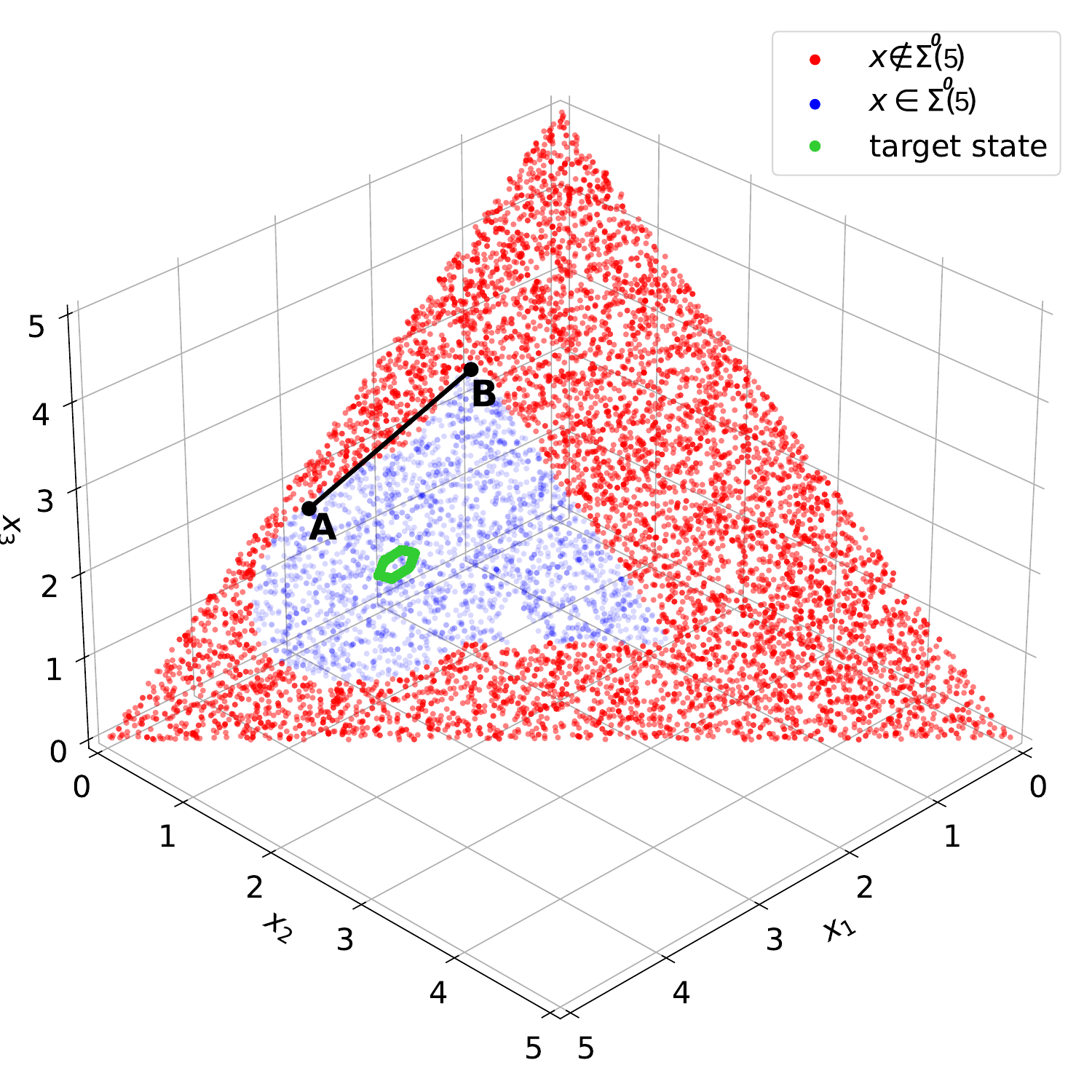}}
        {Non-convexity \label{fig:non_convexity}}{}
    \end{subfigure}}
    {Illustrations of the no-transfer region \label{fig:extensions}}
    {(a): $\lambda_1=\lambda_2=\lambda_3=0.9, \mu_1=\mu_2=\mu_3=1, h_1=h_2=h_3=1, r_{ij}=2, \tilde{K}_{ij}=5, \forall i,j \in \mathcal{N}, i\neq j, \tau=10, M=5$. \\
    (b): $\lambda_1=\lambda_2=0.9, \lambda_3=0.85, \mu_1=1.2, \mu_2=\mu_3=1, h_1=h_2=h_3=1, r_{ij}=1, \forall i,j \in \mathcal{N}, i\neq j, K=3, \tau=5, M=5$.}
\end{figure}

Next, we provide an example in Figure \ref{fig:structure_pairwise_setup} demonstrating that Theorem \ref{thm:structure} may be generalized to the general setup cost function in \eqref{eq:general_setup_cost}. We use the ``pairwise'' setup cost function by setting $\tilde{K}_0=0$ and $\tilde{K}_{ij}=5$ for all $i,j$ for simplicity. For other parameters, we use $\lambda_i=0.9$, $\mu_i=1$, $h_i=1$ for all $i$, $r_{ij}=2$ for all $i,j$, $\tau=10$ and $M=5$. We verify that the no-transfer region is connected and all target states lie in its relative interior.

\subsubsection{Non-Convexity of the No-Transfer Region.}

\jp{Although non-convexity of the no-transfer region is already apparent in Figure \ref{fig:structure_pairwise_setup}, it can be non-convex even for the simpler case of stationary arrival rates, linear holding costs, and joint setup cost, as illustrated in Figure \ref{fig:non_convexity}. Figure \ref{fig:non_convexity} shows two points within the no-transfer region, A and B, whose convex combination (the line segment) is not fully contained within the region. Note that as we go from A to B, the number of customers increases at queue 2 but decreases at queue 1. At some midpoint between A and B, both queues 1 and 2 have sufficiently low numbers of customers and transferring becomes optimal --- in this case queue 3 would send to both queues 1 and 2. However, at either A or B, only one of queues 1 and 2 has a sufficiently low number of customers, but not both, and transferring is not worthwhile given the setup cost. Non-convexity generally makes computing the no-transfer region more challenging, as finding two points in it does not imply that the whole line segment connecting the two is also in the region.}

\section{Approximate Dynamic Programming (ADP) \label{sec:ADP}}
\jp{So far, we have established that the optimal fluid policy is of the \textit{region-of-inaction} type, which partitions the state space into a compact, connected no-transfer region and its complement. In Section \ref{sec:stochastic_system}, we present numerical evidence suggesting that this structure also holds for the optimal policy of the stochastic problem. Motivated by these results, we propose a simulation-based approximate policy iteration (API) algorithm with four key components: (i) a binary classifier to characterize the region-of-inaction, allowing us to bypass target state computations whenever a state is predicted to lie within the region; (ii) use of Common Random Numbers (CRN) and coupling to update value functions; (iii) using the optimal fluid policy for initialization; and (iv) a feasibility check for verifying and preserving the connectedness of the region.}

\jp{We first present a post-decision state DP formulation \citep[Section 6.4.]{powell2007approximate} of the stochastic control problem described in Section \ref{sec:problem_definition}. With a slight abuse of notation, we denote the realizations of the pre- and post-transfer states in period $m$ by $x^m$ and $y^m$, respectively. Let $J^m(x^m)$ be the value function, i.e., the minimum expected cost from period $m$ onward given that the state is $x^m$ at time $t_m^-$. For $m\in \mathcal{M}$, the value function satisfies the optimality equation
\begin{equation} \label{eq:post_decision_opt_eqn}
    J^m(x^m) = \min_{y^m \in\Delta^{\mathbb{Z}}(\mathrm{e}^\top x^m)} \left[C(y^m - x^m) + J^m_a(y^m)\right],
\end{equation}
where $\Delta^{\mathbb{Z}}(n) = \{y \in \mathbb{Z}_+^N: \mathrm{e}^\top y = n\}$, $J^M \equiv 0$, and
\begin{equation} \label{eq:post_decision_value_fn}
    J^m_a(y^m) = \mathbb{E}\left[\int_{t_m}^{t_{m+1}}h\left( X(s) \right)ds + J^{m+1}\left( X(t_{m+1}^-) \right) \bigg| X(t_m)=y^m \right].   
\end{equation}
A conventional method to solve \eqref{eq:post_decision_opt_eqn} is the policy iteration algorithm. However, the large state and action spaces and the difficulty in computing the expectation in \eqref{eq:post_decision_value_fn} make this impractical. Instead, we use simulation to approximate \eqref{eq:post_decision_value_fn} and a classifier to characterize the region-of-inaction. This allows us to solve \eqref{eq:post_decision_opt_eqn} more efficiently, either by skipping it entirely, or restricting the feasible set to the interior of the region-of-inaction.} 

\jp{\subsection{The Proposed API Algorithm}
Denote the state space by $\mathcal{X} \subset \mathbb{Z}_+^{N \times N}$, which represents a truncated system with a maximum size of $n_{\mathrm{max}}$, i.e., $\mathcal{X} = \cup_{n=0}^{n_{\mathrm{max}}} \Delta^{\mathbb{Z}}(n)$. Let $g^m:\mathcal{X} \to [0,1]$ denote a classifier in period $m$ mapping a state to its \textit{probability of belonging} to the region-of-inaction (referred to as state probability). Then an approximate characterization of the region-of-inaction is
\begin{equation} \label{eq:approx_roi}
    \tilde{\Sigma}^m =\{x\in \mathcal{X}: g^m(x) \geq p\},
\end{equation}
representing all states deemed likely to be in the region based on a pre-defined probability threshold $p \in [0,1]$. Let $\mathcal{N}(x) \equiv \{x' \in \mathcal{X}: \|x-x'\|_{\infty} \leq 1\}$ be the neighbourhood of $x$. We then define the ``boundary'' of $\tilde{\Sigma}^m$ as
\begin{equation} \label{eq:boundary}
    \partial \tilde{\Sigma}^m \equiv \{x \in \mathcal{X}: \mathcal{N}(x) \cap \tilde{\Sigma}^m \neq \varnothing \text{ and } \mathcal{N}(x) \cap (\mathcal{X} \setminus \tilde{\Sigma}^m) \neq \varnothing \},
\end{equation}
i.e., all states whose neighbourhood contains states from both inside and outside the region-of-inaction. Finally, denote the policy by $\pi = (\pi^0, \ldots, \pi^{M-1})$ where
\begin{equation}\label{eq:policy}
\pi^m(x) = \begin{cases}
        x, & \mbox{if } x\in \tilde{\Sigma}^m \setminus \partial \tilde{\Sigma}^m; \\
        y(x), & \mbox{otherwise}.
     \end{cases}
\end{equation}
In other words, if $x$ lies in the ``interior'' ($\tilde{\Sigma}^m \setminus \partial \tilde{\Sigma}^m$), then taking no action is optimal; otherwise, the policy selects $y(x)$, which is an optimal solution to,
\begin{equation} \label{eq:solve_target}
    \begin{split}
        \min_{y} &\quad C(y-x) + \Bar{V}^{m}(y) \\
        \mathrm{s.t.} &\quad \mathrm{e}^\top y = \mathrm{e}^\top x, \\
        &\quad y \in \tilde{\Sigma}^m \cup (\{x\} \cap \partial \tilde{\Sigma}^m),
    \end{split}
\end{equation}
where $\Bar{V}^{m}(y)$ is the sample average approximation of the expectation in \eqref{eq:post_decision_value_fn}. The last constraint states that when $x$ is part of the boundary ($\partial \tilde{\Sigma}^m$), we also consider $x$ itself as a candidate. We now provide details of the algorithm below. The pseudocode is available in Appendix \ref{appen:pseudocode}.}

\jp{\textbf{1. Initialization.} For each $x \in \mathcal{X}$ and period $m \in \{0,\ldots,M-1\}$, the value function $\Bar{V}_0^{m}(x)$ is initialized to the optimal cost of the $(M-m)$-period fluid control problem with $x$ as the initial condition (see Appendix \ref{appen:numerical_solution_approach} for the numerical solution approach). Upon obtaining the target state $y$ (optimal post-transfer state), the label of $x$ is initialized as $\mathrm{label}_0^m(x)= 1\{x = y\}$. Using $\{\mathrm{label}_0^m(x): x \in \mathcal{X}\}$ as the target variable, a classifier $g^m_0$ is trained for each $m$. We discuss the advantage of the fluid initialization compared to a more naive method in Appendix \ref{appen:initialization}.}

\jp{\textbf{2. Policy evaluation.} In each iteration $j$, the algorithm performs $B$ simulation runs for each $x \in \mathcal{X}$ and generates $\{l^m(x), \Bar{v}^m(x)\}_{m=0}^{M-1}$, a set of observed labels and value function estimate of $x$ in each period it is encountered. (For $x\in \mathcal{X}$ not encountered, these sets are empty.) The value function estimate $\Bar{v}^m(x)$ is computed as the average of the observed values across all instances where $x$ is encountered in period $m$. Each simulation run follows the policies $\pi_{j-1}^0, \ldots, \pi_{j-1}^{M-1}$, which are obtained after the policy improvement step in the previous iteration $j-1$. Following $\pi^m_{j-1}$ may entail solving \eqref{eq:solve_target} in period $m$ using $\Bar{V}^m_{j-1}$ in the objective function to determine a target state. The system then evolves according to $x^{m+1} = \pi_{j-1}^m(x^m) + a^m - d^m$, where $a^m, d^m \in \mathbb{Z}_+^N$ are realizations of arrivals and departures over period $m$. In period $m$, we assign the label $l^m(x) = 1\{x = \pi_{j-1}^m(x)\}$. Finally, $B$ should be large enough to produce reliable value function estimates but not so large as to slow the algorithm's iterative process excessively. In our numerical experiments in Section \ref{sec:adp_two_queues} and \ref{sec:case_study}, we use $B=10$.} 

\jp{\textbf{3. Policy improvement.} Let $\mathcal{X}^{m,\mathrm{visited}}_j$ denote the set of states encountered at least once in period $m$ in any simulation run. Given the new observations $\{l^m(x), \Bar{v}^m(x)\}_{m=0}^{M-1}$, the algorithm updates the value function in period $m$ by 
\begin{align}
    \Bar{V}_{j+1}^{m}(x) = 
    \begin{cases}
        \frac{u^m(x)}{u^m(x)+1} \Bar{V}_{j}^{m}(x) + \frac{1}{u^m(x)+1}\Bar{v}^m(x), & \mbox{if } x \in \mathcal{X}^{m,\mathrm{visited}}_j, \\[1.5ex]
        \bar{V}^m_j(x), & \mbox{otherwise}.
    \end{cases}
\end{align}
The parameter $u^m(x)$ tracks the number of updates performed at $x$ in period $m$, starting at 1 in iteration 0 (initialization), and increments by at most 1 in each iteration. For each encountered state, the label is updated as $\mathrm{label}_{j+1}^m(x) = l^m(x)$ (otherwise, remains as $\mathrm{label}_{j}^m(x)$). Lastly, using $\{\mathrm{label}_{j+1}^m(x): x\in \mathcal{X}\}$ as the target variable, a new classifier $g_{j+1}^m$ is trained for each $m$. The algorithm terminates after $j_{\mathrm{max}}$ iterations.}

\jp{We note that the algorithm permits label changes only for boundary states. For each boundary state, two outcomes are possible when solving \eqref{eq:solve_target}: (1) $x \notin \tilde{\Sigma}^m_j$ but $y(x) = x$, implying that $x$ should be added to $\tilde{\Sigma}^m_j$; and (2) $x \in \tilde{\Sigma}^m_j$ but $y(x) \neq x$, implying that $x$ should be removed from $\tilde{\Sigma}^m_j$. In contrast, labels are fixed for non-boundary states: we set $l^m(x) = 1$ for all $x \in \tilde{\Sigma}^m_j \setminus \partial \tilde{\Sigma}^m_j$, and $l^m(x) = 0$ for all $x \in \mathcal{X} \setminus (\tilde{\Sigma}^m_j \cup \partial \tilde{\Sigma}^m_j)$. Thus, the algorithm is designed to iteratively refine its approximation of the region-of-inaction's boundary.}

\begin{figure}[]
    {\color{black}
    \FIGURE
    {
    \centering
    \begin{tabular}{ccc}
        \begin{subfigure}{0.3\textwidth}
            \caption{Iteration 0}
            \includegraphics[width=\textwidth]{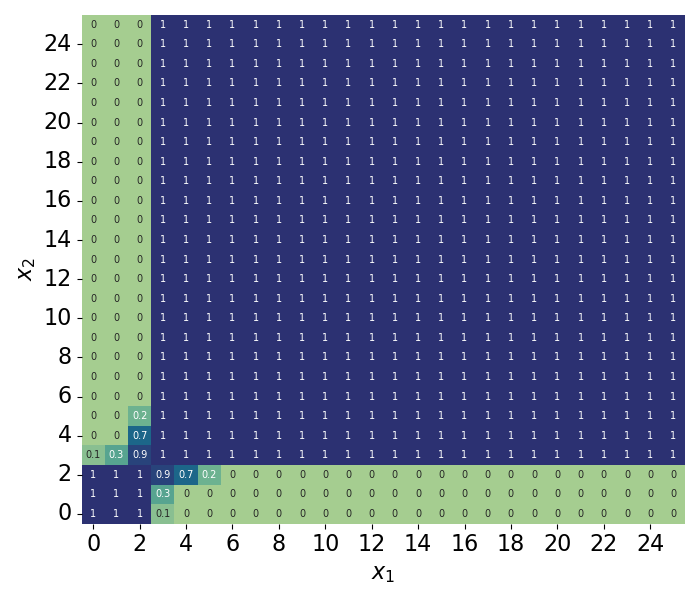}
        \end{subfigure} &
        \begin{subfigure}{0.3\textwidth}
            \caption{Iteration 2}
            \includegraphics[width=\textwidth]{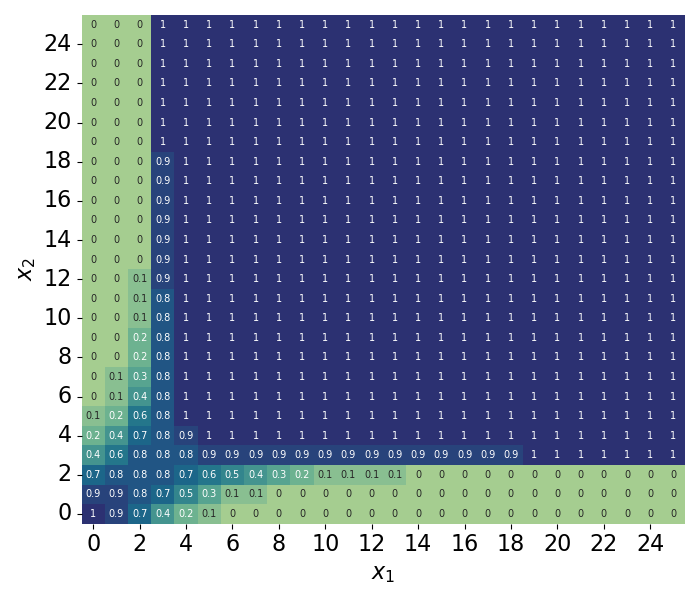}
        \end{subfigure} &
        \begin{subfigure}{0.3\textwidth}
            \caption{Iteration 4}
            \includegraphics[width=\textwidth]{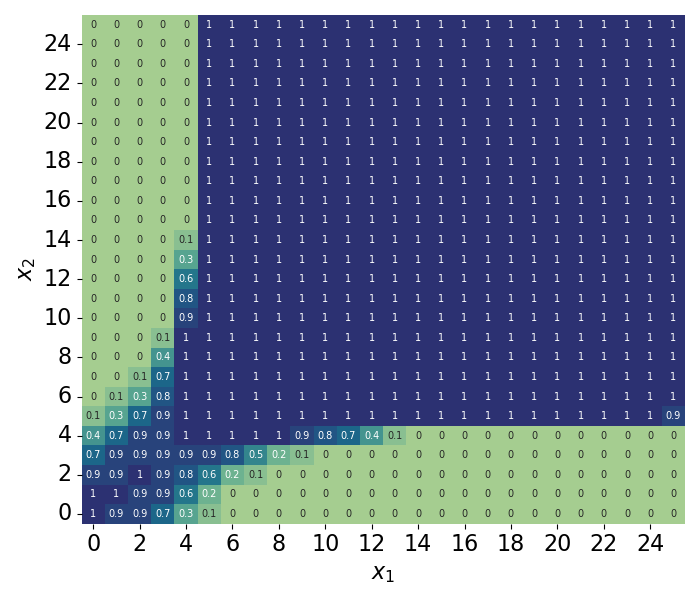}
        \end{subfigure} \\
        \begin{subfigure}{0.3\textwidth}
            \caption{Iteration 6}
            \includegraphics[width=\textwidth]{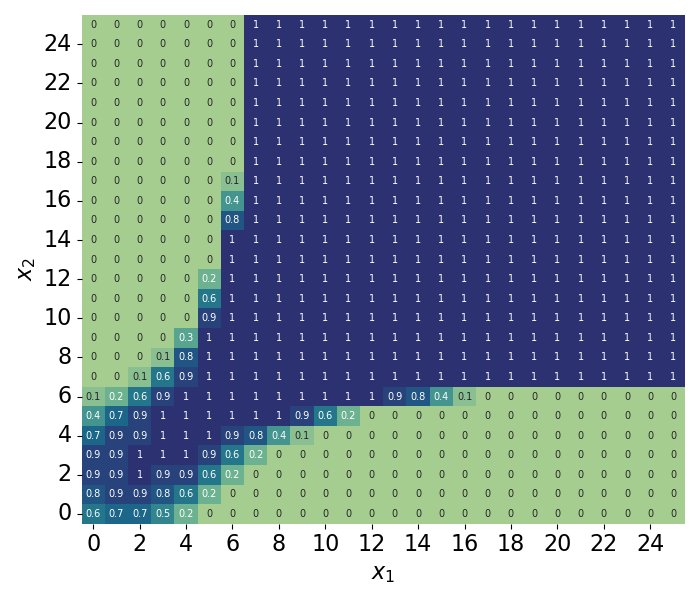}
        \end{subfigure} &
        \begin{subfigure}{0.3\textwidth}
            \caption{Iteration 8}
            \includegraphics[width=\textwidth]{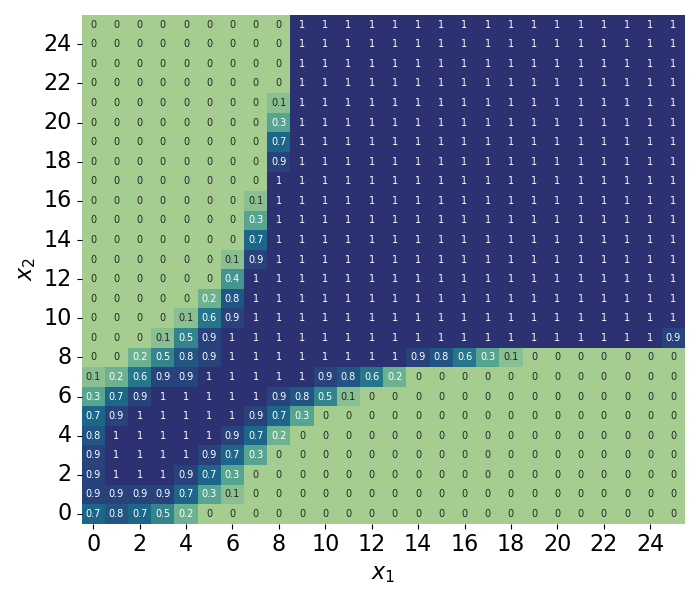}
        \end{subfigure} &
        \begin{subfigure}{0.3\textwidth}
            \caption{Iteration 10}
            \includegraphics[width=\textwidth]{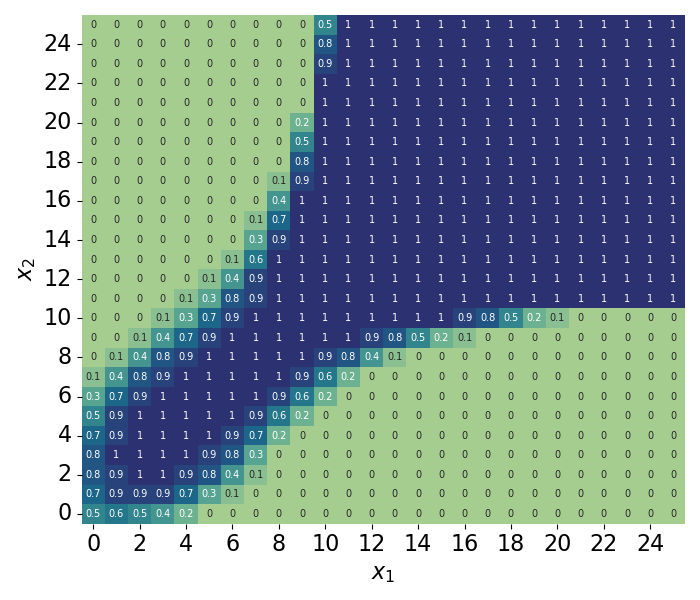}
        \end{subfigure}
    \end{tabular}
    }%
    {Convergence of the API policy for a system of two $M/M/1$ queues \label{fig:adp_convergence}}%
    {$\lambda=(0.9,0.9), \mu=(1,1), \tau=1, M=7, h=(10,10), r_{12}=r_{21}=1, K_{12}=K_{21}=1$. Iteration 0: fluid policy.}
    }
\end{figure}

\jp{In general, convergence guarantees for simulation-based algorithms are difficult to establish \citep{bertsekasApproximatePolicyIteration2011}. For the numerical experiments in Section \ref{sec:adp_two_queues}, where the region-of-inaction can be visualized, we verify that the state probabilities exhibit less oscillations over time and the state labels appear to converge. An example is provided in Figure \ref{fig:adp_convergence} for a two-queue system, starting with the fluid initialization. }


\jp{\subsection{How Does the Algorithm Leverage the Structure?}
We next key outline how the algorithm leverages the structure of the optimal policy.}

\subsubsection{Bypassing computation.} \jp{When $x$ is in the relative interior, $\tilde{\Sigma}^m_j \setminus \partial \tilde{\Sigma}^m_j$, we set $\pi^m_j(x)=x$ and bypass solving \eqref{eq:solve_target}. When $x \notin \tilde{\Sigma}^m_j$, we restrict the feasible set to $\tilde{\Sigma}^m_j \cup (\{x\} \cap \partial \tilde{\Sigma}^m_j)$. Second, we note that if two distinct states $x^1$ and $x^2$ share the same target state $y$, their value functions differ only by their transfer costs in the current period: $J^m(x^1) - J^m(x^2) = C(y-x^1) - C(y-x^2)$. Thus, evaluating one state provides a value for the other without requiring additional sampling or computation in subsequent periods. We enforce this by using CRN, which ensure the two systems follow identical sample paths after coupling at $y$. We numerically illustrate the value of using CRNs numerically in Appendix \ref{appen:value_crns}. Finally, since policies are deterministic and independent of simulation runs, if two sample paths reach the same state in any period, they must share the same target state at that state. We exploit this by reusing target states across simulation runs.}
\subsubsection{Choice of features for the classifier.}
\jp{The structure also informs the choice of features used in training the classifier. The first includes up to third-order polynomials of the queue lengths $x \in \mathbb{Z}_+^N$, along with interaction terms among $N-1$ components of $x$. The second includes up to third-order polynomials of the \textit{distance features} $d_i(x) = \|x - v^i\|$, where $v^i=(0,\ldots,n,\ldots,0)$ is the $i$th vertex of the state space with its $i$th component equal to $n=\mathrm{e}^\top x$ and all others zero. We also include interactions among $N-1$ components of the distance vector $d(x) = (d_1(x), \ldots, d_N(x))$. The distance features help capture the geometry of the region-of-inaction precisely because the region is connected and cannot consist of distinct sub-regions. Indeed, Figure~\ref{fig:feature_comparison} in Appendix \ref{appen:features} illustrates that both feature sets are necessary for high classification accuracy, and using either set alone is generally insufficient.}

\subsubsection{Preserving connectedness.}\label{ssec:preserve_connectedness}
\jp{A key challenge in design of the algorithm is to ensure that the region-of-inaction remains connected over successive iterations. In particular, the choice of the probability threshold $p$ in our characterization \eqref{eq:approx_roi} plays a critical role, as certain choices may lead to a disconnected region; see Appendix \ref{appen:choice_prob_threshold} for additional discussion. In this section, we provide guidance on selecting an appropriate probability threshold and later propose a check for connectedness for the algorithm.}

\jp{First, we propose the following condition on the classifier $g^m_j$.}
\jp{\begin{assumption} \label{ass:classifier}
    For all $m\in \{0,\ldots,M-1\}$ and $j\in \{0,\ldots,j_{\mathrm{max}}\}$, $g^m_j$ is continuous, bounded away from 0, and piecewise-monotone with finitely many pieces.
\end{assumption}}

\noindent\jp{This assumption is not restrictive. Many classifiers can be expressed as $\sigma(\beta^\top f(x))$, where $f(x)$ is a vector of features or (non-linear) transformations of the state, $\beta$ is a vector of model parameters or weights, and $\sigma$ maps to a probability value, with common choices including the sigmoid function. In our application, we consider logistic regression with features that are (Lipschitz) continuous and have finitely many critical points, which satisfies Assumption \ref{ass:classifier}.}

\jp{The key observation is that because our algorithm updates only the boundary states, it guarantees connectedness as long as the boundary labels remain connected with the rest of the region-of-inaction in each iteration.}

\jp{
\begin{proposition} \label{prop:connectedness}
    Under Assumption \ref{ass:classifier}, there exists $p > 0$ such that the set of boundary labels (with label 1) remains connected with the region-of-inaction in each iteration.
\end{proposition}
Note that setting $p=0$ trivially guarantees connectedness by treating the entire state space as the region-of-inaction in every iteration. In contrast, Proposition \ref{prop:connectedness} shows that there exist non-trivial probability thresholds that preserve connectedness while allowing meaningful updates. Its proof and the accompanying discussion in Appendix \ref{appen:proof_connectedness} suggest that relatively small values (e.g., $p \leq 0.5$) should work well, and properties such as Lipschitz continuity or monotonicity of $g^m_j$ offer more flexibility in selecting $p$. In practice, however, the exact value of a suitable threshold may still be difficult to find. In Appendix \ref{appen:check_for_connectedness}, therefore, we provide a practical check for connectedness that can be incorporated into the algorithm.}

\section{Numerical Experiments} \label{sec:numerical_experiments}
\jp{In this section, we first use numerical examples to demonstrate that the optimal policy structure established for the fluid control problem also holds for the stochastic control problem. We then evaluate the performance of the API algorithm for the stochastic system using simulation experiments. Lastly, we conduct a case study on inter-facility patient transfers to quantify the potential benefits of different transfer policies in a practical setting.}

\subsection{Comparison to the Optimal Policy} \label{sec:stochastic_system}
We first examine the structure of the optimal policy for the stochastic control problem. The discrete-time stochastic control problem can be modeled as a Markov decision process (MDP). However, even for small systems, solving the discrete-control MDP is both computationally and analytically hard. This is mainly due to the complexity of computing transition probabilities compared to the continuous-control MDP, where one can apply uniformization to obtain a discrete-time MDP with simple transition probabilities. Our approach here is therefore to solve a continuous-control MDP instead and compare the structure of the continuous-time MDP policy to that of the fluid policy.

\jp{In the uniformization approach (see details in Appendix \ref{appen:MDPsolutions}), we designate $(0,0)$ as an absorbing state, and the MDP policy solves for the quickest way to reach an empty state. Thus, to derive an (approximately) continuous-control fluid problem that is comparable to the continuous-control MDP, we first set the length of each period $\tau$ to be equal to the average time between two successive events (arrival or service completion), and set the length of the horizon to $M(x^0, \tau) = \max_{i=1,2}\{x^0_i/(\tau(\mu_i-\lambda_i))\}$ (where $x^0$ denotes a given initial condition) to ensure a long-enough horizon to empty the system (assuming $\mu_i > \lambda_i$ for $i=1,2$).}

\begin{figure}[h]
    \FIGURE
    {\begin{subfigure}{0.4\textwidth}
        \FIGURE
        {\includegraphics[width=\textwidth]{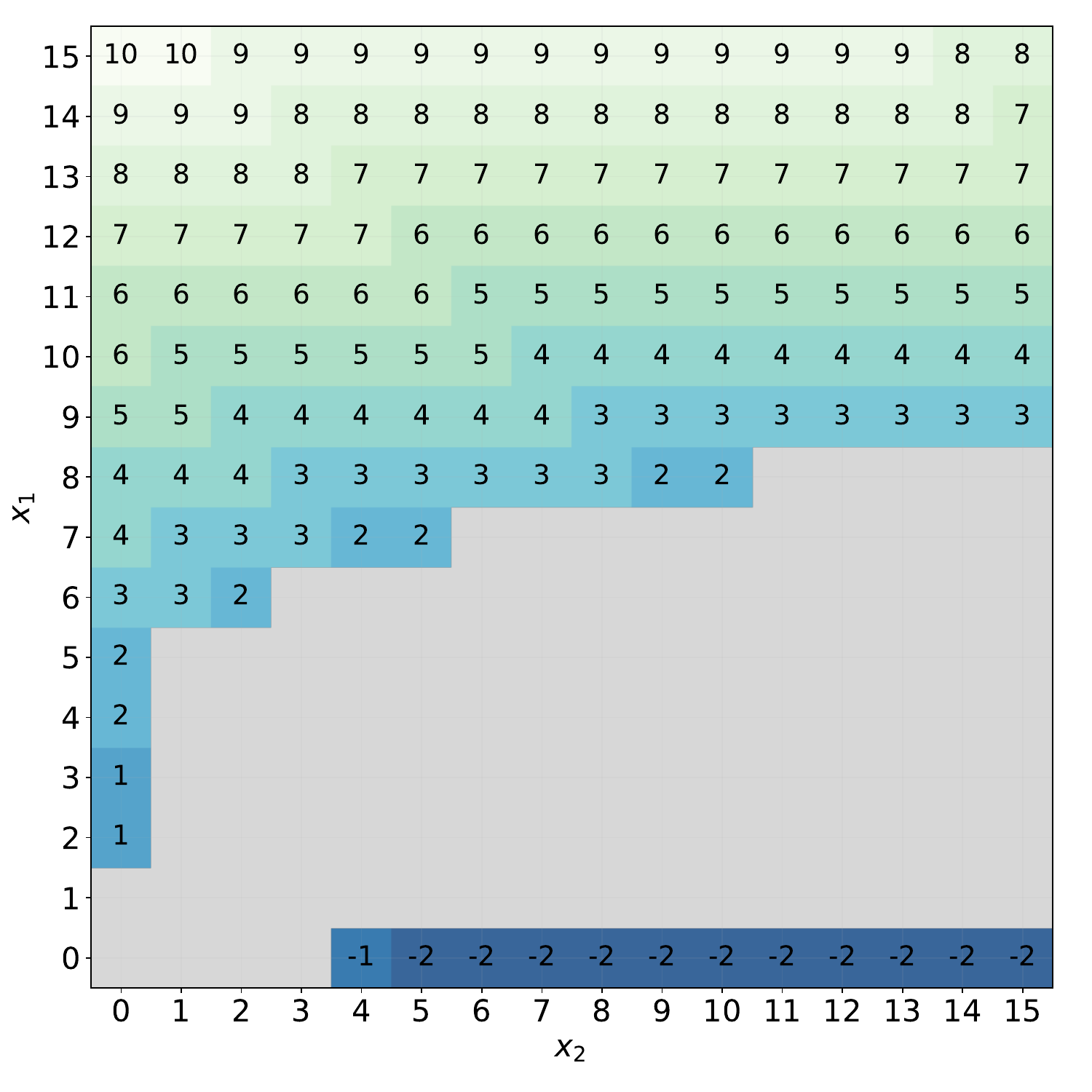}}
        {MDP policy \label{fig:MDP_policy}}{}
    \end{subfigure}
    \hspace{2em}
    \begin{subfigure}{0.4\textwidth}
        \FIGURE
        {\includegraphics[width=\textwidth]{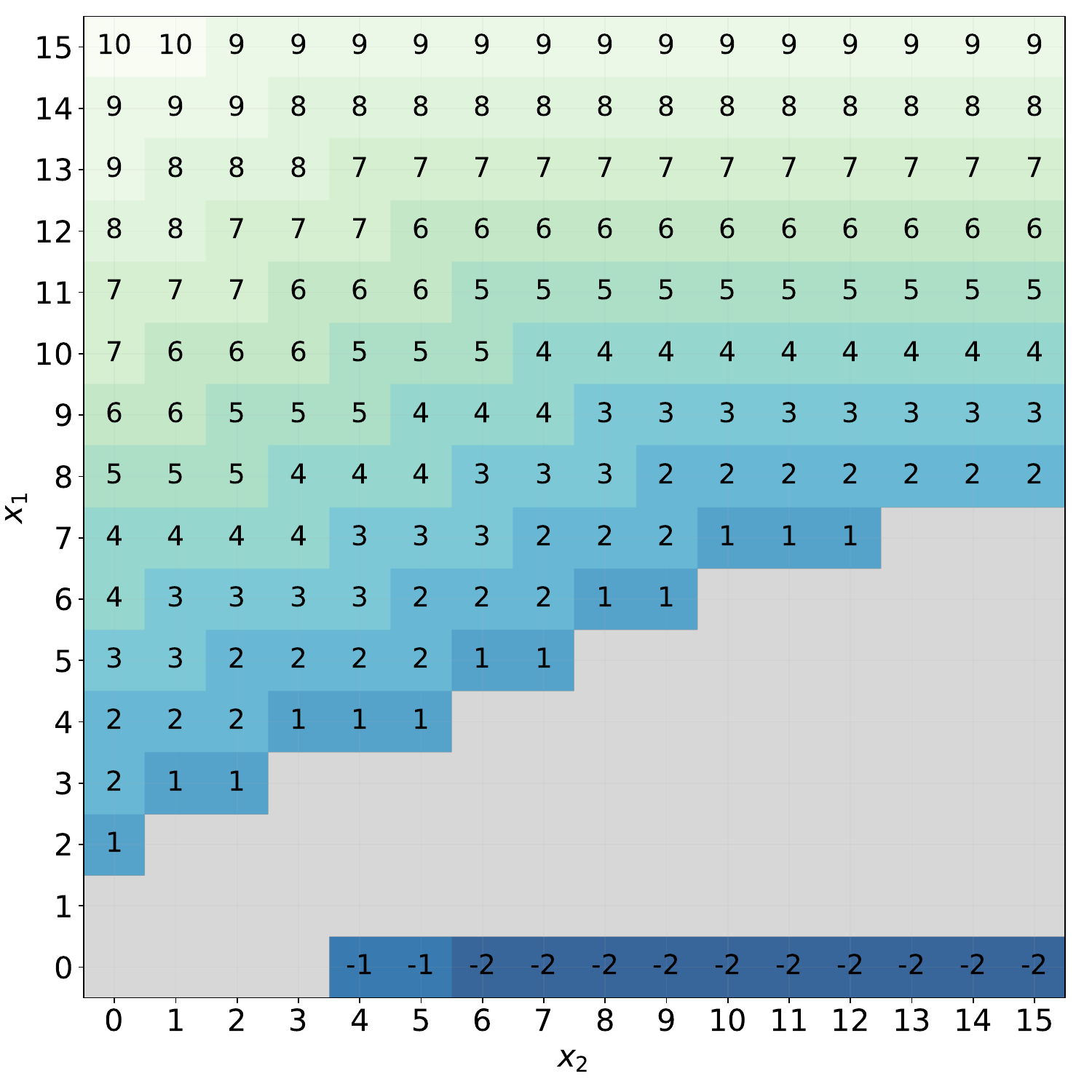}}
        {Fluid policy \label{fig:fluid_policy}}{}
    \end{subfigure}}
    {Example comparison of an MDP policy and a fluid policy  \label{fig:compare_policy_unequal_light_traffic}}
    {$\lambda_1=0.9, \lambda_2=1.5, \mu_1=1.5, \mu_2=2.5, h_1=1.3, h_2=1, r_{12}=r_{21}=1, K=1$.}
\end{figure}

In general, the optimal transfer policy for the stochastic system has the same structure as the fluid policy. One example is provided in Figure \ref{fig:compare_policy_unequal_light_traffic}, which presents an MDP policy (left) and the fluid policy (right), where positive (negative) numbers indicate transfers from queue 1 to 2 (2 to 1). We find that the structure remains consistent, as evident from the connectedness of the no-transfer region (grey). We further note that there are constant re-order and order-up-to points for each fixed $n$, i.e., the total number of customers. \jp{However, the exact values of these parameters can differ. In Appendix \ref{appen:add_results_stochastic_system}, we show that despite this, the optimality gap of the fluid policy is small.} 

\subsection{API Performance for Two-Queue Systems}\label{sec:adp_two_queues}
\jp{To better understand the performance of the policy learned by our API algorithm (API policy), we focus on simple two-queue systems. We compare the API policy against three benchmarks: no-transfer, Myopic, and fluid. Comparing Myopic to fluid provides insights into the value of being forward-looking, while comparing fluid to API highlights the additional value of accounting for stochasticity.}

\jp{The policies are obtained as follows. The Myopic policy, based on Proposition \ref{prop:trade_off}, transfers just enough to bring any queue with $x_i < \tau (\mu_i - \lambda_i)^+$ up to $\tau (\mu_i - \lambda_i)^+$, thereby avoiding idleness for one period. The fluid policy is obtained by solving the fluid control problem in Appendix \ref{appen:numerical_solution_approach}. The API policy is computed by running Algorithm \ref{alg:API} using logistic regression with the features described in Section \ref{sec:ADP}, probability threshold $p=0.1$, $B=10$ simulation runs per state, and $j_{\mathrm{max}}=10$ iterations. We consider the following three systems with increasing variability in system dynamics:
\begin{enumerate}
    \item $M/M/1$ queues under heavy traffic with $\rho = (\lambda_1 + \lambda_2)/(\mu_1 + \mu_2) = 0.9$;
    \item $M/G/1$ queues with $\rho = 0.9$ and log-normal service times having three times the standard deviation of the exponential distribution (but same mean);
    \item $M(t)/G/1$ queues with $\rho = 0.9$, log-normal service times as above, and $\lambda_i(t) = \lambda_i(1+0.5\sin(2\pi t - \pi))$ for $i=1,2$, which ranges between $0.5\lambda_i$ and $1.5\lambda_i$.
\end{enumerate}
The remaining parameters are specified under Table \ref{tab:mm1_policy_performance}.}

\jp{The performance generally depends on the initial condition. We focus on large, imbalanced states where the fluid and API policies disagree --- specifically, states that lie on the boundary of the fluid policy’s no-transfer region but fall outside the API policy’s. As shown in Figure \ref{fig:adp_convergence}, the API algorithm learns a smaller no-transfer region, resulting in more proactive transfers. We present results for the $M/M/1$ system here (Table \ref{tab:mm1_policy_performance}), and relegate the rest to Appendix \ref{appen:add_results_adp_two_queues}.}

\jp{We make several observations. First, all policies, including Myopic, outperform the no-transfer policy, with the API policy consistently achieving the best performance. Second, there is significant value to being forward-looking and accounting for stochasticity. Specifically, the fluid policy improves over Myopic by 5.8--6.7 percentage points (pp) on average, while the API policy achieves an \textit{additional} 6.5--8.3pp improvement over the fluid policy. Lastly, these gains generally increase with system variability and in more critically loaded settings. For example, under log-normal service times, the fluid policy outperforms Myopic by 5.9--7.3pp, and the API policy yields a further 6.1--10.5pp improvement. For the $M(t)/G/1$ system, the respective gains can reach up to 10.3pp and 11.5pp on average.}
\begin{table}
\renewcommand*{\arraystretch}{1.15}
\TABLE
{Performance of Myopic, fluid, and API policies relative to no-transfer for $M/M/1$ system \label{tab:mm1_policy_performance}}
{
{\color{black}\begin{tabular}{clrrl}
\hline
Initial condition & Policy & Holding cost & Transfer cost & Reduction (\%) \\
\hline
\multirow{3}{*}{(1, 15)} & Myopic      & 1048.1 & 3.5 & 3.9 $\pm$ 0.3\% \\
                         & Fluid       & 978.6 & 7.7 & 9.7 $\pm$ 0.8\% \\
                         & API         & 879.0 & 22.8 & \textbf{16.2} $\pm$ \textbf{1.3}\% \\
\hline
\multirow{3}{*}{(1, 17)} & Myopic      & 1170.8 & 3.5 & 3.7 $\pm$ 0.3\% \\
                         & Fluid       & 1090.4 & 8.0 & 10.0 $\pm$ 0.8\% \\
                         & API         & 973.3 & 23.3 & \textbf{17.4} $\pm$ \textbf{1.2}\% \\
\hline
\multirow{3}{*}{(1, 19)} & Myopic      & 1297.9 & 3.4 & 3.5 $\pm$ 0.3\% \\
                         & Fluid       & 1207.8 & 8.2 & 10.0 $\pm$ 0.7\% \\
                         & API         & 1081.7 & 25.5 & \textbf{17.3} $\pm$ \textbf{1.2}\% \\
\hline
\multirow{3}{*}{(1, 21)} & Myopic      & 1427.6 & 3.4 & 3.4 $\pm$ 0.2\% \\
                         & Fluid       & 1328.5 & 8.5 & 10.0 $\pm$ 0.7\% \\
                         & API         & 1185.1 & 24.8 & \textbf{18.3} $\pm$ \textbf{1.1}\% \\
\hline
\multirow{3}{*}{(1, 23)} & Myopic      & 1560.2 & 3.4 & 3.1 $\pm$ 0.2\% \\
                         & Fluid       & 1454.4 & 8.6 & 9.8 $\pm$ 0.7\% \\
                         & API         & 1300.8 & 26.5 & \textbf{17.9} $\pm$ \textbf{1.1}\% \\
\hline
\end{tabular}}
}
{\emph{Note.} $\lambda=(9,9), \mu=(10,10), \tau=1, M=7, h=(10,10), r_{12}=r_{21}=1, K_{12}=K_{21}=1$.}
\end{table}

\subsection{Case Study: Inter-Facility Patient Transfer} \label{sec:case_study}
In this section, we conduct a case study using a simulation model calibrated with data from four hospitals in the Greater Toronto Area during the COVID-19 pandemic; see \cite{chanOptimizingInterHospitalPatient2023a} for additional details on the data. We evaluate and compare \jp{three transfer policies: the Myopic, fluid, and API policies.} Specifically, we evaluate the policies for systems with multiple servers, log-normally distributed service times, with restrictions on the number of transfers, and under both non-stationary arrivals \jp{and prediction errors for the arrival rates.}

\textbf{Simulation model.}
The simulation model consists of four parallel multiserver queues. The servers represent beds in the intensive care units (ICU) and queues represent boarding from the acute ward or emergency department (ED). We note that keeping patients in the ward or ED until ICU capacity becomes available (as opposed to diverting) was common during the pandemic, see, e.g., \cite{bellaniNoninvasiveVentilatorySupport2021, douinICUBedUtilization2021}. Patients arrive to each queue according to a non-homogeneous Poisson process with piecewise-constant rates (varying by day of the week), \jp{and service times are exponentially distributed.} \jp{To approximate multiple patient classes with heterogeneous service rates, we also examine log-normally distributed service times. Finally, we investigate scenarios where arrival rates are subject to prediction errors.}

\jp{\textbf{Transfer policies.}
We compare the no-transfer policy to three transfer policies: Myopic, fluid, and API. The fluid policy is obtained by solving the fluid control problem \eqref{eq:fluid_obj}--\eqref{eq:non_negativity} over $M=7$ days. The Myopic policy is a truncated version of the fluid policy that transfers only up to $(\bar{\lambda_i^m} - \mu_i)^+$ at queue $i$, where $\bar{\lambda_i^m}$ is the average arrival rate on day $m$, thereby avoiding idle capacity (i.e., empty queues) for one day at a time. The API policy is computed by running Algorithm \ref{alg:API} using a logistic regression classifier with the features described in Section \ref{sec:ADP}, probability threshold $p=0.5$, $B=10$ simulation runs per state, and $j_{\mathrm{max}}=5$ iterations. All policies are implemented using a rolling-horizon approach \citep{powell2007approximate}, where decisions are re-computed at the start of each day with a 7-day planning horizon. For practical relevance, all policies are constrained to transfer at most three patients per day per hospital. For the API policy, this is enforced by including the constraint $|y_i - x_i| \leq 3$ for all $i$ in \eqref{eq:solve_target}, while for the fluid and Myopic policies, we impose $|(u[m]^\top-u[m])\mathrm{e}_i| \leq 3$ in \eqref{eq:fluid_obj}--\eqref{eq:non_negativity} for all $i$ and $m$, where $\mathrm{e}_i$ is a standard basis vector. Consistent with practice during the pandemic, we think of these policies as moving only COVID patients.}

\textbf{Calibration of simulation input parameters.}
\jp{We simulate a one-week horizon corresponding to shortly after a surge during the pandemic when the system is recovering from a large and imbalanced distribution of COVID patients.} The daily arrival rates are estimated based on the average number of arrivals for each day of the preceding four weeks of the horizon in our dataset, which includes internal transfers to the ICU from acute wards. The parameters of the service time distributions are estimated using the length-of-stay (LOS) data. We note that COVID patients' particularly long and variable LOS was a significant contributor to hospital congestion \citep{chanOptimizingInterHospitalPatient2023a}. The initial queue lengths are set to the difference between the occupancy at the beginning of the horizon and the number of beds at each ICU \jp{(or set to zero if occupancy is smaller)}. A summary of parameters is given in Table \ref{tab:parameters}.
\begin{table}[]
    {\color{black}\TABLE
    {\jp{Summary of the simulation inputs for the case study\label{tab:parameters}}}
    {\begin{tabular}{ccccc}
    \hline
          & \multicolumn{4}{c}{Queue} \\ \cline{2-5}
          & 1 & 2 & 3 & 4 \\ \hline
        Arrival rate & \multirow{2}{*}{$(2.7, 2.3)$} & \multirow{2}{*}{$(5.5, 4.6)$} & \multirow{2}{*}{$(5.0, 4.7)$} & \multirow{2}{*}{$(3.6, 3.1)$} \\ 
        (weekday, weekend) & \\ \hline
        Service time & \multirow{2}{*}{$(7.2, 13.7)$} & \multirow{2}{*}{$(5.7, 13.1)$} & \multirow{2}{*}{$(6.8, 10.7)$} & \multirow{2}{*}{$(6.1, 9.5)$} \\ 
        (mean, std.) & \\ \hline
        Capacity (beds) & 23 & 33 & 35 & 26 \\ \hline
        Initial condition & 18 & 43 & 45 & 21 \\ \hline
        Variable costs & \multirow{2}{*}{$(0, 0.5, 1.0, 0.7)$} & \multirow{2}{*}{$(0.5, 0, 0.8, 0.2)$} & \multirow{2}{*}{$(1.0, 0.8, 0, 0.7)$} & \multirow{2}{*}{$(0.7, 0.2, 0.7, 0)$} \\ 
        $(r_{i1},r_{i2},r_{i3},r_{i4})$ & \\ \hline
        Fixed (setup) costs & \multirow{2}{*}{$(0, 1, 1, 1)$} & \multirow{2}{*}{$(1, 0, 1, 1)$} & \multirow{2}{*}{$(1, 1, 0, 1)$} & \multirow{2}{*}{$(1, 1, 1, 0)$} \\ 
        $(K_{i1},K_{i2},K_{i3},K_{i4})$ & \\ \hline
    \end{tabular}}
    {}}
\end{table}

\textbf{Cost parameters.}
The unit variable transfer costs $r_{ij}$ reflect the distance between hospital $i$ and $j$ and are normalized to be in $[0,1]$ after dividing by the maximum distance among all hospitals. \jp{We use the general setup cost function $\tilde{\kappa}(u)=\sum_{i\in \mathcal{N}} \sum_{j\in \mathcal{N}} K_{ij} 1\{u_{ij} > 0\}$ with $K_{ij}=1$ for all $i,j$. Lastly, we use a linear holding cost function $h(x) = 2x$. In the Greater Toronto Area, historical transfers during the pandemic focused on equalizing the distribution of COVID patients across hospitals \citep{chanOptimizingInterHospitalPatient2023a}, suggesting a linear cost structure. Therefore, our results aim to compare various transfer policies under this directive in reducing the overall congestion.}

\textbf{Results and discussion.}
The results of the case study are summarized in Table \ref{tab:results}. \jp{All policies outperform the no-transfer policy, with estimated total cost reductions of 5.5\%, 15.7\%, and 27.7\% by the Myopic, fluid, and API policies, respectively. These improvements correspond to reductions of 8.8, 23.1, and 46.1 patient-days over ICU capacity, while requiring fewer than three transfers per day on average across the network. Notably, the fluid policy offers significant improvements over the Myopic policy, and the API policy achieves even greater benefits, reiterating the value of accounting for future costs and stochasticity. While the API policy may incur larger transfer costs due to its more proactive nature, these are more than offset by reductions in holding costs.}

\begin{figure}[]
    \FIGURE
    {\includegraphics[width=0.5\textwidth]{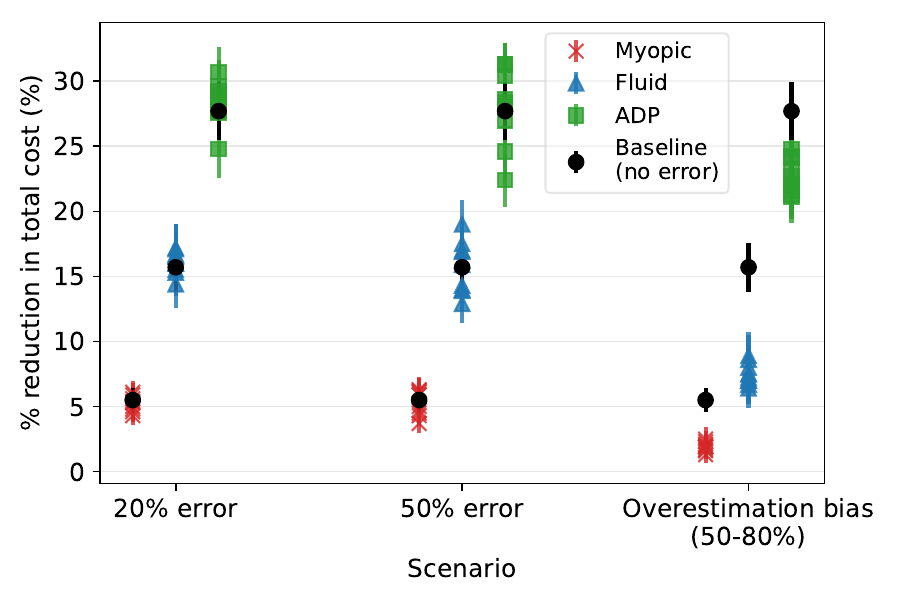}}
    {Performance under different arrival rate prediction error scenarios \label{fig:case_study_pred_err}}{\jp{Under each scenario and policy, the black dot represents the baseline performance from simulating the system with unperturbed arrival rates.}}
\end{figure}

\jp{We conduct two robustness checks. First, when service times are log-normally distributed with parameters from Table \ref{tab:parameters}, our findings remain consistent, as shown in Table \ref{tab:results_logn}. Second, we assess the impact of arrival rate prediction errors by computing policies under fixed arrival rates and then evaluating them in simulations using the ``true'' arrival rates, obtained by perturbing the original values under three different scenarios. In scenarios 1 and 2, we multiply the daily arrival rates by a random value uniformly drawn from $[0.8, 1.2]$ and $[0.5, 1.5]$, respectively, representing up to 20\% and 50\% prediction errors. In scenario 3, we draw uniformly from $[0.2, 0.5]$ to simulate an overestimation bias of 50--80\%. As shown in Figure \ref{fig:case_study_pred_err}, all three policies remain robust in scenarios 1 and 2. In scenario 3, note that the system is capable of clearing on its own relatively quickly, since the true arrival rates are much smaller. Consequently, transfers are generally less valuable, particularly later in the planning horizon, and the Myopic and fluid policies see large performance declines. The API policy remains relatively robust, due to its tendency to balance the system proactively and early in the planning horizon, when imbalances are still large.}

\begin{table}[]
    \TABLE
    {Summary of the simulation outputs for the case study using exponential service times \label{tab:results}}
    {\color{black} \begin{tabular}{lcccc}
    \hline
          & \multicolumn{4}{c}{Policy} \\ \cline{2-5}
          Performance measure & No-transfer & Myopic & Fluid & API \\ \hline
          Expected holding cost & 262.9 $\pm$ 16.3 & 245.3 $\pm$ 16.6 & 216.6 $\pm$ 16.6 & 170.5 $\pm$ 15.0 \\ 
          Patient days over capacity & 131.4 $\pm$ 8.2 & 122.6 $\pm$ 8.3 & 108.3 $\pm$ 8.3 & 85.3 $\pm$ 7.5 \\
          \% Reduction in holding cost & & 7.4 $\pm$ 1.2\% & 18.9 $\pm$ 2.1\% & 36.8 $\pm$ 2.2\% \\ \hline
          Expected transfer cost & & 4.7 $\pm$ 0.6 & 8.0 $\pm$ 0.7 & 22.2 $\pm$ 0.8 \\ 
          Avg. \# of transfers/week & & 3.1 & 8.5 & 18.3 \\ \hline
          Expected total cost & 262.9 $\pm$ 16.3 & 249.9 $\pm$ 16.5 & 224.6 $\pm$ 16.6 & 192.7 $\pm$ 15.4 \\
          \% Reduction in total cost & & 5.5 $\pm$ 0.9\% & 15.7 $\pm$ 1.9\% & 27.7 $\pm$ 2.2\% \\ \hline
    \end{tabular}}
    {\emph{Note.} The number after $\pm$ corresponds to the half-width of the 95\% confidence interval.}
\end{table}

\begin{table}[]
    \TABLE
    {Summary of the simulation outputs for the case study using log-normal service times \label{tab:results_logn}}
    {\color{black} \begin{tabular}{lcccc}
    \hline
          & \multicolumn{4}{c}{Policy} \\ \cline{2-5}
          Performance measure & No-transfer & Myopic & Fluid & API \\ \hline
          Expected holding cost & 259.2 $\pm$ 15.3 & 240.3 $\pm$ 15.4 & 209.1 $\pm$ 14.3 & 178.4 $\pm$ 13.9 \\ 
          Patient days over capacity & 129.6 $\pm$ 7.6 & 120.2 $\pm$ 7.7 & 104.6 $\pm$ 7.2 & 89.2 $\pm$ 7.0 \\
          \% Reduction in holding cost & & 7.8 $\pm$ 1.1\% & 19.7 $\pm$ 2.1\% & 32.0 $\pm$ 2.6\% \\ \hline
          Expected transfer cost & & 4.7 $\pm$ 0.5 & 7.6 $\pm$ 0.6 & 12.7 $\pm$ 0.6 \\ 
          Avg. \# of transfers/week & & 3.4 & 8.8 & 14.8 \\ \hline
          Expected total cost & 259.2 $\pm$ 15.3 & 245.1 $\pm$ 15.4 & 216.7 $\pm$ 14.4 & 191.0 $\pm$ 14.2 \\
          \% Reduction in total cost & & 5.9 $\pm$ 0.9\% & 16.6 $\pm$ 2.0\% & 26.9 $\pm$ 2.6\% \\ \hline
    \end{tabular}}
    {\emph{Note.} The number after $\pm$ corresponds to the half-width of the 95\% confidence interval.}
\end{table}

\section{Conclusion} \label{sec:conclusion}
In this work, we study the problem of transferring customers between parallel queues at discrete time intervals to balance transfer and congestion costs. We study an associated fluid control problem that allows us to obtain transfer policies under fairly general assumptions including time-varying arrivals and convex holding costs. Our analysis of the optimal fluid policy reveals several implications for managing imbalanced load for parallel queueing systems. By a careful trade-off among holding costs, transfer costs, and idleness between periods, we show that effective control policies are characterized by the so-called no-transfer region --- \jp{a compact, connected region of the state-space where no control is optimal --- thus formally establishing the optimality of the region-of-inaction policies}. When holding costs accrue linearly at the same rate at all queues, control is warranted if and only if there will be excessive idleness mid-period. Our results also highlight the impact of fixed costs on the structure of the optimal policy. In the presence of fixed costs, transfers should move the state to the relative interior of the no-transfer region, rather than the boundary. Therefore, frequent, small transfers are not cost effective in the presence of fixed costs. 

\jp{We leverage the structural results to design a simulation-based API algorithm for the original stochastic control problem. Our algorithm computes a region-of-inaction policy by approximating it directly with a classifier that labels each state as inside or outside the region and iteratively refines the classifications. We show that the structural property of the region is preserved when the parameters are set properly and propose a practical procedure to verify the structure. We demonstrate the effectiveness of our algorithm through simulation experiments and a case study based on real data from the COVID-19 pandemic in the Greater Toronto Area.}

\jp{Our model assumes a single class of customers in each queue. An interesting direction for future research is to extend the problem to multiple customer classes, where decisions involve determining the number of customers of each class to transfer. In this case, it may be optimal for each queue to be both a sender and a receiver. Key challenges include establishing the $K$-convexity of the value function and characterizing the state transition function under an appropriate priority discipline. In addition, we assume exogenous arrivals to each queue and focus on optimal transfer policies after customer arrivals. Joint optimization of  routing and transfer decisions would be an interesting direction for future research. In the context of our motivating application, this corresponds to jointly utilizing ambulance diversions and transfers to address imbalances in hospital occupancies.}

\bibliographystyle{informs2014} 
\bibliography{manuscript.bib} 

\begin{thebibliography}{58}
\providecommand{\natexlab}[1]{#1}
\providecommand{\url}[1]{\texttt{#1}}
\providecommand{\urlprefix}{URL }

\bibitem[{Armony(2005)}]{armony2005dynamic}
Armony M (2005) Dynamic routing in large-scale service systems with
  heterogeneous servers. \emph{Queueing Systems} 51:287--329.

\bibitem[{Armony \protect\BIBand{} Ward(2010)}]{armony2010fair}
Armony M, Ward AR (2010) Fair dynamic routing in large-scale
  heterogeneous-server systems. \emph{Operations Research} 58(3):624--637.

\bibitem[{Ata et~al.(2020)Ata, Barjesteh, \protect\BIBand{}
  Kumar}]{ataDynamicDispatchCentralized2020}
Ata B, Barjesteh N, Kumar S (2020) Dynamic {{Dispatch}} and {{Centralized
  Relocation}} of {{Cars}} in {{Ride-Hailing Platforms}}. \emph{Available at
  SSRN 3675888} .

\bibitem[{Avrachenkov et~al.(2015)Avrachenkov, Habachi, Piunovskiy,
  \protect\BIBand{} Zhang}]{avrachenkovInfiniteHorizonOptimal2015}
Avrachenkov K, Habachi O, Piunovskiy A, Zhang Y (2015) Infinite horizon optimal
  impulsive control with applications to {{Internet}} congestion control.
  \emph{International Journal of Control} 88(4):703--716.

\bibitem[{B{\"a}uerle(2000)}]{bauerle2000asymptotic}
B{\"a}uerle N (2000) Asymptotic optimality of tracking policies in stochastic
  networks. \emph{The Annals of Applied Probability} 10(4):1065--1083.

\bibitem[{Bellani et~al.(2021)Bellani, Grasselli, Cecconi
  et~al.}]{bellaniNoninvasiveVentilatorySupport2021}
Bellani G, Grasselli G, Cecconi M, et~al. (2021) Noninvasive {{Ventilatory
  Support}} of {{Patients}} with {{COVID-19}} outside the {{Intensive Care
  Units}} ({{WARd-COVID}}). \emph{Annals of the American Thoracic Society}
  18(6):1020--1026, ISSN 2329-6933,
  \urlprefix\url{http://dx.doi.org/10.1513/AnnalsATS.202008-1080OC}.

\bibitem[{Benjaafar et~al.(2022)Benjaafar, Jiang, Li, \protect\BIBand{}
  Li}]{benjaafarDynamicInventoryRepositioning2022}
Benjaafar S, Jiang D, Li X, Li X (2022) Dynamic {{Inventory Repositioning}} in
  {{On-Demand Rental Networks}}. \emph{Management Science} 68(11):7861--7878.

\bibitem[{Benkherouf \protect\BIBand{}
  Bensoussan(2009)}]{benkheroufOptimalityPolicyCompound2009}
Benkherouf L, Bensoussan A (2009) Optimality of an (s, {{S}}) {{Policy}} with
  {{Compound Poisson}} and {{Diffusion Demands}}: {{A Quasi-Variational
  Inequalities Approach}}. \emph{SIAM Journal on Control and Optimization}
  48(2):756--762.

\bibitem[{Bensoussan et~al.(2005)Bensoussan, Liu, \protect\BIBand{}
  Sethi}]{bensoussanOptimalityPolicyCompound2005}
Bensoussan A, Liu RH, Sethi SP (2005) Optimality of an (s, {{S}}) {{Policy}}
  with {{Compound Poisson}} and {{Diffusion Demands}}: {{A Quasi-variational
  Inequalities Approach}}. \emph{SIAM Journal on Control and Optimization}
  44(5):1650--1676.

\bibitem[{Berry~Jaeker \protect\BIBand{}
  Tucker(2017)}]{berryjaekerPointSpeedingNegative2017}
Berry~Jaeker JA, Tucker AL (2017) Past the {{Point}} of {{Speeding Up}}: {{The
  Negative Effects}} of {{Workload Saturation}} on {{Efficiency}} and {{Patient
  Severity}}. \emph{Management Science} 63(4):1042--1062.

\bibitem[{Bertsekas(2011)}]{bertsekasApproximatePolicyIteration2011}
Bertsekas DP (2011) Approximate policy iteration: A survey and some new
  methods. \emph{Journal of Control Theory and Applications} 9(3):310--335,
  ISSN 1993-0623.

\bibitem[{Cadenillas \protect\BIBand{}
  Zapatero(2000)}]{cadenillasClassicalImpulseStochastic2000}
Cadenillas A, Zapatero F (2000) Classical and {{Impulse Stochastic Control}} of
  the {{Exchange Rate Using Interest Rates}} and {{Reserves}}.
  \emph{Mathematical Finance} 10(2):141--156.

\bibitem[{{Caudillo-Fuentes} et~al.(2010){Caudillo-Fuentes}, Kaufman,
  \protect\BIBand{} Lewis}]{caudillo-fuentesSimpleHeuristicLoad2010a}
{Caudillo-Fuentes} LA, Kaufman DL, Lewis ME (2010) A simple heuristic for load
  balancing in parallel processing networks with highly variable service time
  distributions. \emph{Queueing Systems} 64(2):145--165.

\bibitem[{Chan et~al.(2021)Chan, Huang, \protect\BIBand{}
  Sarhangian}]{chanDynamicServerAssignment2021}
Chan CW, Huang M, Sarhangian V (2021) Dynamic {{Server Assignment}} in
  {{Multiclass Queues}} with {{Shifts}}, with {{Applications}} to {{Nurse
  Staffing}} in {{Emergency Departments}}. \emph{Operations Research}
  69(6):1936--1959.

\bibitem[{Chan et~al.(2023)Chan, Park, Pogacar, Sarhangian, Hellsten, Razak,
  \protect\BIBand{} Verma}]{chanOptimizingInterHospitalPatient2023a}
Chan T, Park J, Pogacar F, Sarhangian V, Hellsten E, Razak F, Verma A (2023)
  Optimizing {{Inter-Hospital Patient Transfer Decisions During}} a
  {{Pandemic}}: {{A Queueing Network Approach}}. \emph{Available at SSRN
  3975839} .

\bibitem[{Chen et~al.(2020)Chen, Dong, \protect\BIBand{}
  Shi}]{chenSurveySkillbasedRouting2020}
Chen J, Dong J, Shi P (2020) A survey on skill-based routing with applications
  to service operations management. \emph{Queueing Systems} 96(1):53--82.

\bibitem[{Chen et~al.(2023)Chen, Dong, \protect\BIBand{} Shi}]{chen2023optimal}
Chen J, Dong J, Shi P (2023) Optimal routing under demand surges: The value of
  future arrival rates. \emph{Operations Research} 0(0).

\bibitem[{Chen et~al.(2009)Chen, Huang, Kulkarni, Unnikrishnan, Zhu, Mehta,
  Meyn, \protect\BIBand{} Wierman}]{chenApproximateDynamicProgramming2009}
Chen W, Huang D, Kulkarni AA, Unnikrishnan J, Zhu Q, Mehta P, Meyn S, Wierman A
  (2009) Approximate dynamic programming using fluid and diffusion
  approximations with applications to power management. \emph{Proceedings of
  the 48h {{IEEE Conference}} on {{Decision}} and {{Control}} ({{CDC}}) Held
  Jointly with 2009 28th {{Chinese Control Conference}}}, 3575--3580, ISSN
  0191-2216.

\bibitem[{Cini et~al.(2023)Cini, Neto, Burrell, Udy, \protect\BIBand{}
  Investigators}]{ciniInterhospitalTransferClinical2023}
Cini C, Neto AS, Burrell A, Udy A, Investigators tSSA (2023) Inter-hospital
  transfer and clinical outcomes for people with {{COVID-19}} admitted to
  intensive care units in {{Australia}}: An observational cohort study.
  \emph{Medical Journal of Australia} 218(10):474--481, ISSN 1326-5377.

\bibitem[{Dai \protect\BIBand{} Gluzman(2022)}]{dai2022queueing}
Dai JG, Gluzman M (2022) Queueing network controls via deep reinforcement
  learning. \emph{Stochastic Systems} 12(1):30--67.

\bibitem[{Dai \protect\BIBand{}
  Shi(2019)}]{daiInpatientOverflowApproximate2019a}
Dai JG, Shi P (2019) Inpatient {{Overflow}}: {{An Approximate Dynamic
  Programming Approach}}. \emph{Manufacturing \& Service Operations Management}
  21(4):894--911, ISSN 1523-4614.

\bibitem[{Dai \protect\BIBand{}
  Yao(2013{\natexlab{a}})}]{daiBrownianInventoryModels2013}
Dai JG, Yao D (2013{\natexlab{a}}) Brownian {{Inventory Models}} with {{Convex
  Holding Cost}}, {{Part}} 1: {{Average-Optimal Controls}}. \emph{Stochastic
  Systems} 3(2):442--499.

\bibitem[{Dai \protect\BIBand{}
  Yao(2013{\natexlab{b}})}]{daiBrownianInventoryModels2013a}
Dai JG, Yao D (2013{\natexlab{b}}) Brownian {{Inventory Models}} with {{Convex
  Holding Cost}}, {{Part}} 2: {{Discount-Optimal Controls}}. \emph{Stochastic
  Systems} 3(2):500--573.

\bibitem[{Dijkstra et~al.(2023)Dijkstra, Baas, Braaksma, \protect\BIBand{}
  Boucherie}]{dijkstraDynamicFairBalancing2023}
Dijkstra S, Baas S, Braaksma A, Boucherie RJ (2023) Dynamic fair balancing of
  {{COVID-19}} patients over hospitals based on forecasts of bed occupancy.
  \emph{Omega} 116:102801, ISSN 0305-0483.

\bibitem[{Dolan et~al.(2022)Dolan, Johnson, Kepler, Lam
  et~al.}]{dolanHospitalLoadBalancing2022}
Dolan E, Johnson N, Kepler T, Lam H, et~al. (2022) Hospital {{Load Balancing}}:
  {{A Data-Driven Approach}} to {{Optimize Ambulance Transports During}} the
  {{COVID-19 Pandemic}} in {{New York City}}.

\bibitem[{Douin et~al.(2021)Douin, Ward, Lindsell, Howell, Hough, Exline, Gong,
  Aboodi, Tenforde, Feldstein, Stubblefield, Steingrub, Prekker, Brown, Peltan,
  Khan, Files, Gibbs, Rice, Casey, Hager, Qadir, Henning, Wilson, Patel, Self,
  \protect\BIBand{} Ginde}]{douinICUBedUtilization2021}
Douin DJ, Ward MJ, Lindsell CJ, Howell MP, Hough CL, Exline MC, Gong MN, Aboodi
  MS, Tenforde MW, Feldstein LR, Stubblefield WB, Steingrub JS, Prekker ME,
  Brown SM, Peltan ID, Khan A, Files DC, Gibbs KW, Rice TW, Casey JD, Hager DN,
  Qadir N, Henning DJ, Wilson JG, Patel MM, Self WH, Ginde AA (2021) {{ICU Bed
  Utilization During}} the {{Coronavirus Disease}} 2019 {{Pandemic}} in a
  {{Multistate Analysis}}\textemdash{{March}} to {{June}} 2020. \emph{Critical
  Care Explorations} 3(3):e0361.

\bibitem[{Down \protect\BIBand{} Lewis(2006)}]{downDynamicLoadBalancing2006}
Down DG, Lewis ME (2006) Dynamic load balancing in parallel queueing systems:
  {{Stability}} and optimal control. \emph{European Journal of Operational
  Research} 168(2):509--519.

\bibitem[{Gallego \protect\BIBand{} Sethi(2005)}]{gallegoKConvexityRn2005a}
Gallego G, Sethi S (2005) K-{{Convexity}} in {{Rn}}. \emph{Journal of
  Optimization Theory and Applications} 127.

\bibitem[{Harrison(2003)}]{harrison2003broader}
Harrison JM (2003) A broader view of brownian networks. \emph{The Annals of
  Applied Probability} 13(3):1119--1150.

\bibitem[{He et~al.(2020)He, Hu, \protect\BIBand{}
  Zhang}]{heRobustRepositioningVehicle2020}
He L, Hu Z, Zhang M (2020) Robust {{Repositioning}} for {{Vehicle Sharing}}.
  \emph{Manufacturing \& Service Operations Management} 22(2):241--256.

\bibitem[{He \protect\BIBand{} Neuts(2002)}]{heTwoQueuesTransfers2002}
He QM, Neuts MF (2002) Two {{M}}/{{M}}/1 {{Queues}} with {{Transfers}} of
  {{Customers}}. \emph{Queueing Systems} 42(4):377--400.

\bibitem[{Henry et~al.(2024)Henry, Funsten, Michealson
  et~al.}]{henryInterfacilityPatientTransfers2024}
Henry MB, Funsten E, Michealson MA, et~al. (2024) Interfacility {{Patient
  Transfers During COVID-19 Pandemic}}: {{Mixed-Methods Study}}. \emph{Western
  Journal of Emergency Medicine} 25(5):758--766, ISSN 1936-900X.

\bibitem[{Hu et~al.(2022)Hu, Chan, \protect\BIBand{}
  Dong}]{huOptimalSchedulingProactive2022}
Hu Y, Chan CW, Dong J (2022) Optimal {{Scheduling}} of {{Proactive Service}}
  with {{Customer Deterioration}} and {{Improvement}}. \emph{Management
  Science} 68(4):2533--2578.

\bibitem[{Korn(1999)}]{kornApplicationsImpulseControl1999}
Korn R (1999) Some applications of impulse control in mathematical finance.
  \emph{Mathematical Methods of Operations Research} 50(3):493--518.

\bibitem[{Kumar \protect\BIBand{} Kumar(2019)}]{kumarIssuesChallengesLoad2019}
Kumar P, Kumar R (2019) Issues and {{Challenges}} of {{Load Balancing
  Techniques}} in {{Cloud Computing}}: {{A Survey}}. \emph{ACM Computing
  Surveys} 51(6):120:1--120:35, ISSN 0360-0300.

\bibitem[{Kuntz et~al.(2015)Kuntz, Mennicken, \protect\BIBand{}
  Scholtes}]{kuntzStressWardEvidence2015}
Kuntz L, Mennicken R, Scholtes S (2015) Stress on the {{Ward}}: {{Evidence}} of
  {{Safety Tipping Points}} in {{Hospitals}}. \emph{Management Science}
  61(4):754--771.

\bibitem[{Luo et~al.(2015)Luo, Rao, \protect\BIBand{}
  Liu}]{luoSpatioTemporalLoadBalancing2015}
Luo J, Rao L, Liu X (2015) Spatio-{{Temporal Load Balancing}} for {{Energy Cost
  Optimization}} in {{Distributed Internet Data Centers}}. \emph{IEEE
  Transactions on Cloud Computing} 3(3):387--397, ISSN 2168-7161.

\bibitem[{Maglaras(2000)}]{maglaras2000discrete}
Maglaras C (2000) Discrete-review policies for scheduling stochastic networks:
  Trajectory tracking and fluid-scale asymptotic optimality. \emph{The Annals
  of Applied Probability} 10(3):897--929.

\bibitem[{Mandelbaum \protect\BIBand{} Massey(1995)}]{mandelbaum1995strong}
Mandelbaum A, Massey WA (1995) Strong approximations for time-dependent queues.
  \emph{Mathematics of Operations Research} 20(1):33--64.

\bibitem[{Maxwell et~al.(2013)Maxwell, Henderson, \protect\BIBand{}
  Topaloglu}]{maxwellTuningApproximateDynamic2013}
Maxwell MS, Henderson SG, Topaloglu H (2013) Tuning {{Approximate Dynamic
  Programming Policies}} for {{Ambulance Redeployment}} via {{Direct Search}}.
  \emph{Stochastic Systems} 3(2):322--361, ISSN 1946-5238.

\bibitem[{Meyn(1997)}]{meyn1997stability}
Meyn S (1997) Stability and optimization of queueing networks and their fluid
  models. \emph{Lectures in applied mathematics-American Mathematical Society}
  33:175--200.

\bibitem[{Meyn(2008)}]{meyn2008control}
Meyn S (2008) \emph{Control techniques for complex networks} (Cambridge
  University Press).

\bibitem[{Mitchell et~al.(2014)Mitchell, Feng, \protect\BIBand{}
  Muthuraman}]{mitchell2014impulse}
Mitchell D, Feng H, Muthuraman K (2014) Impulse control of interest rates.
  \emph{Operations research} 62(3):602--615.

\bibitem[{Moallemi et~al.(2008)Moallemi, Kumar, \protect\BIBand{}
  Van~Roy}]{moallemiApproximateDataDrivenDynamic}
Moallemi CC, Kumar S, Van~Roy B (2008) Approximate and {{Data-Driven Dynamic
  Programming}} for {{Queueing Networks}} .

\bibitem[{Ong et~al.(2025)Ong, Lim, Zahrin
  et~al.}]{ongAmbulanceDiversionIts2025}
Ong JHM, Lim BJW, Zahrin MA, et~al. (2025) Ambulance diversion and its use as
  an {{ED}} overcrowding mitigation strategy: {{Does}} it work? {{A}} scoping
  review. \emph{International Journal of Emergency Medicine} 18(1):125, ISSN
  1865-1380.

\bibitem[{Ormeci et~al.(2008)Ormeci, Dai, \protect\BIBand{}
  Vate}]{ormeci2008impulse}
Ormeci M, Dai JG, Vate JV (2008) Impulse control of brownian motion: The
  constrained average cost case. \emph{Operations Research} 56(3):618--629.

\bibitem[{Powell(2007)}]{powell2007approximate}
Powell WB (2007) \emph{Approximate Dynamic Programming: Solving the curses of
  dimensionality}, volume 703 (John Wiley \& Sons).

\bibitem[{Rockafellar(2015)}]{rockafellarConvexAnalysis2015}
Rockafellar RT (2015) Convex {{Analysis}}. \emph{Convex {{Analysis}}}
  ({Princeton University Press}), ISBN 978-1-4008-7317-3.

\bibitem[{Scarf(1960)}]{scarfOptimalityPoliciesDynamic1960a}
Scarf H (1960) The {{Optimality}} of (s,{{S}}) {{Policies}} in the {{Dynamic
  Inventory Problem}}. \emph{Mathematical Methods in the Social Sciences} .

\bibitem[{Sethi \protect\BIBand{} Thompson(2000)}]{sethi2000economic}
Sethi SP, Thompson GL (2000) Economic applications. \emph{Optimal Control
  Theory: Applications to Management Science and Economics} 289--306.

\bibitem[{Sitaraman(2001)}]{sitaraman2001power}
Sitaraman R (2001) The power of two random choices: A survey of techniques and
  results .

\bibitem[{Sun et~al.(2024)Sun, Dai, \protect\BIBand{} Shi}]{sun2024inpatient}
Sun J, Dai J, Shi P (2024) Inpatient overflow management with proximal policy
  optimization. \emph{arXiv preprint arXiv:2410.13767} .

\bibitem[{Sun \protect\BIBand{} Zhu(2025)}]{sunDynamicControlMaketoOrder2025}
Sun X, Zhu X (2025) Dynamic {{Control}} of a {{Make-to-Order System Under Model
  Uncertainty}}. \emph{Management Science} ISSN 0025-1909.

\bibitem[{Tien et~al.(2020)Tien, Sawadsky, Lewell, Peddle, \protect\BIBand{}
  Durham}]{tienCriticalCareTransport2020}
Tien H, Sawadsky B, Lewell M, Peddle M, Durham W (2020) Critical care transport
  in the time of {{COVID-19}}. \emph{Canadian Journal of Emergency Medicine}
  22(S2):S84--S88, ISSN 1481-8035, 1481-8043.

\bibitem[{{Van der Boor} et~al.(2022){Van der Boor}, Borst, Van~Leeuwaarden,
  \protect\BIBand{} Mukherjee}]{vanderboorScalableLoadBalancing2022}
{Van der Boor} M, Borst SC, Van~Leeuwaarden JSH, Mukherjee D (2022) Scalable
  {{Load Balancing}} in {{Networked Systems}}: {{A Survey}} of {{Recent
  Advances}}. \emph{SIAM Review} 64(3):554--622.

\bibitem[{Zeng et~al.(2018)Zeng, Zhang, Cai, \protect\BIBand{}
  Li}]{zengCostSharingCapacity2018}
Zeng Y, Zhang L, Cai X, Li J (2018) Cost {{Sharing}} for {{Capacity Transfer}}
  in {{Cooperating Queueing Systems}}. \emph{Production and Operations
  Management} 27(4):644--662.

\bibitem[{Zychlinski(2023)}]{zychlinski2023applications}
Zychlinski N (2023) Applications of fluid models in service operations
  management. \emph{Queueing Systems} 103(1):161--185.

\bibitem[{Zychlinski et~al.(2023)Zychlinski, Chan, \protect\BIBand{}
  Dong}]{zychlinski2023managing}
Zychlinski N, Chan CW, Dong J (2023) Managing queues with different resource
  requirements. \emph{Operations Research} 71(4):1387--1413.

\end{thebibliography}

%
%
%

\newpage
\begin{APPENDICES}
\section{Numerical Solution Approach for the Fluid Control Problem} \label{appen:numerical_solution_approach}
\jp{When computing the optimal fluid policy, it is often more convenient to formulate it as a mathematical program. To this end, let $g^m(x,u)$ denote the total single-period cost in period $m$, starting from (pre-transfer) state $x$, and under transfer decision $u$. We have
\begin{align}
    g^m(x,u) = H^m\left(x+ \left(u^\top - u\right)\mathrm{e} \right)  + r \cdot u + \tilde{\kappa}(u).
\end{align}
Then the fluid control problem can be written as follows. Starting with a given initial condition $x^0$, the objective is to find a sequence of control matrices $\{u[m]; m\in\mathcal{M}\}$ to minimize the total cost over the horizon:
\begin{align}
    \min \quad & \sum_{m=0}^{M-1} g^m(x[m],u[m]) \label{eq:fluid_obj}\\
    \mbox{s.t.} \quad & x[m+1] = f^m\left(x[m] + \left(u[m]^\top-u[m] \right)\mathrm{e}, (m+1)\tau \right), \label{eq:dynamics} \quad \forall m \in \mathcal{M}, \\
    & x[0] = x^0, \label{eq:initial_condition}\\
    & u[m] \mathrm{e} \leq x[m],\quad \forall m\in \mathcal{M}, \label{eq:upper_bound_on_u} \\
    & u_{ij}[m] \geq 0, \quad \forall i \in \mathcal{N},j \in \mathcal{N}, m\in \mathcal{M}. \label{eq:non_negativity}
\end{align}
Eq. \eqref{eq:dynamics} relates the state of system at the beginning of the next period to  state and transfer decision in the current period. Eq. \eqref{eq:initial_condition} ensures that we start from the given initial condition. Eqs. \eqref{eq:upper_bound_on_u} and \eqref{eq:non_negativity} parallel the admissibility conditions of the control policy for the stochastic system.}

\jp{For a given time-varying arrival rate function and a convex holding cost function, the fluid control problem \eqref{eq:fluid_obj}--\eqref{eq:non_negativity} can be solved numerically by the following general framework.} We first approximate the continuous fluid dynamics within each period using $L$ discrete intervals of fixed width. We define by
\begin{align*}
    \bar{\lambda}_{lm} = \frac{L}{\tau}\int_{m\tau + l(\tau/L)}^{m\tau + (l+1)(\tau/L)}\lambda(t) dt
\end{align*}
the average arrival rate over interval $l, l \in \{0, 1, \ldots, L-1\}$, within period $m$. Within each interval, the fluid state $x(t)$ is assumed constant, and from interval $l$ to $l+1$ of period $m$, it changes by $(\bar{\lambda}_{lm}-\mu)(\tau/L)$. Given an initial condition $x[m]$ and control $u[m]$ in period $m$, we use $g_i(x[m],u[m])$ to denote the single-period cost at queue $i$ such that $g(x[m],u[m]) = \sum_{i \in \mathcal{N}}g_i(x[m],u[m])$. Then,
\begin{equation*}
    g_i(x[m],u[m]) \approx \frac{\tau}{L}\sum_{l\in\mathcal{L}\setminus \{0\}} h_i\left(\frac{z_{il}[m] + z_{i,l-1}[m]}{2} \right) + \sum_{j \in \mathcal{N}}r_{ij}u_{ij}[m] + \tilde{\kappa}(u[m]),
\end{equation*}
where $\mathcal{L}=\{0, \ldots, L\}$, \jp{$z_{il}[m]$ represents the fluid state at the start of interval $l$ within period $m$,} and $h_i:\mathbb{R}_+ \xrightarrow{} \mathbb{R}_+$ denotes a generic convex function, which can be approximated by the pointwise maximum of $J$ affine functions
\begin{align} \label{eq:holding_cost_approx}
    h_i(a) \approx \max\{h_{i1}a + b_{i1}, \ldots, h_{iJ}a + b_{iJ}\},
\end{align}
which can be linearized using auxiliary variables $w_{il} \in \mathbb{R}_+$ after imposing the constraints $w_{il} \geq h_{ij}a + b_{ij}$ for all $j$ for each $i \in \mathcal{N}$ and $l \in \mathcal{L} \setminus \{0\}$. The goodness of the approximation improves with larger $L$ and $J$. The following constraints replace \eqref{eq:dynamics}-\eqref{eq:upper_bound_on_u} to approximate the fluid dynamics:
\begin{align}
    & y_{il}[m] = y_{i,l-1}[m] + (\bar{\lambda}_{lm} - \mu_i)(\tau/L), \quad \forall  i\in\mathcal{N}, l\in\mathcal{L}\setminus \{0\}, m\in\mathcal{M}, \label{eq:interval_dynamics}\\
    & z_{il}[m] \geq y_{il}[m], \quad \forall i\in\mathcal{N}, l\in\mathcal{L}, m\in\mathcal{M}, m=l\neq0, \label{eq:nonnegative_y1} \\
    & z_{il}[m] \geq 0, \quad \forall i\in\mathcal{N}, l\in\mathcal{L}, m\in\mathcal{M}, \label{eq:nonnegative_y2}\\
    & y_{i0}[m] = z_{iL}[m-1] + \sum_{j\in\mathcal{N}}u_{ji}[m] - \sum_{j\in\mathcal{N}}u_{ij}[m], \quad \forall  i\in\mathcal{N}, m\in\mathcal{M}\setminus \{0\}, \label{eq:post_transfer_y1}\\
    & y_{i0}[0] = z_{i0}[0] + \sum_{j\in\mathcal{N}}u_{ji}[0] - \sum_{j\in\mathcal{N}}u_{ij}[0], \quad \forall i\in\mathcal{N}, \label{eq:post_transfer_y2}\\
    & z_{i0}[0] \geq x^0[0], \quad \forall i\in\mathcal{N}, \label{eq:initial_condition_z}\\
    & \sum_{j \in \mathcal{N}}u_{ij}[m] \leq z_{iL}[m-1], \quad \forall  i \in \mathcal{N}, m\in\mathcal{M}\setminus \{0\}, \label{eq:feasibility_z1}\\
    & \sum_{j \in \mathcal{N}}u_{ij}[0] \leq z_{i0}[0], \quad \forall i \in \mathcal{N}. \label{eq:feasibility_z2}
\end{align}
The variables $y_{il}[m]$ represent the fluid state of queue $i$ at each of the $L$ intervals within a period following the piecewise-constant dynamics, which is enforced in equation \eqref{eq:interval_dynamics}. In equations \eqref{eq:nonnegative_y1} and \eqref{eq:nonnegative_y2}, variables $z_{il}[m]$ take the non-negative part of $y_{il}[m]$ to ensure feasibility. At the beginning of each period, $y_{il}[m]$ is set to the fluid state just after transferring, which is specified through equations \eqref{eq:post_transfer_y1} and \eqref{eq:post_transfer_y2}. Equation \eqref{eq:initial_condition_z} is the initial condition. Finally, we dictate in equations \eqref{eq:feasibility_z1} and \eqref{eq:feasibility_z2} that the total transfers out of any queue is always bounded above by its state just prior to transferring. Using integer variables to model the setup cost function $\tilde{\kappa}(u)$, the resulting optimization problem is a mixed-integer linear program.


\section{Proofs}

\subsection{Proof of Lemma \ref{lem:joint_setup_cost_properties}: Properties of the Joint Setup Cost Function}
\begin{proof}{Proof of Lemma \ref{lem:joint_setup_cost_properties}.}
    The first two properties of the joint setup cost function follow directly from the fact that the indicator function $1\{z\neq0 \}$ satisfies $1\{x+y \neq 0\} \leq 1\{x\neq0\} + 1\{y\neq0\}$ for any $x,y\in \mathbb{R}^N$ and $1\{-z \neq0\} = 1\{z \neq0\}$ 
    for any $z\in\mathbb{R}^N$. For the third property, we note that if $x = y$, the joint setup cost is zero, but equals $K$ otherwise. Thus, for any $x$ and $y$, the joint setup cost can be calculated without the knowledge of the particular transfer matrix in moving the state from $x$ to $y$. Moreover, for all feasible non-zero transfer decision matrices, the joint setup cost is constant. This implies that the joint setup cost cannot affect the optimal transfer decision matrix and vice versa, i.e., $C(y-x) = R(y-x) + \kappa(y-x)$.
\end{proof}

\subsection{Properties of the Value Function} \label{appen:V_properties}
In this section, we provide proofs of the properties of the holding cost function (Lemma \ref{lem:properties_H}) and the value function (Theorem \ref{thm:properties_V}). Upon establishing $K$-convexity of the value function, we conclude by outlining its additional properties which are important in characterizing the optimal policy.

\subsubsection{Proof of Lemma \ref{lem:properties_H}: Properties of the Holding Cost Function} 
\begin{proof}{Proof of Lemma \ref{lem:properties_H}.}
    We show the properties by proving that the state transition function $f^m(\cdot, \tau)$ is convex, continuous, and non-decreasing, which is done by first deriving an equivalent recursive expression for it. From this, the properties of the holding cost function follow. 
    
    Since the holding cost can be analyzed separately by each queue and period, in what follows we will focus on a given queue $i$ and the first period and suppress the dependency of the holding cost function on the period for ease of exposition. Recall that the holding cost function at queue $i$ is defined as 
    \begin{align*} 
        H_i(y_i) = \int_{0}^{\tau} h_i\left(f_i(y_i,s) \right)ds.
    \end{align*}
     We proceed by a recursive expression for $f_i(y_i,t)$. Let $\mathcal{T}_i$ be the union of the set of time points in $(0,\tau)$ such that $\lambda_i(t)$ is monotone between successive points, and the set of all zeros of the fluid dynamics $\lambda_i(t)-\mu_i$. \jp{If $\lambda_i(t) = \mu_i$ over any interval, then we include only the time point marking the end of this interval in $\mathcal{T}_i$.} Then under Assumption \ref{ass:time_varying_arrivals}, $\mathcal{T}_i$ is finite. Denote these points by $t_1, \ldots, t_P$, arranged in increasing order. Additionally, let $t_0=0$ and $t_{P+1}=\tau$. We define $f^p_i:\mathbb{R}_+^2 \xrightarrow{} \mathbb{R}_+$ as
    \begin{align} \label{eq:transition_function_time_varying}
        f^{p+1}_i(y_i,t) = \left( f^{p}_i(y_i,t_p) + \int_{t_p}^{t}\lambda_i(s)ds - \mu_i(t-t_p) \right)^+, \quad t\in[t_p,t_{p+1}) \text{ and } p=0,\ldots,P,
    \end{align}
    with $f^0_i(y_i,t_0) = y_i$. For all $t \in [t_p, t_{p+1}]$, the expression $\lambda_i(t)-\mu_i$ must be either \jp{non-negative or non-positive}. Consequently, the queue length process is monotone in each interval $[t_p, t_{p+1}]$, and as a result, if the queue length reaches zero at any point in $[t_p, t_{p+1}]$, it will remain at zero until $t_{p+1}$. Therefore, the pointwise maximum operator  in \eqref{eq:transition_function_time_varying} correctly calculates the queue length throughout each interval and we arrive at the following equivalent expression for $f_i(y_i,t)$:
    \begin{align} \label{eq:recursive_f}
        f_i(y_i,t) = 
        \begin{cases}
            f^1_i(y_i,t), & 0=t_0 \leq t < t_1, \\
            f^2_i(y_i,t), & t_1 \leq t < t_2, \\
            & \vdots \\
            f^{P+1}_i(y_i,t), & t_P \leq t < t_{P+1}=\tau. \\
        \end{cases}
    \end{align}
    Crucially, each $f^p_i(y_i,t)$ is convex, continuous, and non-decreasing in $y_i$ as the composition of such functions preserves these properties. To calculate the holding cost at queue $i$, we can decompose it as a sum of the holding cost over each interval as follows:
    \begin{align*}
        H_i(y_i) 
        &= \int_0^\tau \Big[ h_i(f^0_i(y_i,s))1\{0\leq s < t_1\} + \cdots + h_i(f^{P+1}_i(y_i,s))1\{t_{P} \leq s < \tau \} \Big]ds \\
        &= \int_0^{t_1} h_i(f^0_i(y_i,s))ds + \cdots + \int_{t_{P}}^{\tau}h_i(f^{P+1}_i(y_i,s))ds \\
        &= \sum_{p=0}^{P} \int_{t_p}^{t_{p+1}}h_i(f^p_i(y_i,s))ds.
    \end{align*}
    Under Assumption \ref{ass:non_linear_holding_cost}, $H_i(\cdot)$ is convex, continuous, and non-decreasing as it is a sum of $P+1$ such functions. Thus, the holding cost of the system $H(y) = \sum_{i \in \mathcal{N}}H_i(y_i)$ is convex, continuous, and non-decreasing.
\end{proof}

\subsubsection{Proof of Theorem \ref{thm:properties_V}.}
Before proving the properties of the value function, we first establish properties of the transfer cost $R(\cdot)$ in \eqref{eq:variable_cost}.
\begin{lemma}\label{lem:R_properties}
Let $\mathcal{Z} = \{z\in \mathbb{R}^N: \mathrm{e}^\top z = 0\}$. The transfer cost function $R(\cdot)$ has the following properties:
\begin{itemize}
    \item (Positive homogeneity): $R(tz) = tR(z)$ for all $z \in \mathcal{Z}$ and $t\geq0$.
    \item (Convexity): $R(\theta z_1 + (1-\theta)z_2) \leq \theta R(z_1) + (1-\theta)R(z_2)$ for all $z_1, z_2 \in \mathcal{Z}$ and $\theta \in [0,1]$.
    \item (Subadditivity): $R(z_1+z_2) \leq R(z_1) + R(z_2)$ for all $z_1, z_2 \in \mathcal{Z}$.
    \item (Continuity): $R(z)$ is continuous in $z \in \mathcal{Z}$.
\end{itemize}
\end{lemma}
\begin{proof}{Proof of Lemma \ref{lem:R_properties}.}
    We note that the $R(z)$ is a bounded and feasible linear program for any $z \in \mathcal{Z}$. Therefore, by strong duality, we can write $R(z) = \max\{p \cdot z: p_j - p_i \leq r_{ij}, \forall i,j \in \mathcal{N} \}$. The rest of the proof follows the same approach from Lemma EC.1 of \cite{benjaafarDynamicInventoryRepositioning2022}. Observe that,
    \begin{equation*}
        R(tz) = \max\{t p \cdot z: p_j - p_i \leq r_{ij}, \forall i,j \in \mathcal{N} \} = t \max\{p \cdot z: p_j - p_i \leq r_{ij}, \forall i,j \in \mathcal{N} \} = tR(z),
    \end{equation*}
    for all $t\geq 0$, i.e., $R(\cdot)$ is positively homogeneous. As $R(\cdot)$ is a pointwise supremum of affine and continuous functions, it is convex and lower semicontinuous. Moreover, $\mathcal{Z}$ is a (\jp{convex}) polyhedron, and hence, a locally simplicial \jp{convex} set. \jp{Since a convex function on a locally simplicial convex set is upper semicontinuous \citep[Theorem 10.2]{rockafellarConvexAnalysis2015}, $R(\cdot)$ is upper semicontinuous, and in turn, continuous}. Finally, using convexity and positive homogeneity, we have
    \begin{equation*}
        R(z_1 + z_2) = 2R\left( \frac{1}{2}z_1 + \frac{1}{2}z_2 \right) \leq 2\left(\frac{1}{2} R(z_1) + \frac{1}{2}R(z_2) \right) = R(z_1) + R(z_2),
    \end{equation*}
    i.e., $R(\cdot)$ is sub-additive.
\end{proof}

\begin{proof}{Proof of Theorem \ref{thm:properties_V}.} Since showing $K$-convexity requires monotonicity, we will first establish monotonicity and then $K$-convexity by induction and lastly continuity. 

    To show monotonicity, consider two initial conditions $y$ and $z$ such that $y\geq z$. We will use the following equivalent representation for the value function:
    \begin{equation*}
        V^{M-1}(x) = \min_{u \in \mathcal{U}(x)}\left[ \jp{H^{M-1}}(x+(u^\top - u)\mathrm{e}) + R((u^\top - u)\mathrm{e}) + \kappa((u^\top - u)\mathrm{e}) \right],
    \end{equation*}
    where $\mathcal{U}(x) = \{u \in \mathbb{R}^{N \times N}_+: u \mathrm{e} \leq x, \forall i \in \mathcal{N}\}$. For convenience, we will use $\phi(u)$ to denote the net transfer $(u^\top - u)\mathrm{e}$ and $\phi_i(u)$ to denote its $i$th component. Intuitively, we now minimize over all feasible transfer decision \emph{matrices}, rather than post-transfer states. Denote by $u^*$ the optimal transfer matrix at $y$. We will construct a feasible transfer matrix $\hat{u}$ at $z$ from $u^*$ in the following way. For each $i \in \mathcal{N}$:
    \begin{itemize}
        \item If $\sum_{j \in \mathcal{N}}u^*_{ij} \leq z_i$, do nothing. This solution is also feasible at $z$.
        \item Otherwise, this means $u^*$ is not feasible at $z$ since more customers are transferred out of queue $i$ than are available. Choose any number of values from $u^*_{i1}, \ldots, u^*_{iN}$ and reduce by some arbitrary amount such that we will ultimately have $\sum_{j \in \mathcal{N}}\hat{u}_{ij} = z_i$ while maintaining $\hat{u}_{ij} \geq 0, \forall j \in \mathcal{N}$.
    \end{itemize}
    We then have,
    \begin{align*}
        V^{M-1}(y) &= \jp{H^{M-1}}(y+ \phi(u^*)) + R(\phi(u^*)) + \kappa(\phi(u^*))\\
        &\geq \jp{H^{M-1}}(z+\phi(\hat{u})) + R(\phi(\hat{u})) + \kappa(\phi(\hat{u}))\\
        &\geq \min_{u \in \mathcal{U}(z)}\left[ \jp{H^{M-1}}(z+\phi(u)) + R(\phi(u)) + \kappa(\phi(u)) \right] = V^{M-1}(z).
    \end{align*}
    The second inequality holds because $\hat{u}$ is feasible at $z$, but not optimal in general. To see that the first inequality holds, we note that by construction, we have $u^*_{ij} \geq \hat{u}_{ij}, \forall i,j \in \mathcal{N}$, which implies $\kappa(\phi(u^*)) \geq \kappa(\phi(\hat{u}))$ and $R(\phi(u^*)) = r \cdot u^* \geq r \cdot \hat{u} = R(\phi(\hat{u}))$. Next, to show $\jp{H^{M-1}}(y+\phi(u^*)) \geq \jp{H^{M-1}}(z+\phi(\hat{u}))$, we note that $\jp{H^{M-1}}(x) = \sum_{i \in \mathcal{N}}\jp{H^{M-1}_i}(x_i)$, so it is sufficient to show this for a fixed $i \in \mathcal{N}$. To this end, we consider two cases, where queue $i$ is a sender or a receiver.
    
    \textbf{Case 1. } Queue $i$ is a ``sender.'' If $\phi_i(u^*) = \phi_i(\hat{u})$, then clearly $y_i + \phi_i(u^*) \geq z_i + \phi_i(\hat{u})$ since $y_i \geq z_i$. Otherwise, we have $z_i + \phi_i(\hat{u}) = 0$ by construction, implying \jp{$y_i + \phi_i(u^*) \geq z_i + \phi_i(\hat{u})$}. Thus, either way, we have $\jp{H^{M-1}_i}(y_i + \phi_i(u^*)) \geq \jp{H^{M-1}_i}(z_i+\phi_i(\hat{u}))$ by monotonicity of $\jp{H^{M-1}}(\cdot)$.
    
    \textbf{Case 2. } Queue $i$ is a ``receiver.'' Since $u^*_{ji} \geq \hat{u}_{ji}$ for all $j$, queue $i$ receives fewer customers in total under $\hat{u}$. Therefore, $y_i + \phi_i(u^*) \geq z_i + \phi_i(\hat{u})$ and $\jp{H^{M-1}_i}(y_i + \phi_i(u^*)) \geq \jp{H^{M-1}_i}(z_i+\phi_i(\hat{u}))$ by monotonicity of $\jp{H^{M-1}}(\cdot)$.
    
    \jp{We do not consider the case where it can be both since Proposition \ref{prop:efficient_policy} guarantees the existence of an optimal policy under which each queue is either one or the other. (Although Proposition \ref{prop:efficient_policy} is stated and proved after Theorem \ref{thm:properties_V}, it does not rely on Theorem \ref{thm:properties_V}.)} This concludes that $V^{M-1}(\cdot)$ is non-decreasing.
    
    Now, suppose the claim holds for period $m+1, \ldots, M-1$. In period $m$, given an initial condition $x$, 
    \begin{equation*}
        V^m(x) = \min_{u \in \mathcal{U}(x)} \left[ \jp{H^m}(x + \phi(u)) + R(\phi(u)) + \kappa(\phi(u)) + V^{m+1}(f^m(x + \phi(u),\tau)) \right].
    \end{equation*}
    Since $V^{m+1}(\cdot)$ is non-decreasing \jp{by the induction hypothesis} and $f^m(\cdot,\tau)$ is non-decreasing under Assumptions \ref{ass:time_varying_arrivals} and \ref{ass:non_linear_holding_cost} based on its recursive definition \eqref{eq:recursive_f}, we can use the same argument as before to show that $V^m(y) \geq V^m(z)$ for any $y \geq z$. Thus, $V^m(\cdot)$ is non-decreasing for all $m\in \mathcal{M}$.

    Using induction, we next show that $V^{m}(\cdot)$ is $K$-convex for all $m\in \mathcal{M}$ by verifying that for any two states $x^1$ and $x^2$, $V^m(\cdot)$ satisfies Definition \ref{def:K_convexity} under the joint setup cost in $\eqref{eq:joint_setup_cost}$ with parameter \jp{$K \geq 0$}.

    Consider the last period $M-1$. Given an initial condition $x$,
    \begin{align*}
        V^{M-1}(x) = \min_{y \in \Delta(\mathrm{e}^\top x)}[\jp{H^{M-1}}(y) + R(y-x) + \kappa(y-x)].
    \end{align*}
    Fix $n\geq0$. Let $x^i \in \Delta(n)$ for $i=1,2$ and assume $x^1\neq x^2$ without loss of generality. Let $\bar{x} = \theta x^1 + (1-\theta)x^2$ where $\theta\in[0,1]$. \jp{If it is not optimal to transfer at $x^1$ and $x^2$, then}
    \begin{align*}
        V^{M-1}(\bar{x}) &\leq \jp{H^{M-1}}(\bar{x}) \\
        &\leq \theta \jp{H^{M-1}}(x^1) + (1-\theta)[\jp{H^{M-1}}(x^2) + \kappa(x^2-x^1)] \\
        &= \theta V^{M-1}(x^1) + (1-\theta)[V^{M-1}(x^2) + \kappa(x^2-x^1)].
    \end{align*}
    \jp{The first equality follows since $\bar{x}$ may not be an optimal solution in general.} The second inequality holds by convexity of $\jp{H^m}(\cdot)$ (Lemma \ref{lem:properties_H}) and because $\kappa(x^2-x^1)\geq0$. \jp{Since it is not optimal to transfer at neither $x^1$ nor $x^2$, the last step follows from $V^{M-1}(x^i)=H^{M-1}(x^i)$ for $i=1,2$. This shows $V^{M-1}(\cdot)$ is $K$-convex.}
    
    \jp{So, assume that it is optimal to transfer for at least one of $x^1$ or $x^2$, and suppose, without loss of generality, that it is $x^1$.} Let $\epsilon > 0$. Then there exist $y^i \in \Delta(n)$ for $i=1,2$, such that
    \begin{equation*}
        \jp{H^{M-1}}(y^i) + R(y^i - x^i) + \kappa(y^i - x^i) \leq V^{M-1}(x^i) + \epsilon, \quad i=1,2.
    \end{equation*}
    In particular, \jp{if it is optimal to transfer at $x^i$}, we may assume that there always exists $y^i\neq x^i$ that satisfies the above inequality for any $\epsilon>0$. Now, let $\theta \in [0,1]$, and for ease of notation, let $\Bar{x} = \theta x^1 + (1-\theta)x^2$ and $\Bar{y} = \theta y^1 + (1-\theta)y^2$. We observe
    \begin{equation*}
    \begin{split}
        V^{M-1}(\Bar{x}) &= \min_{y \in \Delta(n)} \left[ \jp{H^{M-1}}(y) + R(y - \Bar{x}) + \kappa(y - \Bar{x})\right] \\
        &\leq \jp{H^{M-1}}(\Bar{y}) + R(\Bar{y} - \Bar{x}) + \kappa(\Bar{y} - \Bar{x}) \\
        &\leq \theta \left[ \jp{H^{M-1}}(y^1) + R(y^1 - x^1) + \kappa(\Bar{y} - \Bar{x}) \right] + (1-\theta) \left[ \jp{H^{M-1}}(y^2) + R(y^2 - x^2) + \kappa(\Bar{y} - \Bar{x}) \right] \\
        &\leq \theta \left[ \jp{H^{M-1}}(y^1) + R(y^1 - x^1) + \kappa(y^1 - x^1) \right] + (1-\theta) \left[ \jp{H^{M-1}}(y^2) + R(y^2 - x^2) + \kappa(y^2-x^2) + \kappa(x^2 - x^1) \right] \\
        &\leq \theta V^{M-1}(x^1) + (1-\theta)\left[ V^{M-1}(x^2) + \kappa(x^2-x^1) \right] + \epsilon. \\
    \end{split}
    \end{equation*}
    The first inequality holds by $\Bar{y}$ being a feasible solution. The second inequality holds by substituting the definitions of $\Bar{x}$ and $\Bar{y}$ and by convexity of $\jp{H^{m}}(\cdot)$ and $R(\cdot)$ (Lemma \ref{lem:R_properties}). For the third inequality, we use the fact that $\kappa(\bar{y}-\bar{x}) = \kappa(y^1-x^1)$ and $\kappa(\bar{y}-\bar{x}) \leq \kappa(y^2-x^2) + \kappa(x^2 - x^1)$, \jp{which follows from $\kappa(x^2 - x^1) = K$ (we assumed $x_1 \neq x_2$)}. Since this observation holds for any $\epsilon > 0$, we have
    \begin{equation*}
        V^{M-1}(\Bar{x}) \leq \theta V^{M-1}(x^1) + (1-\theta)\left[ V^{M-1}(x^2) + \kappa(x^2 - x^1) \right].
    \end{equation*} This shows that $V^{M-1}(\cdot)$ is $K$-convex.

    Now, suppose the claim holds for periods $m+1, \ldots, M-1$. We note that by monotonicity of $V^m(\cdot)$ for all $m\in \mathcal{M}$, and by convexity of the state transition function $f^m(\cdot,\tau)$ under Assumptions \ref{ass:time_varying_arrivals} and \ref{ass:non_linear_holding_cost} based on its recursive definition \eqref{eq:recursive_f}, $V^t(f^{t-1}(\cdot,\tau))$ is $K$-convex for $t\in \{m+1,\ldots,M-1\}$. In period $m$, again fix $n\geq0$, $\epsilon>0$, and let $x^i \in \Delta(n)$ for $i=1,2$, where $x^1\neq x^2$ without loss of generality. \jp{If it is not optimal to transfer at $x^1$ and $x^2$}, then by a similar argument as before, we see that
    \begin{align*}
        V^m(\bar{x}) &\leq \jp{H^m}(\bar{x}) + V^{m+1}(f^m(\bar{x},\tau)) \\
        &\leq \theta[\jp{H^m}(x^1) + V^{m+1}(f^m(x^1,\tau))] + (1-\theta)[\jp{H^m}(x^2) + V^{m+1}(f^m(x^2,\tau)) +  \kappa(x^2-x^1)] \\
        &= \theta V^m(x^1) + (1-\theta)[V^m(x^2) + \kappa(x^2-x^1)],
    \end{align*}
    i.e., $V^m(\cdot)$ is $K$-convex. The second step holds by convexity of $\jp{H^m}(\cdot)$ and $K$-convexity of $(V^{m+1} \circ f^m)(\cdot)$. 
    
    So, assume that \jp{it is optimal to transfer for at least one of $x^1$ or $x^2$, and suppose, without loss of generality, that it is $x^1$.} Then there exist $y^i \in \Delta(n)$ for $i=1,2$, such that
    \begin{equation*}
        \jp{H^m}(y^i) + R(y^i - x^i) + \kappa(y^i - x^i) \leq V^{m}(x^i) + \epsilon, \quad \forall i=1,2.
    \end{equation*}
    Letting $\theta \in [0,1]$ and defining $\Bar{x} = \theta x^1 + (1-\theta)x^2$ and $\Bar{y} = \theta y^1 + (1-\theta)y^2$, the rest of the proof follows a similar argument as before:
    \begin{equation*}
    \begin{split}
        V^m(\Bar{x}) &= \min_{y \in \Delta(n)} \left[ \jp{H^m}(y) + R(y - \Bar{x}) + \kappa(y - \Bar{x}) + V^{m+1}(f^m(y,\tau)) \right] \\
        &\leq \jp{H^m}(\Bar{y}) + R(\Bar{y} - \Bar{x}) + \kappa(\Bar{y} - \Bar{x}) + V^{m+1}(f^m(\bar{y},\tau)) \\
        &\leq \theta \left[ \jp{H^m}(y^1) + R(y^1 - x^1) + \kappa(y^1 - x^1) + V^{m+1}(f^m(y^1,\tau)) \right] \\
            &\quad + (1-\theta) \left[ \jp{H^m}(y^2) + R(y^2 - x^2) + \kappa(y^2-x^2) + V^{m+1}(f^m(y^2,\tau)) + \kappa(x^2 - x^1) \right] \\
        &\leq \theta V^m(x^1) + (1-\theta)\left[ V^m(x^2) + \kappa(x^2-x^1) \right] + \epsilon. \\
    \end{split}
    \end{equation*}
    Since this observation holds for any $\epsilon > 0$, this shows that $V^m(\cdot)$ is $K$-convex for all $m\in \mathcal{M}$.
    
    Lastly, we show continuity of $V^m(\cdot)$. 
    Define $v^m(y;x) = \jp{H^m}(y) + R(y-x) + V^{m+1}(f^m(y,\tau))$ where $R(y-x)$ is the variable transfer cost of going from a given initial condition $x$ to $y$. Let $n$ be the total number of customers at $x$.
    For period $m=M-1$, $v^{M-1}(\cdot)$ is clearly continuous \jp{in $y$} on $\Delta(n)$.
    Furthermore, since $\Delta(n)$ is compact, $v^{M-1}$ is uniformly continuous on $\Delta(n)$. 
    Then let any two points $x_1, x_2 \in \Delta(n)$ such that $\| x_1-x_2 \| < \delta$ for some $\delta>0$, \jp{where $\|\cdot \|$ is the Euclidean norm}, and let $y^*_1$ and $y^*_2$ denote their \jp{optimal post-transfer states}, respectively. \jp{We can consider two cases: (1) it is optimal to transfer at $x_1$, or (2) it is \textit{not} optimal to transfer at $x_1$. If it is the former case, there must exist $y^*_1 \neq x_1$. Then we can consider some $y_2 \in \jp{\Sigma^m}(n)$ such that $\| y^*_1 - y_2 \| < \delta$. (If there is a unique optimal post-transfer state, then the only choice is $y_2 = y^*_1$, in which case the analysis below becomes trivial by continuity of $v^m(y;x)$ in $x$.)} Then by uniform continuity of $v^{M-1}$ on $\Delta(n)$, this implies that there exists $\epsilon > 0$ such that, 
    \begin{align*}
        V^{M-1}(x_1) = K + v^{M-1}(y^*_1;x_1) \geq K + v^{M-1}(y_2;x_2) - \epsilon \geq V^{M-1}(x_2) - \epsilon,
    \end{align*}
    holds. \jp{In the latter case where it is not optimal to transfer at $x_1$, we must have $y^*_1 = x_1$.} So,
    \begin{align*}
        V^{M-1}(x_1) = v^{M-1}(x_1;x_1) \geq v^{M-1}(x_2;x_2) - \epsilon \geq V^{M-1}(x_2) - \epsilon.
    \end{align*}
    In both cases, we have $V^{M-1}(x_1) \geq V^{M-1}(x_2) - \epsilon$.
    Using the same argument, we can show that $V^{M-1}(x_2) \geq V^{M-1}(x_1) - \epsilon$. This indicates that $| V^{M-1}(x_1) - V^{M-1}(x_2)| \leq \epsilon$ and hence $V^{M-1}(\cdot)$ is continuous. 
    Now, assume $V^{m+1}(\cdot)$ is continuous for some $m \leq M-2$. Using continuity of $f^m(\cdot,\tau)$, it follows that $v^m(\cdot)$ is also continuous. Therefore, following the same analysis as period $M-1$, $V^m(\cdot)$ is continuous.
\end{proof}

\subsubsection{Additional Properties of the Value Function.} \label{appen:add_properties_V}
Finally, we outline additional properties of a $K$-convex function. These properties were first demonstrated by \cite{gallegoKConvexityRn2005a} in $\mathbb{R}^N$ for the $N$-product inventory control problem. We show these results for the value function under the joint setup cost in \eqref{eq:joint_setup_cost}, which will be useful for characterizing the compactness and connectedness of the no-transfer region $\jp{\Sigma^m}(\cdot)$ in the proof of Theorem \ref{thm:structure}.
\begin{lemma}\label{lem:K_convexity_properties}
    For all $m \in \mathcal{M}$, we have:
    \begin{enumerate}
        \item[(i)] $V^m(\cdot)$ is $L$-convex for any $L \geq K$.
        \item[(ii)] If $W(\cdot)$ is $L$-convex, then for any $\alpha\geq0$, $\beta \geq 0$, $\alpha V^m(\cdot) + \beta W(\cdot)$ is $(\alpha K + \beta L)$-convex.
        \item[(iii)] Let $x\in \mathbb{R}^N_+$ and $y \in \Delta(\mathrm{e}^\top x)$. Suppose $V^{m+1}(f^m(\cdot,\tau))$ is $K$-convex. Define $g^m:[0,1] \xrightarrow{} \mathbb{R}_+$ as 
        \begin{equation*}
            g^m(\theta) = H^m(x + \theta(y - x)) + R(\theta(y - x)) + V^{m+1}(f^m(x + \theta(y - x), \tau)).
        \end{equation*}
        Then $g^m(\cdot)$ is $K$-convex in the univariate sense defined by \cite{scarfOptimalityPoliciesDynamic1960a} (Section 3, Equation 8).
    \end{enumerate}
\end{lemma}

\begin{proof}{Proof of Lemma \ref{lem:K_convexity_properties}.}

    The first property follows directly from $K1\{z\neq0\} \leq L1\{z\neq0\}$ for any $z\in \mathbb{R}^N$ and $K\leq L$. To show the second property, we note that for any $x,y\in \mathbb{R}^N_+$, $x \neq y$, and $\theta\in[0,1]$:
    \begin{align*}
        (\alpha V^m + \beta W)(\theta x + (1-\theta)y) &= \alpha V^m(\theta x + (1-\theta)y) + \beta W(\theta x + (1-\theta)y) \\
        &\leq \alpha(\theta V^m(x) + (1-\theta)[V^m(y) + K1\{y-x\neq0\}]) \\
        &\quad + \beta(\theta W(x) + (1-\theta)[W(y) + L1\{y-x\neq0\}]) \\
        &= \theta[\alpha V^m(x) + \beta W(x)] + (1-\theta)[\alpha V^m(y) + \beta W(y) + \alpha K 1\{y-x\neq0\} + \beta L1\{y-x\neq0\}] \\
        &= \theta(\alpha V^m + \beta W)(x) + (1-\theta)[(\alpha V^m + \beta W)(y) + (\alpha K + \beta L)1\{y-x\neq0\}],
    \end{align*}
    i.e., $\alpha V^m(\cdot) + \beta W(\cdot)$ is $(\alpha K + \beta L)$-convex.

    Lastly, we prove the third property by contradiction. Suppose $g^m(\cdot)$ is not $K$-convex. Then there exist $\eta\in[0,1]$ and $\theta_1 \leq \theta_2$ such that,
    \begin{equation} \label{eq:lem4_inequality}
        g^m(\eta\theta_1 + (1-\eta)\theta_2) > \eta g^m(\theta_1) + (1-\eta)[g^m(\theta_2) + K].
    \end{equation}
    For ease of notation, let $\Bar{\theta} = \eta\theta_1 + (1-\eta)\theta_2$ and let $z_1 = x + \theta_1(y-x)$ and $z_2 = x + \theta_2(y-x)$. \jp{First, we expand the left-hand side of \eqref{eq:lem4_inequality}:
    \begin{equation*}
        \begin{split}
            g^m(\Bar{\theta}) &= H^m(x + \Bar{\theta}(y-x)) + R(\Bar{\theta}(y-x)) + V^{m+1}(f^m(x + \Bar{\theta}(y - x),\tau)) \\
            &= H^m(\eta z_1 + (1-\eta)z_2) + \Bar{\theta}R(y-x) + V^{m+1}(f^m(\eta z_1 + (1-\eta)z_2,\tau)).\\
        \end{split}
    \end{equation*}
    Similarly, we expand the right-hand side of \eqref{eq:lem4_inequality}:
    \begin{equation*}
        \begin{split}
            \eta g^m(\theta_1) + (1-\eta)[g^m(\theta_2) + K] &= \eta H^m(z_1) + (1-\eta)H^m(z_2) + \eta R(\theta_1(y-x)) + (1-\eta) R(\theta_2(y-x)) \\
            &\quad + \eta V^{m+1}(f^m(z_1,\tau)) + (1-\eta)[V^{m+1}(f^m(z_2,\tau)) + \kappa(y-x)] \\
            &= \eta H^m(z_1) + (1-\eta)H^m(z_2) + \Bar{\theta}R(y-x) \\
            &\quad + \eta V^{m+1}(f^m(z_1,\tau)) + (1-\eta)[V^{m+1}(f^m(z_2,\tau)) + \kappa(y-x)].
        \end{split}
    \end{equation*}
    Putting both sides together and subtracting $\bar{\theta}R(y-x)$ from both, we conclude
    \begin{equation*}
        \begin{split}
            &H^m(\eta z_1 + (1-\eta)z_2) + V^{m+1}(f^m(\eta z_1 + (1-\eta)z_2,\tau)) \\
            &> \eta H^m(z_1) + (1-\eta)H^m(z_2) + \eta V^{m+1}(f^m(z_1,\tau)) + (1-\eta)[V^{m+1}(f^m(z_2,\tau)) + \kappa(y-x)].
        \end{split}
    \end{equation*}}
    However, we note that, 
    \begin{equation*}
        \begin{split}
            H^m(\eta z_1 + (1-\eta)z_2) &\leq \eta H^m(z_1) + (1-\eta)H^m(z_2),
        \end{split}
    \end{equation*}
    must hold by Lemma \ref{lem:properties_H}, and
    \begin{equation*}
        \begin{split}
            V^{m+1}(f^m(\eta z_1 + (1-\eta)z_2,\tau)) &\leq \eta V^{m+1}(f^m(z_1,\tau)) + (1-\eta)[V^{m+1}(f^m(z_2,\tau)) + \kappa(z_2-z_1)] \\
            &\leq \eta V^{m+1}(f^m(z_1,\tau)) + (1-\eta)[V^{m+1}(f^m(z_2,\tau)) + \kappa(y-x)],
        \end{split}
    \end{equation*}
    must hold by the assumption in the statement. The second inequality holds because $\kappa(z_2-z_1) \leq \kappa(y-x)$, as shown below:
    \begin{equation*}
        \kappa(z_2-z_1) = \kappa((\theta_2-\theta_1) (y-x)) =
        \begin{cases}
            0, & \theta_2-\theta_1=0 \\
            \kappa(y-x), & \theta_2-\theta_1>0
        \end{cases}
        \leq \kappa(y-x),
    \end{equation*}
    Therefore, this is a contradiction, and as a result, $g^m(\cdot)$ must be $K$-convex.
\end{proof}

\subsection{Proofs of the Results in Section \ref{sec:structure}}
\subsubsection{Proof of Proposition \ref{prop:efficient_policy}.}
\begin{proof}{Proof of Proposition \ref{prop:efficient_policy}.}
    Consider an initial condition $x\in\mathbb{R}^N_+$ and a post-transfer state $y\in\Delta(\mathrm{e}^\top x)$, $y\neq x$. For fixed $y$, the holding cost $H^m(y)$, the value function $V^{m+1}(f^m(y,\tau))$, and the joint setup cost $\kappa(y-x)=K$ are constant. Therefore, it suffices to show that there exists an optimal transfer decision matrix $u^*$ such that its cost is $R(y-x)$ and it satisfies $u^*_{ij}u^*_{jl}=0, \forall i,j,l\in \mathcal{N}$. This means that if queue $j$ is receiving customers ($u^*_{ij}>0$ for some $i$), queue $j$ cannot be sending customers to any queue at the same time ($u^*_{jl}=0, \forall l$) and vice versa.

    Suppose that there exist $i,j,l \in \mathcal{N}$ such that $u^*_{ij}u^*_{jl}>0$. We will construct another feasible transfer matrix $\hat{u}$ in the following way. If $i=l$, then we can simply force to zero the smaller of $u^*_{ij}$ and $u^*_{ji}$. Without loss of generality, assume that $u^*_{ij} = \min\{u^*_{ij},u^*_{ji}\}$. Then we can set $\hat{u}_{ij}=0$ and $\hat{u}_{ji}= u^*_{ji} - u^*_{ij} > 0$. By implementing $\hat{u}$, we would reduce the total cost by $r_{ij}u^*_{ij} \geq 0$. So, $\hat{u}$ does equally well, if not better, than $u^*$. If $i\neq l$, we have two cases: $u^*_{ij} \geq u^*_{jl}$ or $u^*_{ij} < u^*_{jl}$. In the first case, we can set $\hat{u}_{jl}=0$, $\hat{u}_{ij} = u^*_{ij} - u^*_{jl} \geq 0$, and $\hat{u}_{il} = u^*_{il} + u^*_{jl}$. Then we would reduce the total cost by $r_{ij}u^*_{jl} + r_{jl}u^*_{jl} - r_{il}u^*_{jl} = u^*_{jl}(r_{ij} + r_{jl} - r_{il}) \geq 0$ (Assumption \ref{ass:triangular_inequality}). So, $\hat{u}$ again does equally well or better. In the second case that $u^*_{ij} < u^*_{jl}$, we can following a similar analysis and reach the same conclusion. Therefore, there always exists an optimal policy under which each queue is either sending or receiving customers, but not both, in the same period.
\end{proof}

\subsubsection{Proof of Theorem \ref{thm:structure}.} \label{appen:proof_structure}
\jp{For clarity, we first show the following auxiliary result which will be used in the proof of Theorem \ref{thm:structure}. Intuitively, the following result states that every state on the boundary has a corresponding alternative state with equivalent costs, whether we stay or move.}
\jp{
\begin{lemma} \label{lem:boundary_is_a_subset}
    Suppose $\kappa(\cdot)$ is the joint setup cost function \eqref{eq:joint_setup_cost}. Then $\partial \jp{\Sigma^m}(n) \subseteq \partial \tilde{\Sigma}$, where 
    \begin{align*}
         \partial \tilde{\Sigma} = \{x \in \Delta(n): g(x,y) = K \text{ for some } y\in\Delta(n), y\neq x\},
    \end{align*}
    and
    \begin{align*}
        g(x,y) = \jp{H^m}(x) + V^{m+1}(f^m(x,\tau)) - \jp{H^m}(y) - R(y-x) - V^{m+1}(f^m(y,\tau)).
    \end{align*}
\end{lemma}
}
\jp{
\begin{proof}{Proof.}
    Let $x \in \partial \Sigma^m(n)$. We want to show $x \in \partial \tilde{\Sigma}$. Since $x \in \partial \Sigma^m(n)$, the open ball $B_\epsilon(x)$ for any $\epsilon > 0$ contains points $z$ such that $g(z,y) \leq K$ for all $y \in \Delta(n)$, $y\neq z$, and points $z'$ such that $g(z',y) \geq K$ for some $y\in\Delta(n)$, $y\neq z'$. Consider a sequence $\{z_n\}$ in $\Sigma^m(n)$ that converges to $x$ and a sequence $\{z'_n\}$ in $\Sigma^m(n)^c$ (complement) that also converges to $x$. Since $g(x,y)$ is continuous in $x$ and $y$ (by Lemma \ref{lem:properties_H}, Lemma \ref{lem:R_properties}, and Theorem \ref{thm:properties_V}), we have $g(z_n, y) \rightarrow g(x,y)$ with each $g(z_n,y) \leq K$ for all $y \in \Delta(n)$, $y\neq z_n$, and $g(z'_n, y) \rightarrow g(x,y)$ with each $g(z'_n, y) \geq K$ for some $y\in\Delta(n)$, $y\neq z'_n$ (possibly with different $y$ for each $z'_n$). Since the values of $g(z_n,y)$ are always at most $K$ and $g(z'_n,y)$ are always at least $K$, the only possibility for the limit $g(x,y)$ is exactly $K$, which implies $x \in \partial \tilde{\Sigma}$.
\end{proof}
}

\begin{proof}{Proof of Theorem \ref{thm:structure}.}
In the following, we first establish the properties of the no-transfer region $\jp{\Sigma^m}(n)$ in the order of non-emptiness, compactness, and connectedness, which is followed by the properties of the optimal policy.

    Fix $n\geq0$. We first show that $\jp{\Sigma^m}(n)$ is non-empty. For any given initial condition $x\in\Delta(n)$, denote its target state by $x^*$, which exists since $\Delta(n)$ is non-empty and compact. It is then easy to see that $x^* \in \jp{\Sigma^m}(n)$. Suppose otherwise. Then, there exists $y \in \Delta(n)$, $y \neq x^*$, such that \jp{$H^m(x^*) + V^{m+1}(f^m(x^*,\tau)) > V^m(x^*) = H^m(y) + R(y - x^*) + \kappa(y - x^*) + V^{m+1}(f^m(y,\tau))$. Thus, we conclude: 
    \begin{align*}
            V^m(x) &= R(x^* - x) + \kappa(x^* - x) + \jp{H^m}(x^*) + V^{m+1}(f^m(x^*,\tau)) \\
            &> R(x^* - x) + \kappa(x^* - x) + \jp{H^m}(y) + R(y - x^*) + \kappa(y - x^*) + V^{m+1}(f^m(y,\tau)) \\
            &> R(y - x) + \kappa(y - x) + H^m(y) + V^{m+1}(f^m(y,\tau)),
    \end{align*}
    where the last inequality holds by subadditivity of $\kappa(\cdot)$ (Lemma \ref{lem:joint_setup_cost_properties}) and $R(\cdot)$ (Lemma \ref{lem:R_properties}).} This contradicts the optimality of $x^*$ as it suggests that moving to $y$ yields a strictly lower cost. Therefore, $x^* \in \jp{\Sigma^m}(n)$, and $\jp{\Sigma^m}(n)$ is non-empty.

    Secondly, we show that $\jp{\Sigma^m(n)}$ is compact. Since $\jp{\Sigma^m(n)} \subseteq \Delta(n)$, and $\Delta(n)$ is compact, it suffices to show that $\jp{\Sigma^m(n)}$ is closed, \jp{which we show by proving that it contains all of its boundary points. Indeed, Lemma \ref{lem:boundary_is_a_subset} shows that $\partial \jp{\Sigma^m(n)} \subseteq \partial \tilde{\Sigma}$, and since $\partial \tilde{\Sigma}$ is clearly a subset of $\jp{\Sigma^m(n)}$, we must have $\partial \jp{\Sigma^m(n)} \subseteq \jp{\Sigma^m(n)}$, and $\jp{\Sigma^m(n)}$ is closed.}

    
    Thirdly, we show that $\jp{\Sigma^m(n)}$ is connected. For a given initial condition $x$, we let 
    \begin{equation*}
    \begin{split}
        y^*(x) = \{y \in \Delta(n) &: \jp{H^m}(y) + R(y - x) + \kappa(y - x) + V^{m+1}(f^m(y,\tau)) \\  
        &\leq \jp{H^m}(z) + R(z- x) + \kappa(z - x) + V^{m+1}(f^m(z,\tau)), \forall z \in \Delta(n), z \neq y \},
    \end{split}
    \end{equation*}
    be the set of all target states corresponding to $x$. Then we have $\jp{\Sigma^m(n)} = \cup_{x \in \Delta(n)} y^*(x)$. To see this, consider any $x\in \Delta(n)$. If $x \notin \jp{\Sigma^m(n)}$, $y^*(x)$ is the set of all equally optimal target states corresponding to $x$, all of which must be contained in $\jp{\Sigma^m(n)}$. If $x \in \jp{\Sigma^m(n)}$, $y^*(x)$ includes $x$ itself and possibly others in $\jp{\Sigma^m(n)}$ to which we are indifferent with regard to transferring or not. Repeating this for every $x\in \Delta(n)$, one can see that the union of $y^*(x)$ must be equal to $\jp{\Sigma^m(n)}$.
    
    Now, consider $y_1$ and $y_2$ in $y^*(x)$ such that $y_1 \neq y_2$. If $y_1=y_2$, we note that the proof below becomes trivial. For ease of notation, let $\bar{y} = \theta y_1 + (1-\theta)y_2$ for some $\theta \in (0,1)$. We show below that $\Bar{y} \in \jp{\Sigma^m(n)}$, which will ultimately be useful in proving that $\jp{\Sigma^m(n)}$ is connected. Observe that,
    \begin{equation*}
        \begin{split}
            \jp{H^m}(\Bar{y}) + R(\Bar{y} - x) + V^{m+1}(f^m(\Bar{y},\tau)) &\leq \theta[\jp{H^m}(y_1) + R(y_1 - x) + V^{m+1}(f^m(y_1,\tau))] \\
            &+ (1-\theta)[\jp{H^m}(y_2) + R(y_2 - x) + V^{m+1}(f^m(y_2,\tau)) + \kappa(y_2 - y_1)] \\
            &\leq \jp{H^m}(z) + R(z - x) + V^{m+1}(f^m(z,\tau)) + K, \quad \forall z \in \Delta(n).
        \end{split}
    \end{equation*}
    The first inequality holds by convexity of $\jp{H^m}(\cdot)$ and $R(\cdot)$ (Lemmas \ref{lem:properties_H} and \ref{lem:R_properties}) and $K$-convexity of $V^{m+1}(f^m(\cdot,\tau))$ (Theorem \ref{thm:properties_V}). The second inequality holds by definition of $y^*(x)$. Then, rearranging $R(\Bar{y} - x)$ to the right-hand side, we observe
    \begin{equation*}
        \begin{split}
            \jp{H^m}(\Bar{y}) + V^{m+1}(f^m(\Bar{y},\tau)) &\leq \jp{H^m}(z) + R(z - x) - R(\bar{y} - x) + V^{m+1}(f^m(z,\tau)) + K, \quad \forall z \in \Delta(n) \\
            &\leq \jp{H^m}(z) + R(z - \bar{y}) + V^{m+1}(f^m(z,\tau)) + K, \quad \forall z \in \Delta(n),
        \end{split}
    \end{equation*}
    where the last inequality holds by subadditivity of $R(\cdot)$ (Lemma \ref{lem:R_properties}). This result indicates that it is better to remain at $\Bar{y}$ than to move to any other states in $\Delta(n)$. Therefore, $\Bar{y} \in \jp{\Sigma^m(n)}$. Intuitively, what we have shown is that any convex combination of two target states corresponding to $x$ must lie in $\jp{\Sigma^m(n)}$.

    Suppose now, for sake of contradiction, that $\jp{\Sigma^m(n)}$ is not connected. Then $\jp{\Sigma^m(n)}$ can be expressed as a union of two non-empty, separated sets, i.e., \jp{$\Sigma^m(n)= \mathcal{V}_1 \cup \mathcal{V}_2$ where $\mathcal{V}_1 \neq \varnothing$, $\mathcal{V}_2 \neq \varnothing$, and $\mathrm{cl}(\mathcal{V}_1) \cap \mathcal{V}_2 = \mathcal{V}_1 \cap \mathrm{cl}(\mathcal{V}_2) = \varnothing$, where $\mathrm{cl}(\cdot)$ represents the closure of a given set.} We note that for all $x\in \Delta(n)$, $y^*(x)$ must lie in either $\mathcal{V}_1$ or $\mathcal{V}_2$, but not both. To see this, suppose that there exist $y_1$ and $y_2$ in $y^*(x)$ such that $y_1 \in \mathcal{V}_1$ and $y_2 \in \mathcal{V}_2$. \jp{Since $\mathcal{V}_1$ and $\mathcal{V}_2$ are separated,} there must exist some $\theta \in (0,1)$ such that $\Bar{y} = \theta y_1 + (1-\theta) y_2 \notin \jp{\Sigma^m(n)}$. However, we have just shown that all convex combinations of two target states corresponding to $x$ must lie in $\jp{\Sigma^m(n)}$, which contradicts that $y_1 \in \mathcal{V}_1$ and $y_2 \in \mathcal{V}_2$. Now, let $\mathcal{U}_1 = y^{*-1}(\mathcal{V}_1)$ and $\mathcal{U}_2 = y^{*-1}(\mathcal{V}_2)$. Then $\mathcal{U}_1$ and $\mathcal{U}_2$ are two non-empty, separated sets such that $\Delta(n) = \mathcal{U}_1 \cup \mathcal{U}_2$. This implies that $\Delta(n)$ is in fact not connected, which is a contradiction. Therefore, $\jp{\Sigma^m(n)}$ must be connected.

    Finally, we establish the properties of the optimal policy. Note that if $x \in \jp{\Sigma^m(n)}$, the existence of a target state $y$ such that $y = x$ is clear from the definition of $\jp{\Sigma^m(n)}$, meaning that it is optimal not to move. So, assume $x \notin \jp{\Sigma^m(n)}$. We consider the three cases from Theorem \ref{thm:structure} in order:
    \begin{itemize}
        \item If $\kappa(\cdot)=0$ and $r_{ij}=0$ for all $i,j$, we note that $V^m(\cdot)$ is convex for all $m\in \mathcal{M}$ (\jp{Corollary \ref{cor:properties_V_no_setup}}) and $V^m(x) = \min_{y \in \Delta(n)}[\jp{H^m}(y) + V^{m+1}(f^m(y,\tau))]$. This suggests that the target state $y$ can be obtained by solving a convex optimization problem over a compact set whose cost is independent of $x$. Thus, there must exist a global target state $y$ such that it is optimal to move to $y$ from any $x \notin \jp{\Sigma^m(n)}$.
        \item Suppose $\kappa(\cdot)=0$ and let $x^*$ denote a target state corresponding to $x$ such that $x^* \in \mathrm{ri}(\jp{\Sigma^m(n)})$. Then for small enough $\theta \in (0,1)$, there must exist $y=x^* + \theta(x-x^*) \in \jp{\Sigma^m(n)}$ and:
        \begin{align*}
            R(y-x) + \jp{H^m}(y) + V^{m+1}(f^m(y,\tau)) &\leq R(y-x) + R(x^*-y) + \jp{H^m}(x^*) + V^{m+1}(f^m(x^*,\tau)) \\
            &\leq R(y-x) + R(x^*-x) + \jp{H^m}(x^*) + V^{m+1}(f^m(x^*,\tau)) \\
            &= (1-\theta)R(x^*-x) + \theta R(x^*-x) + \jp{H^m}(x^*) + V^{m+1}(f^m(x^*,\tau)) \\
            &= R(x^*-x) + \jp{H^m}(x^*) + V^{m+1}(f^m(x^*,\tau)).
        \end{align*}
        The first inequality \jp{follows from $y \in \jp{\Sigma^m(n)}$}. \jp{The second inequality holds since $y$ lies part way on the line segment connecting $x^*$ and $x$, which implies $R(x^*-y) \leq R(x^*-x)$.} The first equality \jp{on the third line} holds by substituting the definition of $y$. This shows that going to $y$ is just as good, if not better, than going to $x^*$ from $x$. Therefore, there must exist a target state $y \in \partial \jp{\Sigma^m(n)}$.
        \item Let $\kappa(\cdot)$ be the joint setup cost function \eqref{eq:joint_setup_cost}. Let $x^*$ denote a target state corresponding to $x$, $x^* \neq x$. Suppose $x^* \in \partial \jp{\Sigma^m(n)}$. Then there must exist $y \in \jp{\Sigma^m(n)}, y\neq x^*$, such that we are indifferent to staying at $x^*$ or to moving from $x^*$ to $y$ \jp{(Lemma \ref{lem:boundary_is_a_subset})}. However, this leads to a contradiction:
        \begin{align*}
                V^m(x) &= \jp{H^m}(x^*) + R(x^*-x) + \kappa(x^* - x) + V^{m+1}(f^m(x^*,\tau)) \\
                &= R(x^*-x) + \kappa(x^* - x) + R(y - x^*) + \kappa(y - x^*) + \jp{H^m}(y) + V^{m+1}(f^m(y,\tau)) \\
                &> R(y - x) + \kappa(y - x) + \jp{H^m}(y) + V^{m+1}(f^m(y,\tau)),
        \end{align*}
        i.e., $x^*$ is not a target state corresponding to $x$. Therefore, all target states belong to $\mathrm{ri}(\jp{\Sigma^m(n)})$. 
    \end{itemize}
\end{proof}

\subsubsection{\jp{Proof of Proposition \ref{prop:optimal_two_queue_policy}.}}
\begin{proof}{\jp{Proof of Proposition \ref{prop:optimal_two_queue_policy}.}}
    \jp{Proposition \ref{prop:optimal_two_queue_policy} follows as a direct consequence of Theorem \ref{thm:structure} except for the claim that the target states $(S_1, n-S_1)$ and $(n-S_2,S_2)$ depend only on the total number of customers $n = x_1+x_2$, not the entire vector $(x_1,x_2)$. We prove this result below.}

    \jp{Consider $x \notin \Sigma^m(n)$. Suppose $x_1 < s_1(n)$, so the optimal policy transfers customers from queue 2 to 1. Let the resulting target state be $(S_1(x), n-S_1(x))$. Then this target state satisfies:
    \begin{align*}
        (S_1(x), n-S_1(x)) \in \argmin_{y:y_1+y_2=n} \left\{ H^m(y) +  V^{m+1}(f^m(y,\tau)) + r_{21}y_1 - r_{21}x_1 + K \right\},
    \end{align*}
    Note that the feasible set of the minimization problem does \textit{not} depend on the entire vector $(x_1,x_2)$, but only on its sum $n = x_1+x_2$. Moreover, the optimal choice of $y$ should stay constant with $x$, as $x$ enters the objective function only as a constant with respect to $y$. Thus, it follows that $(S_1(x), n-S_1(x))$ is the same target state for every $x \notin \jp{\Sigma^m(n)}$ satisfying $x_1 < s_1(n)$. By the same argument, the target state is invariant for every $x \notin \jp{\Sigma^m(n)}$ satisfying $x_2 < s_2(n)$. Since Proposition \ref{prop:efficient_policy} guarantees the existence of an optimal policy where either $x_1 < s_1(n)$ or $x_2 < s_2(n)$ holds, but not both, the proof is complete.}
\end{proof}

\subsection{Proof of the Results in Section \ref{sec:trade_off}}
Since the queueing dynamics are assumed to be stationary in Section \ref{sec:trade_off}, we omit the dependence of the holding cost function on the period.
\subsubsection{Proof of Proposition \ref{prop:trade_off}.}
\begin{proof}{Proof of Proposition \ref{prop:trade_off}.}
    We will first show the result for a two-queue system and extend the argument to a general $N$-queue system. We prove the contrapositive: when $h_j \geq h_i$ and $x_j \geq \tau(\mu_j-\lambda_j)$, a policy which does not transfer customers from queue $i$ to $j$ is better than (or as good as) a policy which does. \jp{Note that the proof of part $(ii)$ follows directly from the proof of part $(i)$.}
    
    Without loss of generality, suppose $h_2 \geq h_1$. We compare the total costs of two processes under two different transfer policies starting from the same initial condition. Consider two processes $x^m=(x^m_1, x^m_2)$ and $y^m=(y^m_1,y^m_2)$, $m\in \mathcal{M}$. These processes represent the queue lengths at the start of each period just \emph{after} transfers. We will denote the states \emph{prior} to transfers by $q^\pi[m]=(q^\pi_1[m],q^\pi_2[m])$, with the initial condition being $q^\pi[0]$. Suppose $q_2^\pi[0] \geq \tau(\mu_2-\lambda_2)^+$. The first process follows a policy $\pi = \{\pi^m\}_{m \in \mathcal{M}}$ which calls for transferring $u>0$ customers from queue 1 to 2 in period 0. This is denoted by $\pi^0 = u$, where $\pi^m > 0$ indicates that customers are transferred from queue 1 to 2. The second process follows another policy $\tilde{\pi}$ which is the same as $\pi$ except that in period 0, the policy $\tilde{\pi}$ does \emph{not} move customers from queue 1 to 2, i.e., $\tilde{\pi}^0 = 0$. Then at the start of period 1, just prior to transferring, we are in one of two scenarios: either (1) $q^\pi_1[1] >0$, or (2) $q^\pi_1[1]=0$.
    
    \textbf{Case 1. } In the first case, it implies that $x^0_1 > \tau(\mu_1-\lambda_1)^+$ (after transferring $\pi^0=u$ to queue 2). \jp{We can then set $\tilde{\pi}^1=\pi^1+u$} and the two processes will coincide in period 1. We then let $\pi=\tilde{\pi}$ thereafter. Denote by $\Delta_{\pi-\tilde{\pi}} \in \mathbb{R}$ the total cost of the first process minus that of the second process. We observe 
    \begin{align*}
        \Delta_{\pi-\tilde{\pi}} = \tau(h_2-h_1)u \geq 0,
    \end{align*}
    i.e., policy $\tilde{\pi}$ performs equally well or better. In calculating $\Delta_{\pi - \tilde{\pi}}$, the term $\tau h_2u$ represents the holding cost at queue 2 over period 0 when we follow policy $\pi$. This assumes that policy $\pi$ does not involve any transfers out of queue 2. Since Proposition \ref{prop:efficient_policy} ensures the existence of an optimal policy where no queues are both sending and receiving customers in the same period, one can assume this without loss of optimality.

    \textbf{Case 2. } The second case implies  $x^0_1 \leq \tau(\mu_1-\lambda_1)^+$ (after transferring $\pi^0=u$ to queue 2). This means that queue 1 under policy $\pi$ will empty before period 1 and $q^\pi_1[1]=0$ must hold. Moreover, $0 \leq q^{\tilde{\pi}}_1[1] \leq u$ holds, i.e., by the start of period 1, just prior to transferring any customers, the state of queue 1 in the second process cannot be larger than $u$. Thus, set $\pi^1=0$ and $\tilde{\pi}^1 = \hat{u} \equiv q^{\tilde{\pi}}_1[1]$, where $0 \leq \hat{u} \leq u$. We note that $\hat{u}=0$ may be the only feasible policy. Then $x^1_1=y^1_1=0$, i.e., the states of queue 1 under the two processes coincide. We let the two processes follow the respective optimal (fluid) trajectories thereafter. We observe that
    \begin{align*}
        \Delta_{\pi-\tilde{\pi}} &\geq ru - r\hat{u} + \tau h_2u - \tau h_1u + V^1(0,x^1_2) - V^1(0,y^1_2) \\
        &\geq r(u-\hat{u}) + \tau(h_2-h_1)u \\
        &\geq 0,
    \end{align*}
    where $V^1(\cdot)$ is the minimum cost-to-go starting from period 1. The first inequality holds because the fourth term, $\tau h_1u$, is the maximum difference in the holding costs at queue 1 between the first and the second processes over the course of period 0; by using $u$, which is the largest difference in the queue lengths by the start of period 1, we have established a lower bound on $\Delta_{\pi-\tilde{\pi}}$. The second inequality follows from $V^1(0,x^1_2) - V^1(0,y^1_2) \geq 0$, which holds by monotonicity of $V^m(\cdot)$ for all $m\in \mathcal{M}$ (Theorem \ref{thm:properties_V}). Note that $x^1_2 \geq y^1_2$ holds because under policy $\tilde{\pi}$, queue 2 in the second process receives $\hat{u} \leq u$ customers at the start of period 1. The last inequality follows from $u \geq \hat{u}$ and $h_2 \geq h_1$. Therefore, $\tilde{\pi}$ performs equally well or better. This shows that for a two-queue system, there is always an optimal policy which does not involve transferring to a more expensive queue when its state is already large enough to last a period without emptying.
    

     To extend this result to a general $N$-queue system, consider again two processes $x^m = (x_1^m, \ldots, x_N^m)$ and $y^m = (y_1^m, \ldots, y_N^m)$, $m \in \mathcal{M}$, which start from the same initial condition but follow policies $\pi$ and $\tilde{\pi}$, respectively. Suppose that $\pi$ involves transferring $\pi^0_{ij}=u_{ij} > 0$ customers from queue $i$ to $j$, $i\neq j$, at the start of period 0 when $h_j \geq h_i$ and $q^{\pi}_j[0] \geq \tau(\mu_j - \lambda_j)^+$. The policy $\tilde{\pi}$ is identical to $\pi$ except that in period 0, it does \emph{not} move customers from queue $i$ to $j$, i.e., $\tilde{\pi}^0_{ij} = 0$. Since the two processes are identical other than at queues $i$ and $j$, we can follow the same analysis above with the two-queue system (where we replace queue 2 with $j$ and queue 1 with $i$) and show that $\tilde{\pi}$ performs equally well or better than $\pi$. We can thus think of a sequence of policies $\{\tilde{\pi}_n\}$ where policy $\tilde{\pi}_n$ improves upon policy $\tilde{\pi}_{n-1}$ in the same manner until there are no pairs of queues $(k,l)$ under $\tilde{\pi}_n$ with $u_{kl} > 0$ in period 0 when $h_l \geq h_k$  and $q^{\pi}_l[0] \geq \tau(\mu_l - \lambda_l)^+$. This shows that there always exists an optimal policy which does not transfer customers to a more expensive queue when it already has enough customers to last a period without emptying.
\end{proof}

\subsubsection{Preliminaries for Proving Proposition \ref{prop:trade_off_no_transfer_setup_costs}.}
In this section, we formally define the concept of directional derivative, derive the closed-form expression for the derivative of the single-period holding cost function for a two-queue system, and lastly prove monotonicity of the value function in the \emph{total} number of customers when there are no transfer and setup costs.

We first define the concept of directional derivatives. Let $z$ be a feasible direction at $x$, i.e., $\sum_{i \in \mathcal{N}}z_i = 0$ and $x + z \geq 0$. The condition $\sum_{i \in \mathcal{N}}z_i = 0$ ensures that any new state along the feasible direction preserves the total number of customers. Define $W^m(x) = H^m(x) + V^{m+1}(f^m(x,\tau))$ for all $m \in \mathcal{M}$, and define the directional derivative of $W^m(x)$ at $x$ along the feasible direction $z$ as
\begin{equation*}
    \nabla_zW^m(x) \equiv \lim_{t \xrightarrow{} 0^+}\frac{W^m(x + tz) - W^m(x)}{t}.
\end{equation*}
If $\kappa(\cdot)=0$ (no setup costs), $W^m(\cdot)$ is convex and continuous and $\nabla_zW^m(\cdot)$ is well-defined, i.e., \jp{it always exists and is finite \citep[Theorem 23.1]{rockafellarConvexAnalysis2015}}.

\begin{lemma}\label{lem:directional_derivative}
    Suppose $\kappa(\cdot) = 0$. For all $m\in \mathcal{M}$, it is optimal not to transfer if and only if $\nabla_zW^m(x) \geq -R(z)$ for all feasible direction $z$ at a given initial condition $x$.
\end{lemma}
\begin{proof}{Proof.}
The proof approach is available in \cite{benjaafarDynamicInventoryRepositioning2022}, which we include here. Suppose that it is optimal not to transfer at $x$. Then, based on the optimality equation, we must have
\begin{equation*}
    W^m(x + tz) + tR(z) \geq W^m(x) \iff \frac{W^m(x + tz) - W^m(x)}{t} \geq -R(z)
\end{equation*} 
for all $t > 0$. Taking the limit as $t \xrightarrow{} 0^+$, we obtain $\nabla_zW^m(x) \geq -R(z)$.

Now, suppose that $\nabla_zW^m(x) \geq -R(z)$ holds for all feasible direction $z$ at $x$. Define $\omega_m(t) = W^m(x + tz)$. \jp{Then $\omega_m(t)$ is convex (by Lemma \ref{lem:properties_H} and Corollary \ref{cor:properties_V_no_setup})}, $\omega_m(0) = W^m(x)$, and $\nabla_zW^m(x)$ can be expressed as $\omega_m'(0^+)$. By the subgradient inequality, we have $\omega_m(t) \geq \omega_m(0) + t\omega_m'(0^+) \geq \omega_m(0) - tR(z)$. Thus, $\omega_m(t) - \omega_m(0) = W^m(x + tz) - W^m(x) \geq -tR(z)$. So, it is optimal not to transfer customers at $x$.
\end{proof}

Next, we explicitly characterize \jp{the holding cost function for a two-queue system} under the assumption of stationary arrival rates and linear holding costs \jp{(hence we remove the dependence of $H$ on the period $m$)}. For any $y\in\mathbb{R}_+^N$, the holding cost at queue $i$ is given by
\begin{equation*}
    H_i(y_i) = \int^\tau_0 h_i(y_i + \lambda_i s - \mu_i s)^+ds = 
    \begin{cases}
    \frac{h_i}{2(\mu_i - \lambda_i)}y_i^2, & \mbox{if } 0 \leq y_i \leq \tau(\mu_i - \lambda_i)^+, \\
    h_i[y_i\tau + \frac{1}{2}(\lambda_i - \mu_i)\tau^2], & \mbox{if } y_i \geq \tau(\mu_i - \lambda_i)^+,
    \end{cases}
\end{equation*}
for all $i \in \mathcal{N}$ and thus
\begin{equation}\label{eq:H_two_queue}
    H(y) = H_1(y_1) + H_2(y_2) = 
    \begin{cases}
    \frac{h_1}{2(\mu_1-\lambda_1)}y_1^2 + \frac{h_2}{2(\mu_2-\lambda_2)}y_2^2, & y \in \mathcal{A}_1, \\
    h_1\left[y_1\tau + \frac{1}{2}(\lambda_1 - \mu_1)\tau^2 \right] + \frac{h_2}{2(\mu_2-\lambda_2)}y_2^2, & y \in \mathcal{A}_2, \\
    \frac{h_1}{2(\mu_1-\lambda_1)}y_1^2 + h_2\left[y_2\tau + \frac{1}{2}(\lambda_2 - \mu_2)\tau^2 \right], & y \in \mathcal{A}_3, \\
    h_1\left[y_1\tau + \frac{1}{2}(\lambda_1 - \mu_1)\tau^2 \right] + h_2\left[y_2\tau + \frac{1}{2}(\lambda_2 - \mu_2)\tau^2 \right], & y \in \mathcal{A}_4,
    \end{cases}
\end{equation}
where
\begin{align*}
    \mathcal{A}_1 &\equiv \{y: y_1 \leq \tau(\mu_1-\lambda_1)^+ \text{ and } y_2 \leq \tau(\mu_2-\lambda_2)^+ \}, \\
    \mathcal{A}_2 &\equiv \{y: y_1 \geq \tau(\mu_1-\lambda_1)^+ \text{ and } y_2 \leq \tau(\mu_2-\lambda_2)^+ \}, \\
    \mathcal{A}_3 &\equiv \{y: y_1 \leq \tau(\mu_1-\lambda_1)^+ \text{ and } y_2 \geq \tau(\mu_2-\lambda_2)^+ \}, \\
    \mathcal{A}_4 &\equiv \{y: y_1 \geq \tau(\mu_1-\lambda_1)^+ \text{ and } y_2 \geq \tau(\mu_2-\lambda_2)^+ \}.
\end{align*}

Finally, we show that when there are no transfer and setup costs, the value function is in fact non-decreasing in the \emph{total} number of customers. This means that even when a state is not component-wise smaller than another, if the total number of customers at that state is smaller, its value function must be smaller. 
\begin{lemma}\label{lem:monotonicity_total}
    Suppose $\kappa(\cdot) = 0$ and $r_{ij}=0$ for all $i,j \in \mathcal{N}$. Then for all $m \in \mathcal{M}$, $V^m(x) \geq V^m(z)$ for any $x$ and $z$ such that $\mathrm{e}^\top x \geq \mathrm{e}^\top z$, i.e., $V^m(\cdot)$ is non-decreasing in the total number of customers.
\end{lemma}
\begin{proof}{Proof.}
    The proof is by induction. Assume period $M-1$. Denote by $y$ the target state corresponding to a given initial condition $x$. Consider another state $z$ such that $\mathrm{e}^\top x \geq \mathrm{e}^\top z$. Suppose we construct a state $\hat{z} \in \Delta(\mathrm{e}^\top z)$ in the following way. Prescribe all customers at queue $1$ until $\hat{z}_1 = y_1$. If this is impossible (because $\mathrm{e}^\top z < y_1$), stop; otherwise, continue on to queue 2 and prescribe all remaining customers ($\mathrm{e}^\top z - y_1$) until there are none left or until $\hat{z}_2 = y_2$. Proceeding in this way with queues 3, $\ldots, N$, we must have that $y \geq \hat{z}$. Then given no transfer and setup costs, we observe
    \begin{align*}
        V^{M-1}(x) &=  \jp{H^{M-1}}(y) \geq \jp{H^{M-1}}(\hat{z}) \geq V^{M-1}(z).
    \end{align*}
    This shows that $V^{M-1}(\cdot)$ is non-decreasing in the total number of customers. 
    
    Now, suppose the claim holds for period $m+1, \ldots, M-1$. In period $m$, for an arbitrary initial condition $x$, we have
    \begin{equation*}
        V^m(x) = \min_{y \in \Delta(\mathrm{e}^\top x)} \left[ \jp{H^m}(y) + V^{m+1}(f^m(y,\tau)) \right].
    \end{equation*}
    Consider again a state $z$ such that $\mathrm{e}^\top x \geq \mathrm{e}^\top z$. Denote by $y$ the target state corresponding to the initial condition $x$ and a state $\hat{z} \in \Delta(\mathrm{e}^\top z)$ which we construct in the same manner as before. Due to the monotonicity of $\jp{H^m}(\cdot)$, $V^{m+1}(\cdot)$, and $f^m(\cdot,\tau)$, it follows by the same argument that $V^m(x) \geq V^m(z)$.
\end{proof}

\subsubsection{Proof of Proposition \ref{prop:trade_off_no_transfer_setup_costs}.}
\begin{proof}{Proof of Proposition \ref{prop:trade_off_no_transfer_setup_costs}.} ($i$)
    We prove the bounds on the target state by showing that in any period, if $x_i > \tau(\mu_i-\lambda_i)^+$ for some $i \geq 2$, we can always find a feasible direction along which the total cost improves. The contrapositive of this statement states that if it is optimal not to transfer at a state (i.e., no feasible directions improve the total cost), it must satisfy $x_i \leq \tau(\mu_i-\lambda_i)^+$ for all $i \geq 2$. \jp{Given the stationarity of arrival rates, we remove the dependence of $H$ on the period $m$ below.}
    
    Consider any period $m \in \mathcal{M}$. Suppose a given initial condition $x$ satisfies $x_i > \tau(\mu_i - \lambda_i)^+$ for some $i \geq 2$. Consider a policy which transfers customers from queue $i$ to 1 while preserving the total number of customers between the two, in such a way that no other queues are affected. This policy can be represented by a feasible direction $z \in \mathbb{R}^N$ such that $z_i = -\delta$ and $z_1 = \delta$ for some $\delta > 0$ and $z_l=0$ for all other $l$. Define $W^m(x) = H(x) + V^{m+1}(f^m(x,\tau))$. Then, 
    \begin{align*}
        \nabla_zW^{m}(x) &= \lim_{t \xrightarrow{} 0^+} \frac{H(x+tz) - H(x) + V^{m+1}(f^m(x+tz,\tau)) - V^{m+1}(f^m(x,\tau))}{t} \\
        &\leq \lim_{t \xrightarrow{} 0^+} \frac{1}{t}[th_1 - th_i]\tau\delta + \nabla_zV^{m+1}(f^m(x,\tau))\\
        &\leq (h_1 - h_i)\tau\delta \leq 0.
    \end{align*}
    This shows that $z$ is a feasible direction that leads to an equally good or better state. \jp{The \jp{second} inequality holds since $R(\cdot)=0$ \jp{by the assumption in the statement of the proposition} and $\nabla_zV^{m+1}(f^m(x,\tau)) \leq 0$.} \jp{Indeed, the feasible direction $z$ ensures $f^m(x+tz,\tau)^\top \mathrm{e} \leq f^m(x,\tau)^\top \mathrm{e}$, which indicates that by the end of the period, the total number of customers starting from $x+tz$ is less than or equal to that starting from $x$; therefore, by Lemma \ref{lem:monotonicity_total}, we must have $V^{m+1}(f^m(x+tz,\tau)) - V^{m+1}(f^m(x,\tau)) \leq 0$ for small enough $t > 0$}. Dividing both sides by $t$ and taking the limit $t\rightarrow 0^+$, we have $\nabla_z V^{m+1}(f^m(x,\tau)) \leq 0$. \jp{Lastly, from the assumption that $h_1 \leq h_i$, we conclude $\nabla_zW^{m}(x) \leq 0$. Thus, invoking Lemma \ref{lem:directional_derivative}, it is optimal to transfer customers.} Repeating this argument for all $i \geq 2$, we conclude that a target state $y$ must satisfy $y_i \leq \tau(\mu_i - \lambda_i)^+$ for $i \geq 2$.

    ($ii$) We first want to show that it is optimal not to transfer \jp{only if} an initial condition $x$ satisfies $x \geq \tau(\mu - \lambda)^+$. The reverse direction follows directly from Proposition \ref{prop:trade_off}. To prove the forward direction, assume that it is not optimal to transfer at $x$. For sake of contradiction, suppose $x_i < \tau(\mu_i - \lambda_i)^+$ for some $i$. Consider any period $m \in \mathcal{M}$ and a policy which transfers $\delta > 0$ customers to queue $i$ from $j$, where $j\neq i$ and \jp{$x_j \geq \tau(\mu_i - \lambda_i)^+$}. \jp{Similar to the proof of part $(i)$, this policy is represented by a feasible direction $z \in \mathbb{R}^N$ where $z_i=\delta$ and $z_j=-\delta$.} Define $W^m(x) = H(x) + V^{m+1}(f^m(x,\tau))$. Then
    \begin{align*}
        \nabla_zW^{m}(x) &= \lim_{t \xrightarrow{} 0^+} \frac{H(x+tz) - H(x) + V^{m+1}(f^m(x+tz,\tau)) - V^{m+1}(f^m(x,\tau))}{t} \\
        &= \lim_{t \xrightarrow{} 0^+} \frac{1}{t} \jp{\left[\frac{h_i}{\mu_i - \lambda_i}tx_i\delta + \frac{h_i}{2(\mu_i - \lambda_i)}t^2\delta^2 \right] - h_j\tau\delta} + \nabla_zV^{m+1}(f^m(x,\tau)) \\
        &\leq \jp{h_i \left(\frac{x_i}{\mu_i - \lambda_i} - \tau \right)\delta} < 0.
    \end{align*}
    \jp{The second line follows from \eqref{eq:H_two_queue}, case $\mathcal{A}_3$: the first two terms in the brackets are equal to $H_i(x_i+t\delta) - H_i(x_i)$ while the third term, $-h_j\tau\delta$, is equivalent to $\lim_{t \xrightarrow{} 0^+}H_j(x_j-\delta) - H_j(x_j)$. By Lemma \ref{lem:directional_derivative}, the final inequality indicates \jp{($< 0$)} that it is \emph{strictly} optimal to transfer at $x$.} Contradiction. Therefore, if it is not optimal to transfer at $x$, it must satisfy $x_i \geq \tau(\mu_i - \lambda_i)^+$ for all $i$.

    Next, we show that \emph{any} $y \geq \tau(\mu - \lambda)^+$ is a target state. Consider an initial condition $x$ and suppose its corresponding target state is $y \geq \tau(\mu - \lambda)^+$. Consider another candidate target state $\hat{y}$ such that $\hat{y} \neq y$, $\mathrm{e}^\top \hat{y} = \mathrm{e}^\top y$, and $\hat{y} \geq \tau(\mu - \lambda)^+$. Since there are no transfer and setup costs, we note that the net transfer $\hat{y} - y$ can be represented as a series of vectors $z_1, \ldots, z_K$ in $\mathbb{R}^N$ involving two queues at a time such that $\sum_{k=1}^Kz_k = \hat{y} - y$. Following a similar analysis as above, we observe that for all $k=1,\ldots,K$,
    \begin{align*}
        \nabla_{z_k}W^{m}(y) = \lim_{t \xrightarrow{} 0^+} \frac{1}{t}(th_i - th_j)\tau\delta + \nabla_{z_k}V^{m+1}(f^m(y,\tau)) = (h_i - h_j)\tau\delta = 0.
    \end{align*}
    This indicates that we are indifferent to the choice of $y$ and $\hat{y}$ when moving the state from $x$, and as a result, both are optimal.
\end{proof}

\subsection{Proof of Proposition \ref{prop:connectedness} in Section \ref{ssec:preserve_connectedness} and Additional Discussion} \label{appen:proof_connectedness}
\begin{proof}{\jp{Proof.}}
    \jp{Since $p=0$ trivially guarantees connectedness, we state a definition of $p > 0$ below which leads to a connected region-of-inaction. For readability, we omit the dependence of the classifier $g_j$ on the period $m$.}

    \jp{Consider iteration $j$, a fixed $n \in \mathbb{Z}_+$, and define $\mathcal{X}(n) = \{x\in \mathcal{X}: \mathrm{e}^\top x = n\}$. Intuitively, $\mathcal{X}(n)$ is a ``slice'' of the state space corresponding to the simplex where the sum of the components is $n$ for all states. Consider two states $x_1, x_2 \in \mathcal{X}(n)$ such that they achieve the minimum and maximum probabilities according to the classifier $g_j$, i.e., $g_j(x_1) = \min_{x \in \mathcal{X}(n)} \{g_j(x)\}$ and $g_j(x_2) = \max_{x \in \mathcal{X}(n)} \{g_j(x)\}$. Define the path function $z_n(t) = x_1 + t(x_2-x_1)$, $t\in [0,1]$. Since $g_j$ is continuous and piecewise-monotone with finitely many pieces (Assumption \ref{ass:classifier}), we can partition $[0,1]$ into finitely many ordered sub-intervals $I_1^j, \ldots, I_R^j$ such that $g_j(z_n(t))$ is monotone on each $I_r^j$. Define 
    \begin{align*}
        p_{n,j} \equiv \max \left\{p: \exists t \in I_1^j \text{ with } g_j \left(z_n(t) \right)=p \text{, and } g_j \left(z_n(t)\right) \geq p \text{ for all } t \in I_2^j \cup \cdots \cup I_R^j \right\},
    \end{align*}
    with the convention that $p_{n,j} = \min_{t\in[0,1]} \{g_j(z_n(t))\}$ if no such $p$ exists. Let $p^* = \min_{n,j}\{p_{n,j}\}$. Since $g_j$ is bounded away from 0, it follows that $p^* > 0$. Additionally, define the label function as $\mathrm{label}_j(x,p)=1\{g_j(x) \geq p\}$. Since $g_j(z_n(t))$ is monotone on the first sub-interval $I_1^j$, and $g_j(z_n(t)) \geq p_{n,j}$ for all $t$ in subsequent sub-intervals, $\mathrm{label}_j(z_n(t),p^*)$ must be non-decreasing in $t$. Thus, the transition from label 0 to 1 happens at most once, and the set of labels with 1 is connected with the rest of the region-of-inaction.}
\end{proof}

\jp{\textbf{Discussion.} Based on the definition of $p_{n,j}$, relatively small thresholds (e.g., $p \leq 0.5$) should work well. We can also see how different properties of $g_j$ may impact the choice of the threshold $p^*$. For instance, if $g_j$ is Lipschitz continuous, $g_j(z_n(t))$ cannot drop very quickly as we move along the path $z_n(t)$ from $x_1$ to $x_2$, especially since $g_j(z(t))$ must eventually increase from $g_j(x_1)$ to $g(x_2)$. As such, Lipschitz continuity should allow for higher values of $p_{n,j}$ (and hence $p^*$). Additionally, if $g_j(z_n(t))$ is monotone on $[0,1]$, then there is only a single sub-interval $I_1^j=[0,1]$ for all $j$, and our algorithm guarantees connectedness of the region-of-inaction for \textit{any} choice of $p^*$.}

\section{Supplementary Material for the API Algorithm in Section \ref{sec:ADP}} \label{appen:supp_material_ADP}

\subsection{Pseudocode}\label{appen:pseudocode}
\jp{The main algorithm is outlined in Algorithm \ref{alg:API} while policy evaluation is outlined in Algorithm \ref{alg:policy_evaluation}. In Algorithm \ref{alg:policy_evaluation}, Steps \ref{step:follow_policy}--\ref{step:update} describe the ``forward pass'' in which costs and states (and their labels) are sampled one period at a time from the beginning of the horizon. Step \ref{step:backward_pass} collects sampled costs and labels, and is often easier to implement as a ``backward pass,'' i.e., starting with $m=M-1$ and proceeding to $m=0$.}
\begin{algorithm}[t]
\caption{\jp{Sampling-Based Approximate Policy Iteration}}
\label{alg:API}
{\color{black}
\begin{algorithmic}[1]
\STATE Input initial value function $\bar{V}^m_0$, initial classifier $g^m_0$, and number of iterations $j_{\mathrm{max}}$. Set $j=0$ and $u^m(x) = 0$ for all $x\in \mathcal{X}$ and $m\in \mathcal{M}$.
\STATE \textbf{Policy Evaluation:} Call Algorithm \ref{alg:policy_evaluation}. Denote the outputs as $\{l^m(x), \Bar{v}^m(x)\}_{m,x}$, $\{u^m(x)\}_{m,x}$, and $\{\mathcal{X}_j^{m, \mathrm{visited}}\}_m$.  \label{PE}
\STATE \textbf{Policy Improvement:} For each $m \in \mathcal{M}$ and $x \in \mathcal{X}_j^{m,\mathrm{visited}}$, update the value function by $\Bar{V}_{j+1}^{m}(x) = \frac{u^m(x)}{u^m(x)+1} \Bar{V}_{j}^{m}(x) + \frac{1}{u^m(x)+1}\Bar{v}^m(x)$ and the labels as $\mathrm{label}_{j+1}^m(x) = l^m(x)$. Using $\{\mathrm{label}_{j+1}^m(x): x\in \mathcal{X}\}$ as the target variable, train a new classifier $g_{j+1}^m$ for each $m$.
\STATE Set $j = j + 1$. If $j < j_{\mathrm{max}}$, go to Step \ref{PE}.
\RETURN $\{\bar{V}^m_j\}_{m,j}$ and $\{g^m_j\}_{m,j}$.
\end{algorithmic}}
\end{algorithm}

\begin{algorithm}[t]
\caption{\jp{Policy Evaluation}}
\label{alg:policy_evaluation}
{\color{black}
\begin{algorithmic}[1]
\STATE Input value functions $\{\bar{V}^m\}_m$, classifiers $\{g^m\}_m$, probability threshold $p$, number of simulation runs $B$, and number of cumulative updates $u^m(x)$.
\STATE Set $\mathcal{X}^{m, \mathrm{visited}}=\varnothing$ and $\bar{v}^m(x)=0$, $c^m(x) = 0$ for $x\in \mathcal{X}, m\in \mathcal{M}$.
\FOR{each point $x \in \mathcal{X}$}
    \FOR{$b=1$ to $B$}
        \STATE Set $m=0$ and $x^m = x$.
        \STATE Follow policy $\pi^m$ which is constructed as in \eqref{eq:policy}. Denote the target state for $x^m$ as $y^m$. \label{step:follow_policy}
        \STATE Simulate one period to sample number of arrivals $a^m$ and departures $d^m$. Denote the single-period holding and transfer costs by $H^m(y^m)$ and $C(y^m-x^m)$.
        \STATE Set $x^{m+1} = y^m + a^m + d^m$, $\mathcal{X}^{m,\mathrm{visited}} = \mathcal{X}^{m,\mathrm{visited}} \cup \{x^m, y^m\}$, $c^m(x^m) = c^m(x^m) + 1$, $c^m(y^m) = c^m(y^m) + 1$, and $m = m+1$. If $m < M$, go to Step \ref{step:follow_policy}. \label{step:update}
        \STATE Set $\bar{v}^m(x^m) = \bar{v}^m(x^m) + \sum_{k=m}^{M-1} (H^k(x^k) + C(y^k - x^k))$ and $l^m(x^m) = 1\{x^m = y^m\}$ for $m\in \mathcal{M}$.\label{step:backward_pass}
    \ENDFOR
\ENDFOR
\STATE For $m\in \mathcal{M}$ and $x \in \mathcal{X}^{m, \mathrm{visited}}$, set $\bar{v}^m(x) = \bar{v}^m(x) / c^m(x)$ and $u^m(x) = u^m(x) + c^m(x)$.
\RETURN $\{l^m(x), \bar{v}^m(x)\}_{m,x}, \{u^m(x)\}_{m,x}$, and $\{\mathcal{X}^{m, \mathrm{visited}}\}_m$.
\end{algorithmic}}
\end{algorithm}

\subsection{Impact of Initialization} \label{appen:initialization}
\jp{Figure \ref{fig:impact_initialization} illustrates the convergence of the API algorithm when initialized with the fluid policy (top row), compared to a naive policy (bottom row) that assumes the entire state space is the region-of-inaction. In general, while both approaches appear to converge to the same policy, it is faster under the fluid initialization.}
\begin{figure}[h]
    \FIGURE
    {\includegraphics[width=01\textwidth]{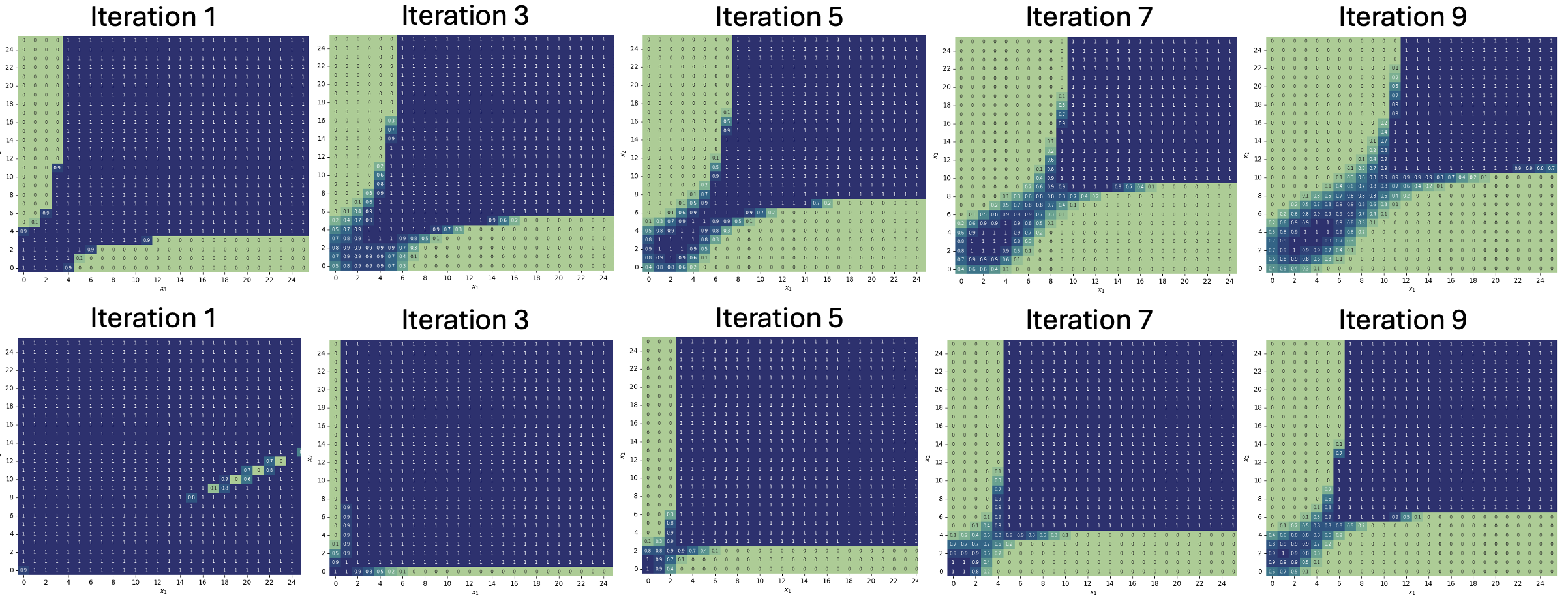}}
    {\jp{Estimated no-transfer region in the first nine iterations of API algorithm for a two-queue system using the fluid initialization (top row) versus naive initialization (bottom row)}\label{fig:impact_initialization}}
    {\jp{$\lambda_1=\lambda_2=7, \mu_1=\mu_2=10, h_1=h_2=5, r_{12}=r_{21}=1, K_{12}=K_{21}=1, \tau=1, M=7$.}}
\end{figure}

\subsection{Value of Common Random Numbers and Coupling} \label{appen:value_crns}
\jp{As described in Section \ref{sec:ADP}, Common random numbers (CRN) allow us to reduce the variance in the estimated value function differences while speeding up policy evaluation via coupling. Figure \ref{fig:impact_crns} illustrates the convergence of the API algorithm with (bottom row) and without CRNs (top row). Given the symmetric parameters, we note that the optimal policy must also be symmetric, i.e., $\pi(x_1,x_2)=\pi(x_2,x_1)$. The figure shows that not using CRNs can eventually result in an asymmetric, thus incorrect, policy.}
\begin{figure}[h]
    \FIGURE
    {\includegraphics[width=1\textwidth]{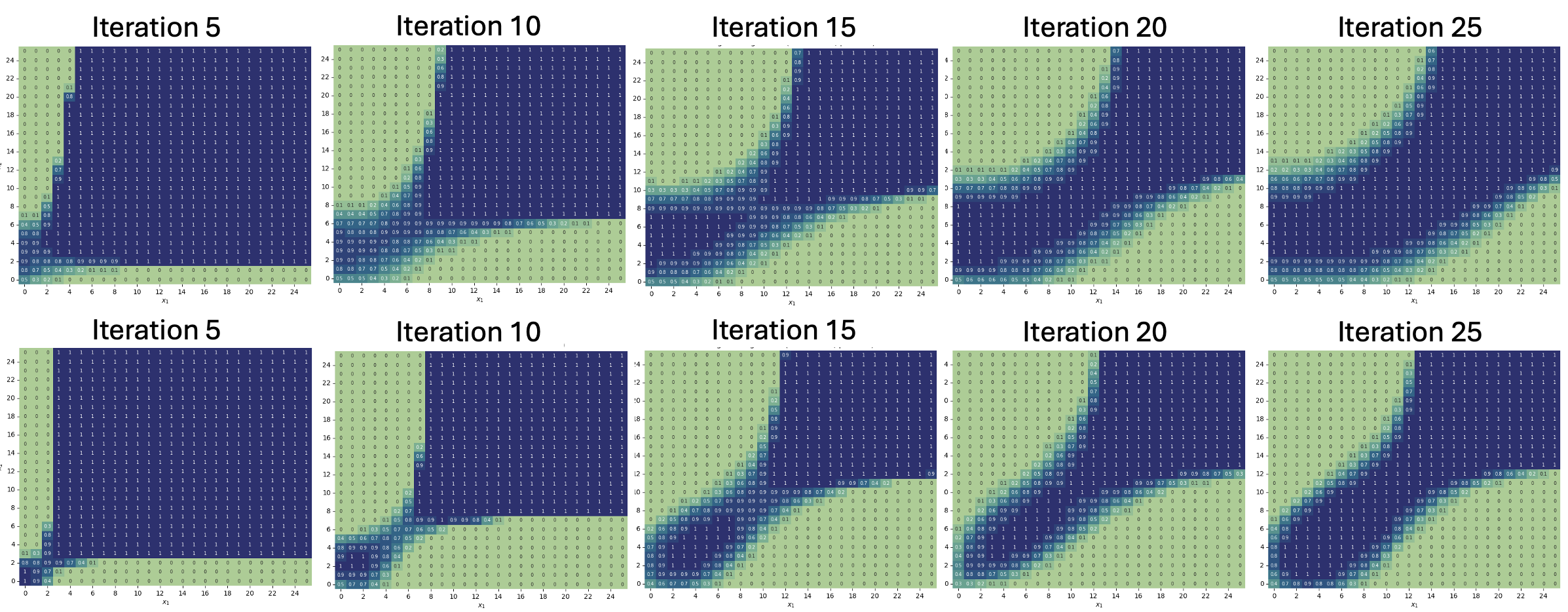}}
    {\jp{Estimated no-transfer region by the API algorithm for a two-queue system with (bottom row) and without CRNs (top row)}\label{fig:impact_crns}}
    {\jp{$\lambda_1=\lambda_2=9, \mu_1=\mu_2=10, h_1=h_2=5, r_{12}=r_{21}=1, K_{12}=K_{21}=1, \tau=1, M=7$.}}
\end{figure}


    


\subsection{Feature Importance} \label{appen:features}
\jp{We illustrate the strength of the distance features described in Section \ref{sec:ADP} in Figure \ref{fig:feature_comparison}, which shows the accuracy of logistic regression in approximating the no-transfer region for an example two-queue system, after running one iteration of Algorithm \ref{alg:API} using a probability threshold of $p=0.1$ and $B=10$ simulation runs per state. The system parameters are specified in Figure \ref{fig:feature_comparison}. The classifier is trained using only the polynomial features (Figure \ref{fig:polynomial_features}), only the distance features (Figure \ref{fig:distance_features}), or using both sets of features (Figure \ref{fig:all_features}). The target variable corresponds to the state labels from the policy evaluation step. Note that a high classification accuracy is attained only when using both sets of features (100\% accuracy).}
\begin{figure}[h]
    \FIGURE
    {\begin{subfigure}{0.33\textwidth}
        \FIGURE
        {\includegraphics[width=\textwidth]{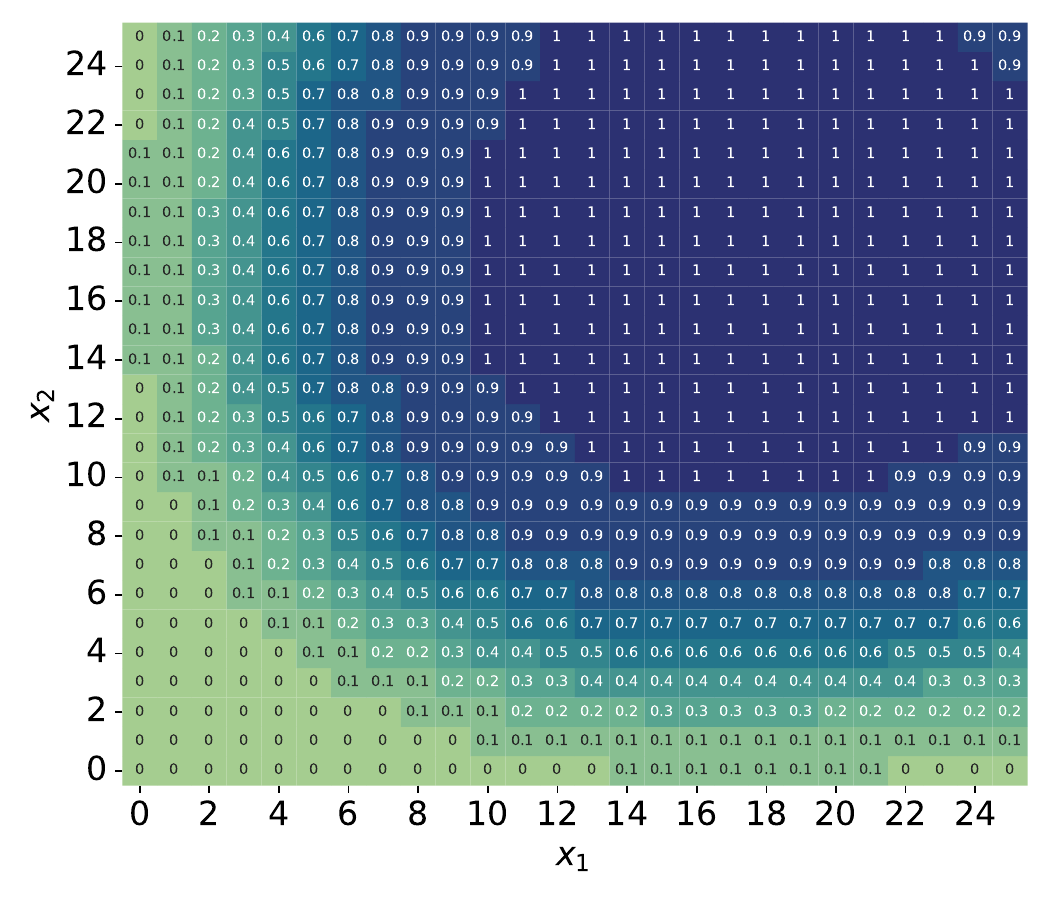}}
        {\jp{Polynomial features: 84\%}\label{fig:polynomial_features}}{}
    \end{subfigure}
    \hfill
    \begin{subfigure}{0.33\textwidth}
        \FIGURE
        {\includegraphics[width=\textwidth]{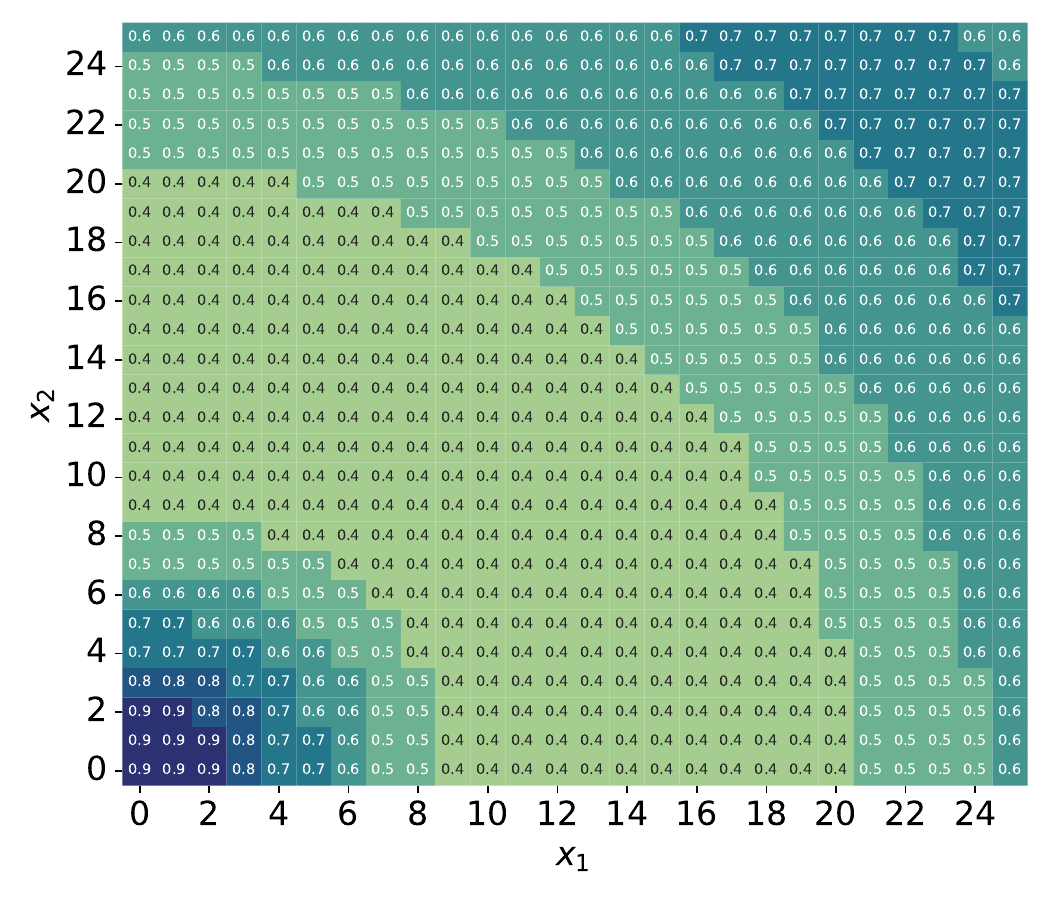}}
        {\jp{Distance features: 63\%} \label{fig:distance_features}}{}
    \end{subfigure}
    \hfill
    \begin{subfigure}{0.33\textwidth}
        \FIGURE
        {\includegraphics[width=\textwidth]{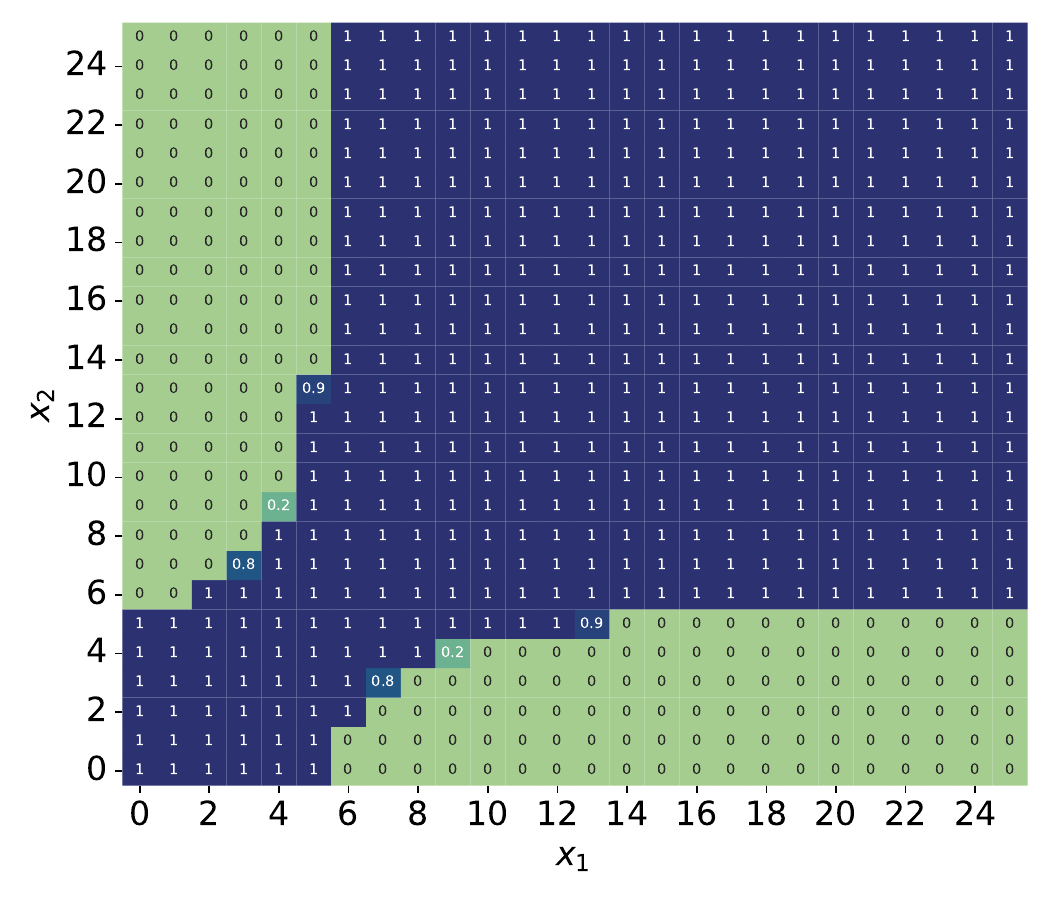}}
        {\jp{All features: 100\%} \label{fig:all_features}}{}
    \end{subfigure}}
    {\jp{Accuracy of logistic regression in estimating the no-transfer region for a two-queue system when the classifier is trained using different feature sets. State probabilities are plotted.} \label{fig:feature_comparison}}
    {\jp{$\lambda_1=\lambda_2=5, \mu_1=\mu_2=6, h_1=h_2=6, r_{12}=r_{21}=1, K_{12}=K_{21}=1, \tau=1, M=7$, $p=0.1$.}}
\end{figure}

\subsection{Impact of the Choice of Probability Threshold} \label{appen:choice_prob_threshold}
\jp{Certain choices of the probability threshold lead to a disconnected region-of-inaction because the classifier $g^m$ is generally not monotone. For example, consider $x$ and $x'$ such that $\mathrm{e}^\top x = \mathrm{e}^\top x'$ and $g^m(x) = g^m(x') = p_1$, and suppose there exists $\lambda \in (0,1)$ such that $g^m(\lambda x + (1-\lambda)x') = p_0 < p_1$. Then any $p \in (p_0, p_1]$ results in a disconnected region-of-inaction from $x$ to $x'$.}

\jp{In Figure \ref{fig:disconnected_roi}, we illustrate this phenomenon for a two-queue system after 5 iterations of the API algorithm. We use a probability threshold of $p=0.9$ throughout the algorithm. To obtain the final deterministic policy, we again apply $p=0.9$ to the state probabilities shown in Figure \ref{fig:disconnected_roi_probs} (probability of belonging to the region-of-inaction). The resulting no-transfer region in Figure \ref{fig:disconnected_roi_labels} consists of two distinct sub-regions, demonstrating that not all choices of the probability threshold will guarantee the connected structure.}

\begin{figure}[h]
    \FIGURE
    {\begin{subfigure}{0.4\textwidth}
        \FIGURE
        {\includegraphics[width=\textwidth]{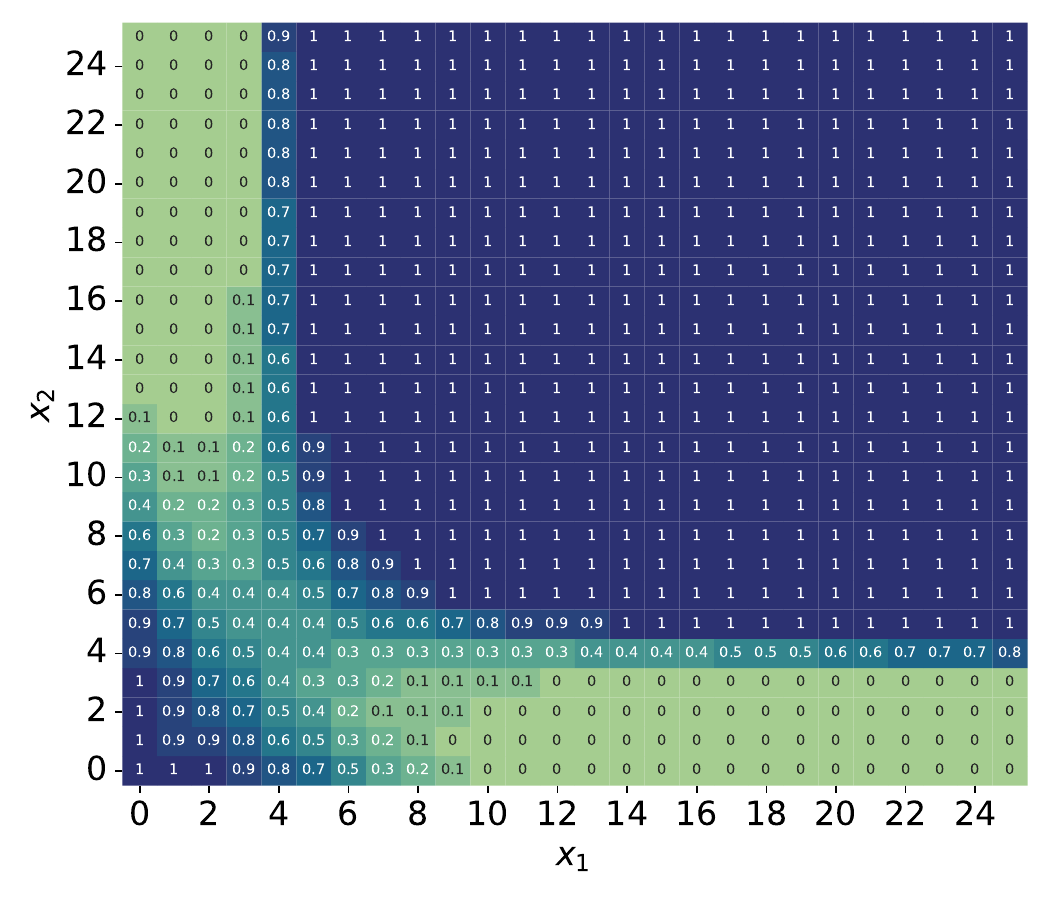}}
        {\jp{State probabilities $g^0(x)$} \label{fig:disconnected_roi_probs}}{}
    \end{subfigure}
    \hspace{2em}
    \begin{subfigure}{0.4\textwidth}
        \FIGURE
        {\includegraphics[width=\textwidth]{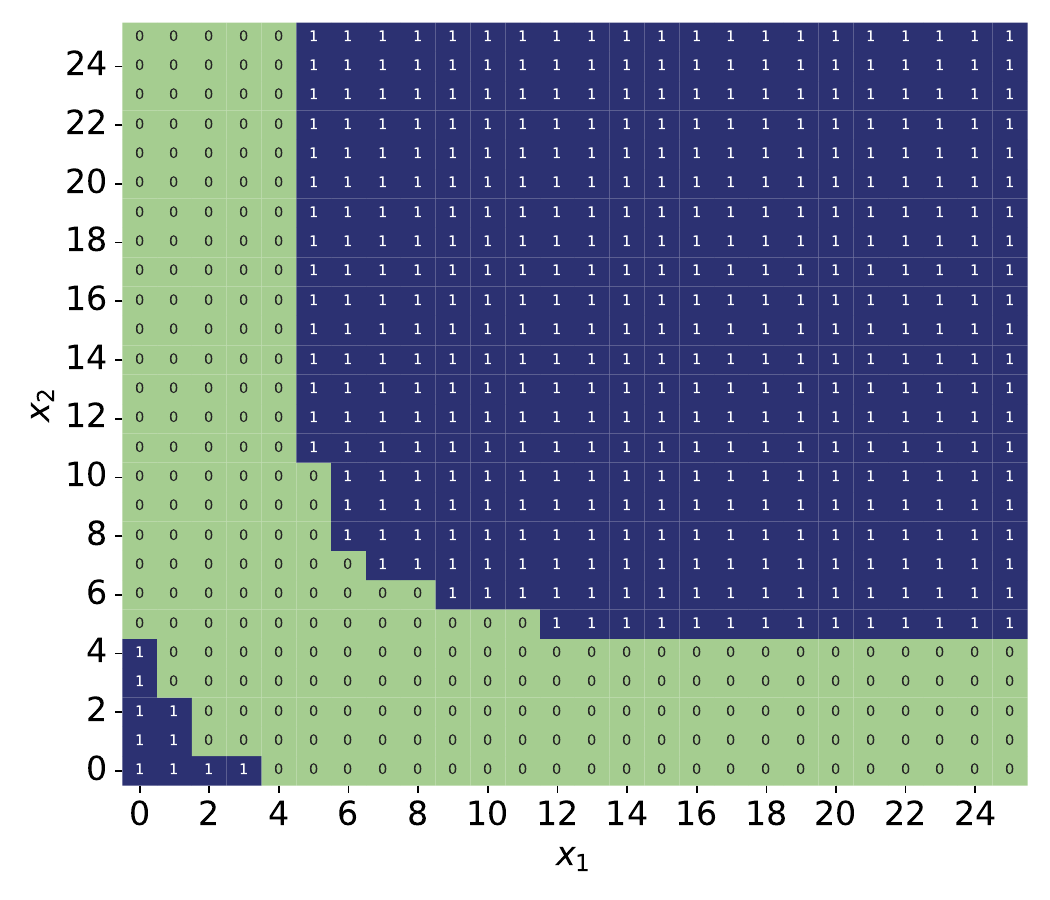}}
        {\jp{State labels $1\{g^0(x) \geq 0.9\}$} \label{fig:disconnected_roi_labels}}{}
    \end{subfigure}}
    {\jp{Estimated no-transfer region after 5 iterations of API algorithm for a two-queue system using a probability threshold of $p=0.9$}\label{fig:disconnected_roi}}
    {\jp{$\lambda_1=\lambda_2=5, \mu_1=\mu_2=6, h_1=h_2=6, r_{12}=r_{21}=1, K_{12}=K_{21}=1, \tau=1, M=7$.}}
\end{figure}

\subsection{A Check for Connectedness} \label{appen:check_for_connectedness}
\jp{We formulate the task of verifying connectedness of the region-of-inaction as a feasibility problem. For a given probability threshold $p$, let $\partial \tilde{\Sigma}^{\mathrm{outer}} \equiv \{x \in \partial \tilde{\Sigma}: g(x) < p\}$, where we have omitted the dependence on the period and iteration for ease of exposition. Denote by $y(x)$ the optimal solution to \eqref{eq:solve_target} for some $x \in \partial \tilde{\Sigma}^{\mathrm{outer}}$. Since $g(x) < p$, we note that $y(x) \neq x$. Partition $[0,1]$ into $D$ intervals of equal width, and define the path function $z_d(x)= x + (d/D)(y(x) - x)$ for $d=0,\ldots,D$. The level of discretization $D$ should be large enough to allow $z_d(x)$ to capture all integer-valued states in going from $x$ to $y(x)$. Then the following feasibility problem checks for connectedness:} 
\jp{\begin{equation} \label{eq:check_connected}
    \begin{split}
        \min &\quad 0 \\
        \mathrm{s.t.} &\quad 1\{g(z_d(x)) \geq p\} \leq 1\{g(z_{d+1}(x)) \geq p\}, \quad \forall d \in \{0,\ldots,D-1\}, x \in \partial \tilde{\Sigma}^{\mathrm{outer}}.
    \end{split}
\end{equation}
For example, using a logistic regression classifier, $g(z_d(x)) \geq p$ is simplified to $\beta^\top f(z_d(x)) \geq \theta$, where $\theta=\log(\frac{p}{1-p})$. The inequality in \eqref{eq:check_connected} forces monotonicity in the sequence of boundary labels in going from $x$ to $y(x)$. Thus, if an optimal solution is found, then the set of boundary labels indeed forms a connected set with the rest of the region-of-inaction. Additionally, since the expressions within $1\{\cdot\}$ are deterministic parameters, \eqref{eq:check_connected} can be reformulated as a mixed-integer linear program (MILP) and solved using an off-the-shelf solver. This check for connectedness can be incorporated at the start of the policy improvement step just prior to updating the value function.}

\section{Additional Details and Computations for Section \ref{sec:stochastic_system}} 
\subsection{MDP Solutions} \label{appen:MDPsolutions}
We provide details of the two-queue MDP formulation and solution method which are used to produce the figures in Section \ref{sec:stochastic_system}. 

The original state space of the two-queue system is described by $\mathcal{S} = \{(n,i) \in \mathbb{Z}_+ \times \mathbb{Z}_+:i \leq n \}$, where $n$ is the total number of customers in the system and $i$ is the number of customers at queue 2. This is truncated such that the maximum number in the system is at most $\Bar{n}=40$. To implement this, we impose the arrival rates to be $\lambda_1=\lambda_2=0$ when $n=40$. The action space is described by $\mathcal{A} = \{a \in \mathbb{Z}: -i \leq a \leq n-i\}$, which represents the number of customers to transfer from queue 1 to 2. When $a$ is negative, it signifies the opposite direction. Let $\Lambda = \lambda_1+\lambda_2+\mu_1+\mu_2$. We solve the following equations.
\begin{align*}
    V(n,i)= \min_{a\in\{-i,\ldots,n-i\}} \left[  K1\{a\neq0\} + r|a| + \frac{1}{\Lambda}\left\{ (n-i-a)h_1 + (i+a)h_2 + W(n,i+a) \right\} \right],
\end{align*}
where
\begin{align*}
    W(n,i)= \lambda_1 V(n+1,i)+\lambda_2 V(n+1,i+1) + \mu_1 V((n-1)^+,i) + \mu_2 V((n-1)^+,(i-1)^+),
\end{align*} 
and $W(0,0) = 0$ to enforce a terminal cost of zero when the system reaches an empty (absorbing) state.

\subsection{Additional Numerical Experiments for Section \ref{sec:stochastic_system}} \label{appen:add_results_stochastic_system}
\jp{Although the structure of the MDP policy matches that of the fluid policy, the exact parameter values can vary. For example, in Figure \ref{fig:compare_policy_unequal_light_traffic}, at $n=13$, the MDP policy suggests $s_2(13)=6$ and $S_2(13)=8$ while the fluid policy suggests $s_2(13)=8$ and $S_2(13)=9$. In this section, we examine the sub-optimality of the fluid policy by evaluating the performance of the fluid policy from a sequence of control problems indexed by $\eta$, as detailed in Section \ref{sec:DP}, for $\eta=1,2,3,4,5$. Because the fluid model becomes more accurate as the arrival and service rates increase, we expect the optimality gaps to decrease with $\eta$.}

We consider four two-queue systems in Table \ref{tab:convergence_gap} with parameters $\mu=(1,1)$, $h=(1,1)$, $r_{12}=r_{21}=2$, and $K_{12}=K_{21}=5$. We vary the traffic intensity of the system \jp{$\rho \equiv (\lambda_1 + \lambda_2)/(\mu_1 + \mu_2)$} between 0.6 and 0.8 in the experiments. For each $\rho$, the systems further differ by the arrival rates between the two queues. We set $\lambda_2-\lambda_1=0.2$ in one case and $\lambda_2-\lambda_1=0.1$ in the other. For each system, we use a common set of 20 randomly sampled initial conditions from $\mathcal{I} = \{x^0\in \mathbb{R}^2_+: 10 \leq x^0_1 + x^0_2 \leq 20\}$. Using 1,000 sample paths starting from each sampled initial condition, we compute the optimality gap, defined as the mean expected relative difference between the system costs under the fluid policy and the MDP policy until the first time its state reaches $(0, 0)$. The optimality gaps are then averaged across all initial conditions. 

\jp{Table \ref{tab:convergence_gap} shows that at $\eta=5$, the mean optimality gap is small in all cases. Moreover, the maximum optimality gap is always less than 5\% for $\eta=5$. To ensure that the small optimality gaps are not due to cases with negligible transfers, we also demonstrate that the fluid policy performs significantly better than the no-transfer policy. Additionally, even for small systems under $\eta=1$ or $\eta=2$, the mean gaps of the fluid policy are near or less than 5\% while still performing significantly better than the no-transfer policy.}
\begin{table}[h]
\renewcommand*{\arraystretch}{1.15}
\TABLE
{Fluid policy's optimality gap to MDP policy (optimal policy) under increasing scaling parameter $\eta$ \label{tab:convergence_gap}}
{
{\color{black}\begin{tabular}{ccccc}
\multicolumn{5}{c}{\up Case (a): $\rho=0.8$, $\lambda_1=0.7$, $\lambda_2=0.9$} \\
\hline
$\eta$ & Mean & Min. & Max. & Mean gap to no-transfer \\
\hline
1 & 6.5\% & 5.0\% & 9.2\% & -30.5\% \\
2 & 4.1\% & 1.8\% & 7.5\% & -33.7\% \\
3 & 3.0\% & 1.4\% & 5.1\% & -33.7\% \\
4 & 2.5\% & 1.0\% & 4.4\% & -33.6\% \\
5 & 1.8\% & 0.6\% & 3.5\% & -32.0\% \\
\hline
\multicolumn{5}{c}{\up Case (c): $\rho=0.6$, $\lambda_1=0.5$, $\lambda_2=0.7$} \\
\hline
$\eta$ & Mean & Min. & Max. & Mean gap to no-transfer \\
\hline
1 & 0.2\% & -0.4\% & 1.7\% & -17.1\% \\
2 & 0.4\% & 0.1\% & 1.2\% & -19.7\% \\
3 & 0.1\% & -0.6\% & 0.8\% & -20.1\% \\
4 & 0.0\% & -0.3\% & 0.4\% & -20.2\% \\
5 & 0.1\% & -0.3\% & 0.1\% & -20.3\% \\
\hline
\end{tabular}
\quad
\begin{tabular}{ccccc}
\multicolumn{5}{c}{\up Case (b): $\rho=0.8$, $\lambda_1=0.75$, $\lambda_2=0.85$} \\
\hline
$\eta$ & Mean & Min. & Max. & Mean gap to no-transfer \\
\hline
1 & 2.0\% & 0.8\% & 3.4\% & -28.9\% \\
2 & 2.3\% & 1.0\% & 4.3\% & -29.8\% \\
3 & 1.7\% & 0.7\% & 3.1\% & -28.8\% \\
4 & 1.3\% & 0.6\% & 2.3\% & -27.6\% \\
5 & 1.2\% & 0.3\% & 2.2\% & -26.6\% \\
\hline
\multicolumn{5}{c}{\up Case (d): $\rho=0.6$, $\lambda_1=0.55$, $\lambda_2=0.65$} \\
\hline
$\eta$ & Mean & Min. & Max. & Mean gap to no-transfer \\
\hline
1 & 0.1\% & -0.4\% & 0.7\% & -16.6\% \\
2 & 0.5\% & 0.0\% & 1.0\% & -19.6\% \\
3 & 0.4\% & 0.1\% & 0.8\% & -19.4\% \\
4 & 0.3\% & 0.1\% & 0.5\% & -19.7\% \\
5 & 0.3\% & 0.1\% & 0.5\% & -19.8\% \\
\hline
\end{tabular}}
}
{\emph{Note.} The last column represents the relative difference to the no-transfer policy, where negative numbers indicate improvement (reduction) in system cost.}
\end{table}

\subsection{Additional Results for Section \ref{sec:adp_two_queues}} \label{appen:add_results_adp_two_queues}
\jp{Tables \ref{tab:mg1_policy_performance} and \ref{tab:mtg1_policy_performance} present the results of the simulation experiments for the two-queue $M/G/1$ and $M(t)/G/1$ systems, respectively.}

\begin{table}[h]
\renewcommand*{\arraystretch}{1.15}
\TABLE
{Performance of Myopic, fluid, and API policies relative to no-transfer for $M/G/1$ system \label{tab:mg1_policy_performance}}
{
{\color{black}\begin{tabular}{clrrl}
\hline
Initial condition & Policy & Holding cost & Transfer cost & Reduction (\%) \\
\hline
\multirow{3}{*}{(1, 15)}  & Myopic & 1057.2 & 4.6 & 4.0 $\pm$ 0.5\% \\
                          & Fluid  & 975.4  & 10.3 & 9.9 $\pm$ 1.1\% \\
                          & API    & 870.2  & 26.3 & \textbf{16.0} $\pm$ \textbf{2.4}\% \\
\hline
\multirow{3}{*}{(1, 17)}  & Myopic & 1149.5 & 4.6 & 3.8 $\pm$ 0.5\% \\
                          & Fluid  & 1058.6 & 10.6 & 10.0 $\pm$ 1.2\% \\
                          & API    & 933.5  & 27.5 & \textbf{17.8} $\pm$ \textbf{2.1}\% \\
\hline
\multirow{3}{*}{(1, 19)}  & Myopic & 1248.2 & 4.6 & 3.7 $\pm$ 0.4\% \\
                          & Fluid  & 1145.8 & 11.0 & 10.3 $\pm$ 1.1\% \\
                          & API    & 1004.8 & 27.8 & \textbf{19.4} $\pm$ \textbf{1.8}\% \\
\hline
\multirow{3}{*}{(1, 21)}  & Myopic & 1350.5 & 4.6 & 3.5 $\pm$ 0.4\% \\
                          & Fluid  & 1235.1 & 11.5 & 10.8 $\pm$ 1.1\% \\
                          & API    & 1079.1 & 30.0 & \textbf{20.1} $\pm$ \textbf{1.8}\% \\
\hline
\multirow{3}{*}{(1, 23)}  & Myopic & 1454.1 & 4.6 & 3.5 $\pm$ 0.4\% \\
                          & Fluid  & 1333.7 & 11.6 & 10.7 $\pm$ 1.0\% \\
                          & API    & 1152.7 & 30.3 & \textbf{21.2} $\pm$ \textbf{1.7}\% \\
\hline
\end{tabular}}
}
{\emph{Note.} $\lambda=(0.9,0.9), \mu=(1,1), \tau=1, M=7, h=(10,10), r_{12}=r_{21}=1, K_{12}=K_{21}=1$.}
\end{table}

\begin{table}[h]
\renewcommand*{\arraystretch}{1.15}
\TABLE
{Performance of Myopic, fluid, and API policies relative to no-transfer for $M(t)/G/1$ system \label{tab:mtg1_policy_performance}}
{
{\color{black}\begin{tabular}{clrrl}
\hline
Initial condition & Policy & Holding cost & Transfer cost & Reduction (\%) \\
\hline
\multirow{3}{*}{(1, 15)} & Myopic & 597.7 & 5.4 & 5.1 $\pm$ 0.9\% \\
                         & Fluid  & 527.8 & 10.3 & 11.8 $\pm$ 1.4\% \\
                         & API    & 454.0 & 19.4 & \textbf{18.9} $\pm$ \textbf{2.5}\% \\
\hline
\multirow{3}{*}{(1, 17)} & Myopic & 669.8 & 5.5 & 5.0 $\pm$ 0.8\% \\
                         & Fluid  & 586.6 & 11.0 & 12.7 $\pm$ 1.3\% \\
                         & API    & 498.7 & 19.8 & \textbf{21.8} $\pm$ \textbf{2.2}\% \\
\hline
\multirow{3}{*}{(1, 19)} & Myopic & 747.2 & 5.7 & 5.1 $\pm$ 0.7\% \\
                         & Fluid  & 652.0 & 11.6 & 13.3 $\pm$ 1.3\% \\
                         & API    & 547.3 & 21.2 & \textbf{23.5} $\pm$ \textbf{2.2}\% \\
\hline
\multirow{3}{*}{(1, 21)} & Myopic & 832.6 & 5.8 & 5.0 $\pm$ 0.7\% \\
                         & Fluid  & 719.0 & 12.3 & 14.5 $\pm$ 1.2\% \\
                         & API    & 599.2 & 22.7 & \textbf{25.2} $\pm$ \textbf{2.1}\% \\
\hline
\multirow{3}{*}{(1, 23)} & Myopic & 921.7 & 5.9 & 5.2 $\pm$ 0.6\% \\
                         & Fluid  & 791.1 & 13.0 & 15.5 $\pm$ 1.2\% \\
                         & API    & 654.6 & 24.1 & \textbf{27.0} $\pm$ \textbf{1.9}\% \\
\hline
\end{tabular}}
}
{}
\end{table}

\end{APPENDICES}








\end{document}